\documentclass{siamltex}
\usepackage{amsmath}
\usepackage{amssymb}
\usepackage{dsfont}
\usepackage{color}
\usepackage{graphicx}
\usepackage{subfig}
\usepackage{multirow}
\usepackage{morefloats}
\usepackage{pgfplots}
\usetikzlibrary{plotmarks}
\pgfplotsset{compat=newest}
\usepackage{tikz}

\renewcommand{\S}{\mathbb{S}}

\newtheorem{remark}[theorem]{Remark}

\begin{document}
\title{Mean field models for interacting ellipsoidal particles}

\author{R. Borsche \footnotemark[1]\and A. Klar\footnotemark[1] \footnotemark[2]
       \and A. Meurer \footnotemark[1]\and O. Tse\footnotemark[1] }
\footnotetext[1]{Technische Universit\"at Kaiserslautern, Department of Mathematics, Erwin-Schr\"odinger-Stra{\ss}e, 67663 Kaiserslautern, Germany 
  (\{borsche,klar,meurer,tse\}@mathematik.uni-kl.de)}
\footnotetext[2]{Fraunhofer ITWM, Fraunhoferplatz 1, 67663 Kaiserslautern, Germany}

\maketitle

% % % % % % % % % % % % %
% % % Introduction  % % %
% % % % % % % % % % % % %

\begin{abstract}
 We consider a mean field hierarchy of models for large systems of interacting ellipsoids suspended in an incompressible  fluid. The models range from microscopic to macroscopic mean field models. The microscopic model is based on three ingredients. Starting from a Langevin type model for rigid body interactions, we use a Jefferys type term to model the influence of the fluid on the ellipsoids and a simplified interaction potential between the ellipsoids to model the interaction between the ellipsoids. A mean field equation and corresponding equations for the marginals of the distribution function are derived and a numerical comparison between the different levels of the model hierarchy is given. The results clearly justify the suitability of the proposed approximations for the example cases under consideration.
\end{abstract}

\section{Introduction}

Large systems of interacting ellipsoidal shaped particles have attracted the attention of researchers from many different fields. For example, such systems are used to describe polymers and liquid crystals in the chemical sciences \cite{doi1988theory, kuzuu1983constitutive, walter2010stochastic}. The movement of ellipsoidal particles suspended in a fluid is also used in process engineering to describe the physics inside a liquid-liquid extraction column \cite{atwi2013three, mortensen2008orientation, ouchene2013drag, tavakol2015dispersion, zhang2001ellipsoidal}. A recent application of such models may be found in paper production processes \cite{lindstrom2008modelling,lindstrom2007simulation,  lindstrom2008simulation, lindstrom2009numerical}.

Mathematically the movement of ellipsoidal particles can be described on a microscopic level by large systems of ordinary differential equation based on Newtonian laws of mechanics for translational and rotational motion of the ellipsoids. For ellipsoidal particles suspended in an incompressible viscous fluid, one can use the model of Jeffery (e.g. \cite{jeffery1922motion, junk2007new, walter2010stochastic}), which describes the forces exerted by the fluid on an ellipsoid. In our work, the inter-particle interaction forces between ellipsoids are described via pairwise potentials for the particles and a random force. For the interaction potentials, we use potentials common in the literature for polymers \cite{berardi1995generalized, berne1972gaussian, cleaver1996extension, everaers2003interaction, gay1981modification, perram1996ellipsoid}, where the form of the ellipsoids are modeled with the help of Gaussian type functions. This leads to a Langevin-type microscopic model similar to the ones described in \cite{coffey1996langevin,han2006brownian, sun2008langevin,walter2010stochastic}. For the numerical treatment of large systems of hard interacting ellipsoids, we refer, for example, to \cite{DTS2005}.

For a very large number of particles, macroscopic equation for density, mean velocity, and other statistical quantities are expected to be a more efficient approximation of these models.
In the present work we derive, via a mean field approximation, corresponding kinetic equation, 
which can be used in turn to derive  hydrodynamic and diffusive limit equations. 
This procedure has also been used for example in the case of self-organizing systems of particles or for the description of pedestrian or granular flows \cite{canizo2011well,carrillo2009double,carrillo2010self,etikyala2014particle}.

The paper is organized as follows:
In section \ref{sec:micro} the microscopic model is introduced and the ellipsoidal interaction potential is constructed. 
From this, we derive the mean field limit equation in Section \ref{sec:Limit} and use different moment closure procedures for the derivation of different hydrodynamic limit equations in Section \ref{sec:hydro}.
In Section \ref{sec:num ex}, we numerically compare the derived models with the microscopic model for several different examples and flow fields.

% % % % % % % % % % % % % % %
% % % Microscopic Model % % %
% % % % % % % % % % % % % % %
\section{The Microscopic Model}
\label{sec:micro}
We consider a microscopic Langevin-type model as in \cite{walter2010stochastic} to describe the motion of ellipsoidal particles suspended in an incompressible  fluid in two space dimensions. 
The interaction of the fluid with the ellipsoids is described by a Jefferys type term \cite{jeffery1922motion,junk2007new}.
Interactions of the ellipsoids with each other are described by a many-particle interaction potential similar to \cite{berne1972gaussian}.

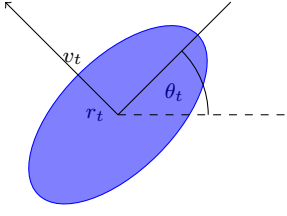
\begin{figure}
\centering
\begin{tikzpicture}[scale=3]
			\node (A) at (0.25,0.1){\footnotesize{$\theta_t$}};
			\node (B) at (-0.1,0){\footnotesize{$r_t$}};
			\node (C) at (-0.2,0.25){\footnotesize{$v_t$}};
			\draw[draw=blue,fill=blue,fill opacity=0.5] (0,0) ellipse [rotate=45, x radius=0.5, y radius=0.25];
			\draw[black] (0,0)--(0.5,0.5);
			\draw[dashed,black] (0,0)--(0.75,0);
			\draw[->,black] (0,0)--(-0.5,0.5);
			\draw (0.4,0) arc (0:45:0.4cm);
			\end{tikzpicture}
\caption{Sketch of an ellipsoidal particle with position $r_t$, orientation angle $\theta_t$ and velocity $v_t$.}
\label{fig:ellipse}
\end{figure}

As illustrated in Figure~\ref{fig:ellipse}, each ellipsoidal particle is described by its position $r_t\in\mathbb{R}^2$, velocity $v_t\in\mathbb{R}^2$, orientation 
angle $\theta_t\in[0,~2\pi)$ and angular velocity $\omega_t\in\mathbb{R}$. 
The angle $\theta_t$ is given by the relative angle between the horizontal axis and the main axis of the ellipsoidal particle such 
that  the angle $\theta_t=0$ corresponds to the orientation $(1,0)^\top$.
The equations of motion for $N$ particles $i=1,\ldots,N$ are
\begin{align}\label{eq:micro} 
\left\{\begin{aligned}
dr_t^i &= v_t^i dt\\
dv_t^i &= \gamma(u-v_t^i)dt-\frac{1}{m}\frac{1}{N}\sum_{i\neq j} \nabla_{r_i}U(r_t^i,r_t^j,\theta_t^i,\theta_t^j)dt\\
&\hspace*{10em}-\nabla_{r} V_1(r_t^i)dt - (A^2/2)v_t^idt+A dW_t^{A,i}\\
d\theta_t^i &= \omega_t^idt\\
d\omega_t^i &= \bar{\gamma}(g(\theta_t^i,u)-\omega_t^i)dt-\frac{1}{I_c}\frac{1}{N}\sum_{i\neq j} \nabla_{\theta_i}U(r_t^i,r_t^j,\theta_t^i,\theta_t^j)dt\\
&\hspace*{10em}-\nabla_{\theta} V_2(\theta_t^i) dt - (B^2/2)\omega_t^idt+Bd W_t^{B,i},
\end{aligned}\right.
\end{align}
with appropriate initial conditions.
Here $u$ is the velocity of a stationary surrounding fluid and $g(\theta,u)$ is given by
\begin{align*}
g(\theta,u)&=\frac{1}{2}\text{rot}(u)+\lambda\left(
\begin{array}{r}
-\sin\theta\\
\cos\theta
\end{array} \right)^\top\left(\frac{1}{2}(\nabla u+\nabla u^\top)\right)\left(
\begin{array}{c}
\cos\theta\\
\sin\theta
\end{array} \right).
\end{align*}
The first terms on the right hand side of  the velocity and angular velocity equations describe
the relaxation  of the particles to  the velocity of the fluid and to the rotation resulting from the velocity field, respectively.
The speed of relaxation is determined by the parameters $\gamma$ and $\bar{\gamma}$. 
The second terms model the repulsive interaction between the particles.
The parameters $m$ and $I_c$ are the mass and the moment of inertia of the particles.
The functions $V_1,V_2$ model an outer potential like for example gravitation or a magnetic field. 
The parameters $A,B$ are nonnegative diffusion constants and $W^{A,i},W^{B,i}$ are independent standard Brownian motions. 
The interaction potential is given by the following considerations.

There exist many different interaction potentials for ellipsoidal particles \cite{berardi1995generalized,berne1972gaussian,cleaver1996extension,everaers2003interaction,gay1981modification,perram1996ellipsoid}.
We use the soft potential as proposed by Berne \cite{berne1972gaussian}.
It is obtained by overlapping two ellipsoidal Gaussians representing the mutual repulsion of two particles. 
This leads to
\begin{align*}
\tilde{U}(r,\bar r,\theta,\bar \theta)&=a(\theta,\bar \theta)
\exp\left(-\left(\bar r-r\right)\left(\gamma (\theta)+\gamma(\bar \theta) \right)^{-1}\left(\bar r-r\right)\right),
\end{align*}
where $a$ and $\gamma$ are defined by
\begin{gather*}
a(\theta,\bar \theta)=\epsilon_0\left(1-\lambda^2
(
\eta(\theta)\cdot\eta(\bar\theta))^2\right)^{-\frac{1}{2}},
\qquad \eta(\theta)=(\cos\theta, \sin\theta)^\top,\\
\gamma (\theta) =\left(l^2-d^2\right)
\eta(\theta)\otimes
\eta(\bar\theta)
+d^2\mathds{1},\qquad
\lambda=\frac{l^2-d^2}{l^2+d^2}.
\end{gather*}
Here, $l=2L$ and $d=2D$ where $L$ is the length and the $D$ the width of the particle. 
The parameter $\epsilon_0$ models the strength of the potential.
To have compact support we slightly modify the potential and define
\begin{align}
U(r,\bar r, \theta,\bar \theta)&=a(\theta, \bar \theta)
\exp\left(-\frac{\left(\bar r-r\right)\left(\gamma(\theta)+\gamma(\bar \theta)\right)^{-1}\left(\bar r-r\right)}{1-\left(\bar r-r\right)\left(\gamma(\theta)+\gamma(\bar\theta)\right)^{-1}\left(\bar r-r\right)}\right).
\label{eq:interaction}
\end{align}
\begin{remark}
Note that although these interaction potentials prevent overlapping of particles, collisions are still possible. In order to completely avoid collisions, one has to consider hard core potentials or exclusions, which are nontrivial for ellipsoids.
\end{remark}
% % % % % % % % % % % % % % % % % % % % % % %
% % % Mean field and Macroscopic Limits % % %
% % % % % % % % % % % % % % % % % % % % % % %
\section{Mean field equations}
\label{sec:Limit}
For the number of ellipsoids $N$ tending to infinity, the microscopic system may be described in a probabilistic sense, where a partial differential equation is used to describe the evolution of the distribution of the particles in phase space. Using the {\em weak coupling} scaling limit 
\cite{braun1977vlasov,carrillo2009double,dobrushin1979vlasov,spohn1991Large}, i.e., the rescaling of the interaction 
potential with $1/N$ in the microscopic equations, one can  derive a kinetic mean field equation following for example \cite{carrillo2009double}. 

The empirical measure $f^N$ of the stochastic process $z_t^i=(r_t^i,\theta_t^i,v_t^i,\omega_t^i)\in \S$ on the state space $\S=\mathbb{R}^2\times [0,2\pi)\times \mathbb{R}^2\times \mathbb{R}$ is given by
%given by equation \eqref{eq:micro} 
$$
f^N(t,z)=\frac{1}{N}\sum_{i=1}^N \delta(z-z_t^i),
$$
where $\delta$ denotes the usual Dirac distribution and $z=(r,\theta,v,\omega)\in \S$.
The mean field limit describes  the convergence as $N \rightarrow \infty$ of the stochastic empirical measure $f^N$ towards the deterministic distribution $f$ of the stochastic process $z_t=(r_t,\theta_t,v_t,\omega_t)$ governed by the so-called nonlinear McKean--Vlasov equation
%where $z_t=(r_t,\theta_t,v_t,\omega_t)$ satisfies the 
\begin{align}
\left\{\begin{aligned}
dr_t &= v_tdt\\
dv_t &= \gamma(u-v_t)dt-\frac{1}{m}\int \nabla_{r}U(r_t,\bar{r},\theta_t,\bar{\theta}) f(t,\bar{z}) d\bar zdt\\
&\hspace*{10em} -\nabla_{r} V_1(r_t)dt - (A^2/2)v_tdt+AdW_t^{A}\\
d\theta_t &= \omega_t dt\\
d\omega_t &= \bar{\gamma}(g(\theta_t,u)-\omega_t)dt-\frac{1}{I_c}\int\nabla_{\theta}U(r_t,\bar{r},\theta_t,\bar{\theta})f(t,\bar{z}) d\bar zdt\\
&\hspace*{10em} -\nabla_{\theta} V_2(\theta_t) dt - (B^2/2)\omega_t dt+Bd W_t^{B}.
\end{aligned}\right.
\label{eq:ODE_mf} 
\end{align}
The corresponding differential equation for the evolution of the distribution on state space $f\colon\mathbb{R}^+ \times \S\rightarrow \mathbb{R}$, which is determined using It\^o's formula, is called the mean field equation. It is given by
\begin{align}\label{eq:meanfield}
 \partial_t f+v\cdot\nabla_r f+\omega\cdot \nabla_\theta f + \nabla_v\cdot(S_1[f]f)  + \nabla_\omega\cdot(S_2[f]f) = L_1[f] + L_2[f],
\end{align}
where the operators are given by
\begin{align*}
 S_1[f] &= \gamma(u-v) - \nabla_r V_1(r) - \frac{1}{m}\int \nabla_{r}U(r,\bar{r},\theta,\bar{\theta}) f(\bar z)\,d\bar z, \\
 S_2[f] &= \bar{\gamma}(g(\theta,u)-\omega) - \nabla_\theta V_2(\theta) - \frac{1}{I_c}\int \nabla_{\theta}U(r,\bar{r},\theta,\bar{\theta}) f(\bar z)\,d\bar z, \\
 L_1[f] &= (A^2/2)\nabla_v\cdot (vf+\nabla_v f), \\
 L_2[f] &= (B^2/2)\nabla_\omega\cdot (\omega f+\nabla_\omega f).
\end{align*}
 Due to conservation of mass, we normalize the initial condition to have that
 \[
 \int f(t,z)\, dz = 1 \quad \text{for all\; $t\ge 0$}.
 \]
 \begin{remark}
  Notice that the local equilibria in velocity and angular velocity are of Maxwellian type. Indeed, any distribution of the form $f_{\mathcal{M}}=q(r,\theta)\mathcal{M}(v,\omega)$, with
  \begin{align}\label{eq:maxwellian}
   \mathcal{M}(v,\omega) = Z^{-1}e^{-|v|^2/2-\omega^2/2}, %\qquad Z = \iint e^{-|v|^2/2-\omega^2/2}\,dv d\omega
  \end{align}
  where $Z$ is the normalizing factor for $\mathcal{M}$, satisfies $L_i[f_{\mathcal{M}}]=0$ for $i=1,2$.
 \end{remark}
 
In the next section we consider hydrodynamic limits of the kinetic mean field equation \eqref{eq:meanfield} and describe its derivations.

\section{Hydrodynamic Limit}
\label{sec:hydro}
In this section we formally derive closed equations for two marginals of the distribution function $f$ using different closure procedures. For further details on hydrodynamic limiting procedures in the case of self-propelled particles, we refer the reader to \cite{carrillo2009double, carrillo2010self, chuang2007state}.

In this setting, we define the following macroscopic quantities that describe the moments of marginals corresponding to the distribution $f$:
\begin{gather*}
q =\int f(t,r,\theta,v,\omega)d\Omega_1,\qquad
q\tilde{v} =\int v f(t,r,\theta,v,\omega)d\Omega_1, \\
q\tilde{\omega} = \int \omega f(t,r,\theta,v,\omega)d\Omega_1,
\end{gather*}
where $d\Omega_1=dvd\omega$. Obviously, we have conservation of probability, i.e.,
\[
\int q(t,r,\theta) drd\theta = 1\qquad\text{for all\; $t\ge 0$}.
\]
Integrating the mean field equation \eqref{eq:meanfield} against $(d\Omega_1, vd\Omega_1,\omega d\Omega_1)$, we obtain set of the balance equations with the continuity equation on $(r,\theta)\in\mathbb{R}^2\times\mathbb{R}$, given by
 \begin{align}
 \label{balance1}
	\partial_t q+\nabla_r\cdot(q\tilde{v})+\nabla_\theta\cdot(q\tilde{\omega})=0,
\end{align}
and the balance equations for $\tilde v$ and $\tilde\omega$
\begin{align}\label{balance2}
\begin{aligned}
\partial_t(q \tilde{v})+\nabla_r\cdot\left( \int v\otimes v f\, d\Omega_1\right) +\nabla_\theta \cdot\left(\int v\,\omega f\,d\Omega_1\right)
&= qS_1[q] -(A^2/2)q\tilde{v},\\
\partial_t(q\tilde{\omega})+\nabla_r\cdot\left(\int v\,\omega f\,d\Omega_1\right) + \nabla_\theta\cdot\left(\int \omega^2 f\,d\Omega_1\right)
&= qS_2[q] - (B^2/2)q\tilde{\omega},
\end{aligned}
\end{align}
where the operators are given by
\begin{align*}
 S_1[q] &= \gamma(u-\tilde{v})-\nabla_r V_1(r)-\frac{1}{m}\iint \nabla_{r}U(r,\bar{r},\theta,\bar{\theta})q(\bar{r},\bar{\theta})d\bar{r}d\bar{\theta},\\
 S_2[q] &= \bar{\gamma}(g(\theta,u)-\tilde{\omega}) - \nabla_\theta V_2(\theta) -\frac{1}{I_c}\iint \nabla_{\theta}U(r,\bar{r},\theta,\bar{\theta})q(\bar{r},\bar{\theta})d\bar{r}d\bar{\theta}.
\end{align*}

\subsection{$(r,\theta)$-marginals and the Maxwellian closure}
%\subsection{$(r,\theta)$-marginals and the Maxwellian closure} 
We now consider a moment closure based on the local equilibria of the velocities, which approximates the dependence of the distribution on $v$ and $\omega$ by
\[ 
 f\sim q\mathcal{M}(v-\tilde{v},\omega-\tilde{\omega}),
\]
where $\mathcal{M}$ is the Maxwellian given in \eqref{eq:maxwellian}. Inserting this Ansatz function into the integral terms above and dropping the tilde, we obtain the hydrodynamic equations
%in 2-D with zeroth and second moment equal to $1$ and first moment equal to $0$, which is the equilibrium distribution of the stochastic terms.
%We get
\begin{align}\label{eq:macro_q_max}
\left\{\begin{aligned}
\partial_t q+\nabla_r\cdot(qv)+\nabla_\theta\cdot(q\omega) &=0, \\
\partial_t(qv)+\nabla_r\cdot(q v\otimes v + q) +\nabla_\theta\cdot(qv\omega) &= qS_1[q] -(A^2/2)qv,\\
\partial_t(q \omega)+\nabla_r\cdot(qv\omega)+\nabla_\theta(q\omega^2 + q)  &= qS_2[q] -(B^2/2)q\omega.
\end{aligned}\right.
\end{align}
Notice the resemblance of the system with a dissipative isothermal Euler type model under the influence of an external and interaction potential, which suggests a trend towards equilibrium that we will numerically investigate in the next section (cf.~\cite{klartse14}).

\subsection{$(r,\theta)$-marginals and the mono-kinetic closure}
%\subsection{$(r,\theta)$-marginals and the mono-kinetic closure}
Starting again from the balance equations (\ref{balance1}) and (\ref{balance2}), we now use a moment closure with a mono-kinetic closure distribution in both velocity and angular velocity, i.e.,
\[ 
 f\sim q\delta(\tilde{v}-v)\otimes\delta(\tilde{\omega}-\omega),
\]
which suggests that the distributions in velocity and angular velocity are deterministic and evolves according to the moments $\tilde v$ and $\tilde\omega$. Using the mono-kinetic Ansatz in the balance equations and dropping the tilde leads to the system
\begin{align}\label{eq:macro_q_mono}
\left\{\begin{aligned}
\partial_t q+\nabla_r\cdot(q v)+\nabla_\theta\cdot(q \omega) &= 0, \\
\partial_t(q\tilde{v}) +\nabla_r(q \tilde{v}\otimes\tilde{v})+ \nabla_\theta(q\tilde{v}\tilde{\omega}) &= qS_1[q] - (A^2/2)q\tilde{v}, \\
\partial_t(q \tilde{\omega}) +\nabla_r(q\tilde{v}\tilde{\omega})+\nabla_\theta(q\tilde{\omega}^2) &= qS_2[q] - (B^2/2)q\tilde{\omega}.
\end{aligned}\right.
\end{align}
In both equations \eqref{eq:macro_q_max} and \eqref{eq:macro_q_mono} all quantities depend on position and angle.

\subsection{$r$-marginals and mono-kinetic closure}
The above descriptions can be further reduced by integrating the mean field equation over $v,w$ and $\theta$. 
The resulting macroscopic quantities are 
\begin{align*}
\rho &=\int f(t,r,\theta,v,\omega)\,d\Omega_2, &
\rho \tilde{v}&=\int v f(t,r,\theta,v,\omega)\,d\Omega_2,\\
\rho\phi&=\int \theta f(t,r,\theta,v,\omega)\,d\Omega_2, &
\rho \tilde{\omega}&=\int \omega f(t,r,\theta,v,\omega)\,d\Omega_2,
\end{align*}
where $d\Omega_2=dvd\theta d\omega$. Again the distribution function is normalized and we have
\[
\int \rho(t,r) dr = 1\qquad\text{for all\; $t\ge 0$}.
\]
Integrating the balance equations \eqref{balance1} and \eqref{balance2} over $dw$, we
obtain 
\begin{align*}
\partial_t\rho+\nabla_r(\rho \tilde{v}) &= 0,\\
\partial_t(\rho \tilde{v}) + \nabla_r\cdot\left( \int v\otimes v f\,d\Omega_2\right) &= \int qS_1[q]\,d\theta -(A^2/2)\rho\tilde{v},\\
\partial_t(\rho \tilde{\omega}) + \nabla_r\cdot\left(\int v\omega f\,d\Omega_2\right) &= \int qS_2[q]\,d\theta -(B^2/2)\rho\tilde{\omega}.
\end{align*}
Integrating further the mean field equation \eqref{eq:meanfield} with respect to $\theta d\Omega_2$, we obtain
\[
 \partial_t(\rho\phi) + \nabla_r\cdot\left(\int\theta vf\,d\Omega_2\right) = \rho\tilde{\omega}.
\]
We now use a mono-kinetic and mono-angular closure function of the form
\[ 
 f\sim \rho\delta(\tilde{v}-v)\otimes\delta(\phi-\theta)\otimes \delta(\tilde{\omega}-\omega),
\]
to obtain the following hydrodynamic limit equations (without tilde), given by
\begin{align}\label{eq:macro_rho_mono}
\hspace*{-.7em}\left\{\begin{aligned}
\partial_t\rho+\nabla_r\cdot(\rho v)&=0, \\
\partial_t(\rho\phi)+\nabla_r\cdot(\rho \phi v)&=\rho \omega, \\
\partial_t(\rho v)+\nabla_r\cdot(\rho v\otimes v)&=-\gamma \rho (u-v) + K_1[\rho]\rho
+\rho\nabla_r V_1- (A^2/2)\rho v, \\
\partial_t(\rho \omega)+\nabla_r\cdot(\rho v\omega)&= -\bar{\gamma}\rho(g(\phi,u) - \omega) + K_2[\rho]\rho + \rho\nabla_\theta V_2(\phi) - (B^2/2)\rho\omega.
\end{aligned}\right.
\end{align}
with the interaction terms
\begin{align*}
 K_1[\rho] &= \frac{1}{m}\int\nabla_{1} U(r,\bar{r},\phi(r),\phi(\bar{r}))\rho(\bar{r})\,d\bar{r}, \\
 K_2[\rho] &= \frac{1}{I_c}\int\nabla_{3}U(r,\bar{r},\phi(r),\phi(\bar{r}))\rho(\bar{r})\,d\bar{r},
\end{align*}
where $\nabla_k U$ denotes the derivative of $U$ with respect to the $k$-th component.

\begin{remark}
Using a suitable scaling, one can also derive the diffusive limit of the hydrodynamic equations \eqref{eq:macro_rho_mono}. In this case, one obtains the system
\begin{align}\label{eq:difflimit}
\begin{aligned}
\partial_t\rho + \nabla_r\cdot\left(\sigma_1\rho\big(\gamma u-K_1[\rho]- \nabla_r V_1\big)\right) &= 0,\\
\partial_t(\rho\phi)+\nabla_r\cdot\left(\sigma_1\rho\phi\big(\gamma u-
K_1[\rho]- \nabla_r V_1\big)\right) &= \sigma_2\rho\big(\bar{\gamma} g_u(\phi)-K_2[\rho] - \nabla_\theta V_2(\phi) \big).
\end{aligned}
\end{align}
with $\sigma_1=2/(2\gamma-A^2)$ and $\sigma_2=2/(2\bar{\gamma}-B^2)$.
\end{remark}

% % % % % % % % % % % % % % % % % % % % % %
% % % Numerical Methods and Examples  % % %
% % % % % % % % % % % % % % % % % % % % % % 
\section{Numerical Methods and Examples}
\label{sec:num ex}
In this section we consider the numerical solution of the different models derived above and compare these models in various test cases. 
% % % Numerical Methods % % %
\subsection{Numerical Methods}
To solve the microscopic system \eqref{eq:micro} we use a first-order leapfrog algorithm for the deterministic part and the Euler-Maruyama scheme for the stochastic terms \cite{kloeden1992numerical}. As for the hydrodynamic systems \eqref{eq:macro_q_max}, \eqref{eq:macro_q_mono}, \eqref{eq:macro_rho_mono} and the diffusive limit \eqref{eq:difflimit} we use a Godunov splitting scheme to split conservative part of the system and source terms \cite{toro1999riemann}. The conservation law is solved with the FORCE scheme that may be found for example in \cite{toro1999riemann}, and a dimensional Godunov splitting in each dimensions.
The source terms, on the other hand, are treated with a simple first-order implicit Euler method.

% % % Numerical Example % % %
\subsection{Convergence towards a stationary solution}
In this test case we neglect the surrounding fluid by setting $\gamma \equiv \bar \gamma \equiv 0$ and use suitable exterior potentials $V_1$, $V_2$, such that the resulting equations have a stationary distribution.

Note that the stationary equation for the mean field equation \eqref{eq:meanfield} and equations \eqref{eq:macro_q_max} with Maxwellian closure for $m\equiv I_c\equiv 1$ may be given by the integral equation
\[ 
 q_{stat}=Z_q^{-1}\exp\big(-V_1-V_2-K[q_{stat}]\big),
\]
where $Z_q$ is a normalizing factor, such that $\int q_{stat}\,dr d\theta = 1$, and
\[
 K[q] = \frac{1}{m}\iint \nabla_{r}U(r,\bar{r},\theta,\bar{\theta})q(\bar{r},\bar{\theta})d\bar{r}d\bar{\theta}.
\]

In the following, we investigate and compare the convergence to a stationary solution of the microscopic system \eqref{eq:micro} and of the 
$q$-equation with Maxwellian closure \eqref{eq:macro_q_max}. We choose the computational spatial domain  $\Omega=[-6,~6]\times[-6,~6]$, and the exterior potentials $V_1(r)=\vert r\vert^2/2$ and $V_2(\theta)=\sin(\theta)$. 
Moreover,  we choose
\[
 L=1,\qquad D=0.5,\qquad \epsilon_0=1.
\]
As initial conditions for the microscopic case we choose randomly distributed positions of the particles inside 
$\Omega_0=[1.5,~3.5]\times[1.5,~3.5]$ and
\[
\theta_i^o=0,\qquad v_i^o=(0,0)^\top, \qquad \omega_i^o=0,\qquad \text{for~} i=1,\ldots,200.
\]
The corresponding hydrodynamic initial states for equations \eqref{eq:macro_q_max} are
\begin{gather*}
q(0,r,\theta)=\begin{cases}
\frac{1}{4 \Delta \theta}, &\text{if~} (r,\theta)\in\Omega_0\times [ 0,\Delta \theta ],\\
0, &\text{else},
\end{cases}\\
v(0,r,\theta)=(0,0)^\top,\qquad
\omega(0,r,\theta)=0,\qquad
\forall (r,\theta)\in\Omega\times[0,~2\pi].
\end{gather*}
The discretization of the spatial domain for the hydrodynamic problem  is given by $60\times60\times8$ grid 
cells. We compare the difference to the stationary solution in the $L^2$-norm and show the convergence for 
$A=B=0,0.2,1,2,3,5,7,10$. 

\begin{figure}
\includegraphics[width=0.49\textwidth]{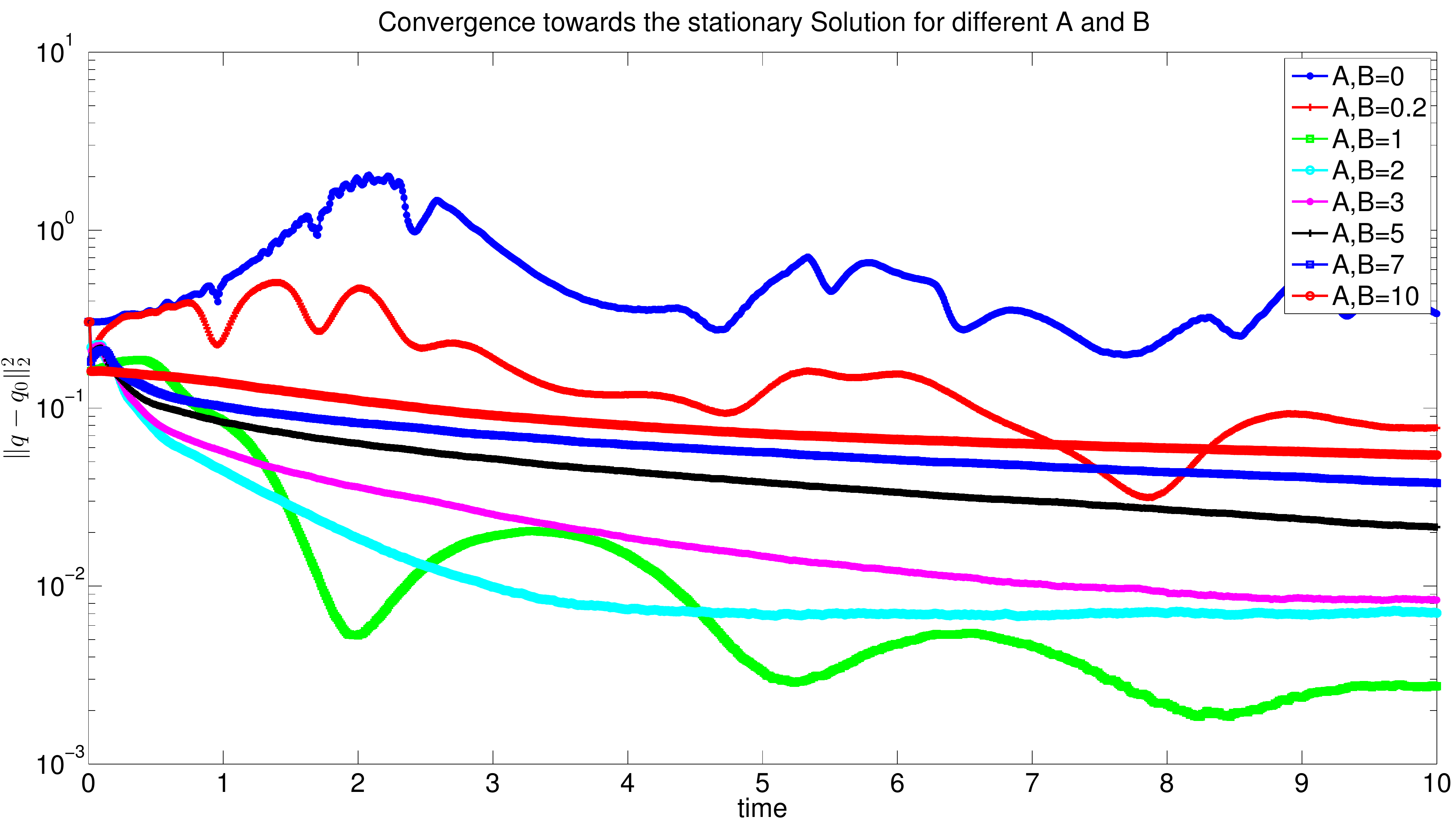}
\includegraphics[width=0.49\textwidth]{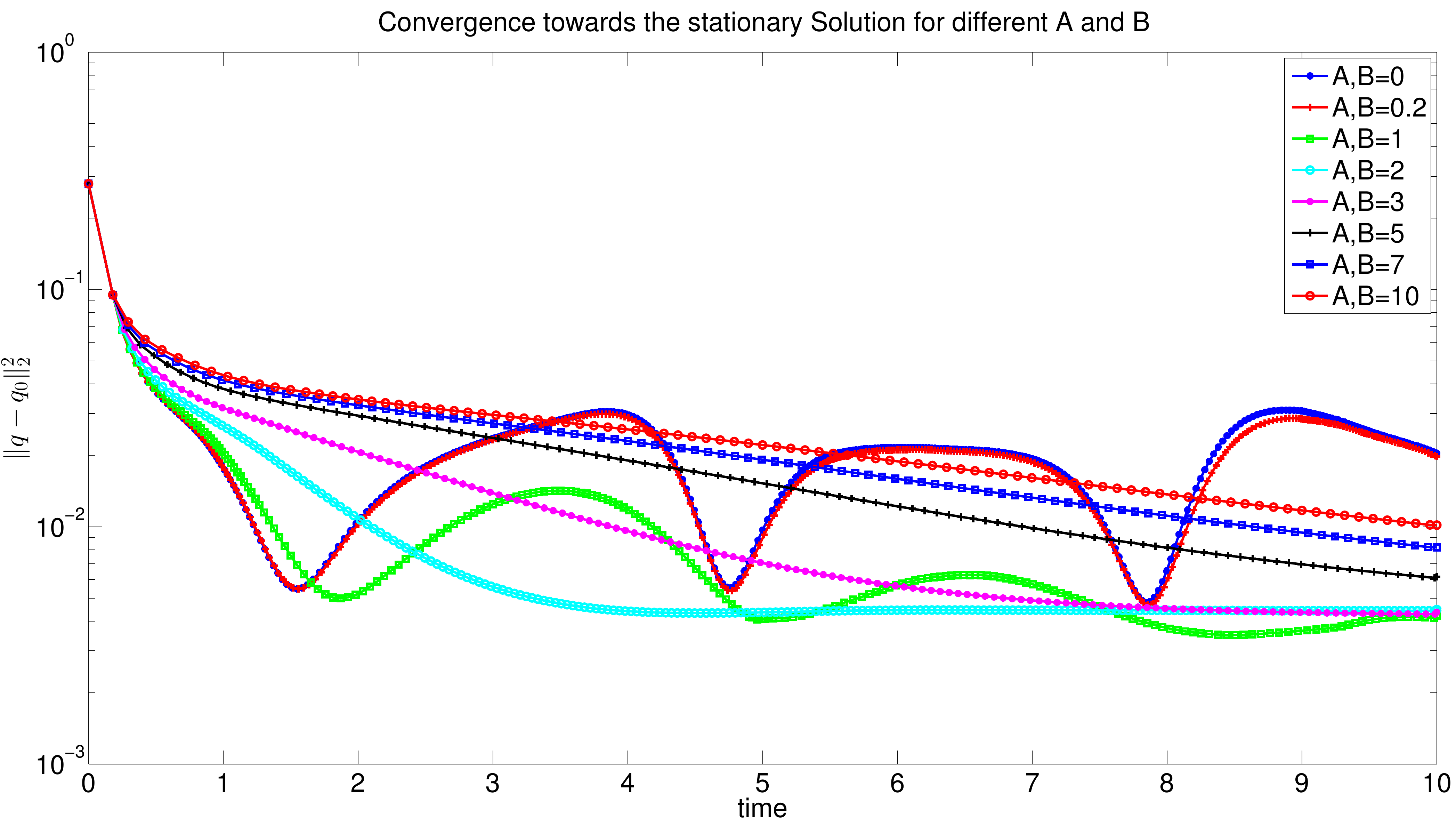}
\caption{Convergence towards the stationary solution for microscopic equation and  q-equation with Maxwellian closure.  Left: System \eqref{eq:micro}; Right:  \eqref{eq:macro_q_max}.}
\label{fig:GG_konv_q_inter}
\end{figure}

The results are displayed in Figure~\ref{fig:GG_konv_q_inter}. For values of $A=B$ less than one, we observe an oscillatory behavior of the trend to equilibrium in the $L^2$-norm, while for values larger than one, the decay happens monotonically. Also, for $A=B$ larger than one, the trend to equilibrium for both models exhibit similar behavior. Moreover, for intermediate values of $A=B$, we obtain fastest decay, which suggests signs of hypocoercivity, compare \cite{klar2014approximate} for similar results for a simpler case. 

Note that the influence of the interaction term is relatively small in the present case. In fact, the numerical results neglecting the interaction show a similar behavior.

\begin{figure}
\center
\includegraphics[width=0.5\textwidth]{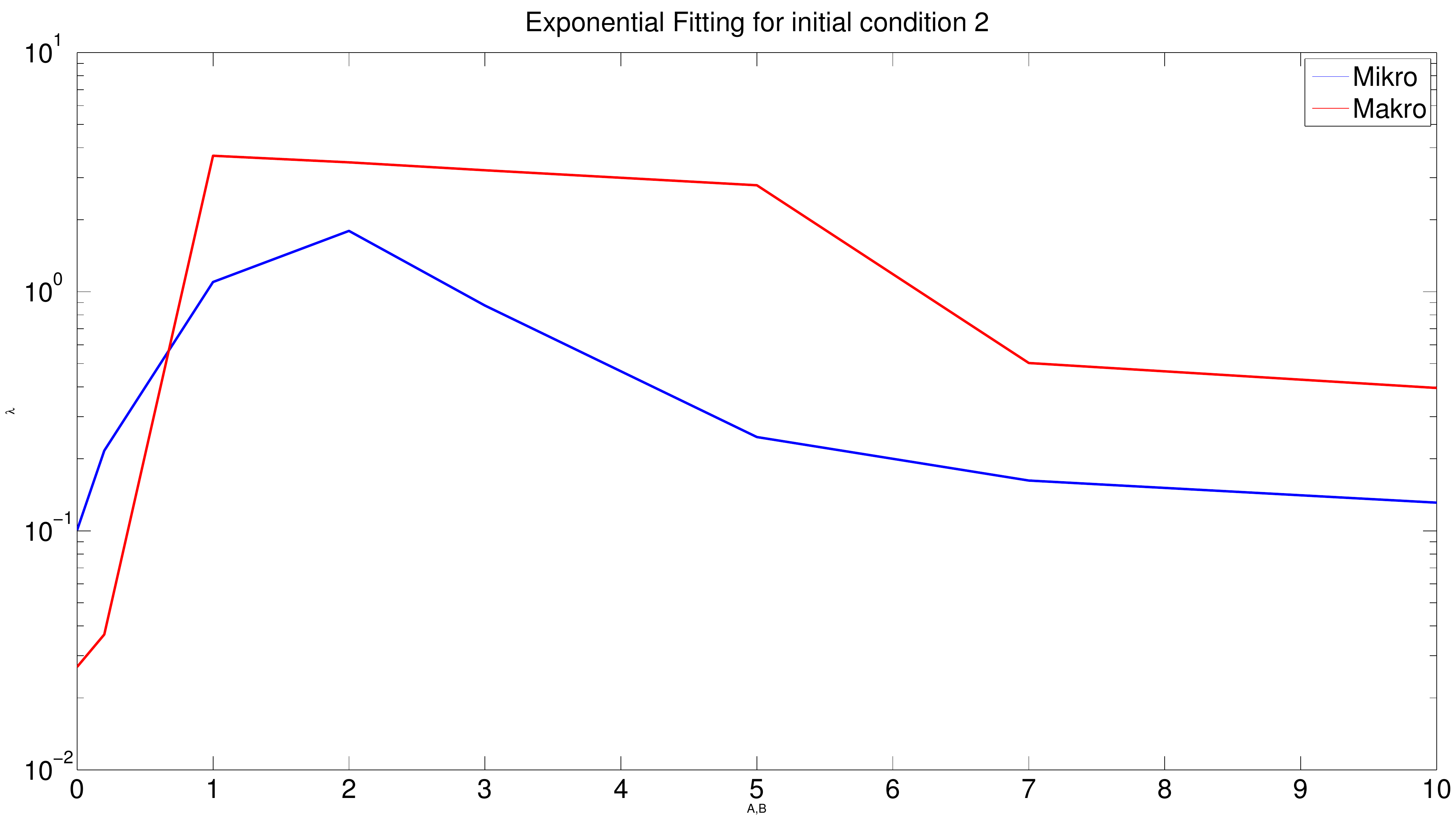}
\caption{Decay rate $\lambda$ for different $A,B$ using exponential fitting of decay curves.}
\label{fig:GG_lam}
\end{figure}

The corresponding decay rates of the solutions to the stationary solution are determined by a fitting procedure assuming an exponential decay, i.e., $e^{-\lambda t}$, with decay rate $\lambda>0$, and fitting the exponent $\lambda$ to the experimentally determined decay functions. The decay rate $\lambda$ as a function of $A,B$ are depicted in Figure \ref{fig:GG_lam}.

\subsection{Top-bottom flow}
\begin{figure}[h!]
	\centering
\includegraphics[width=0.5\textwidth]{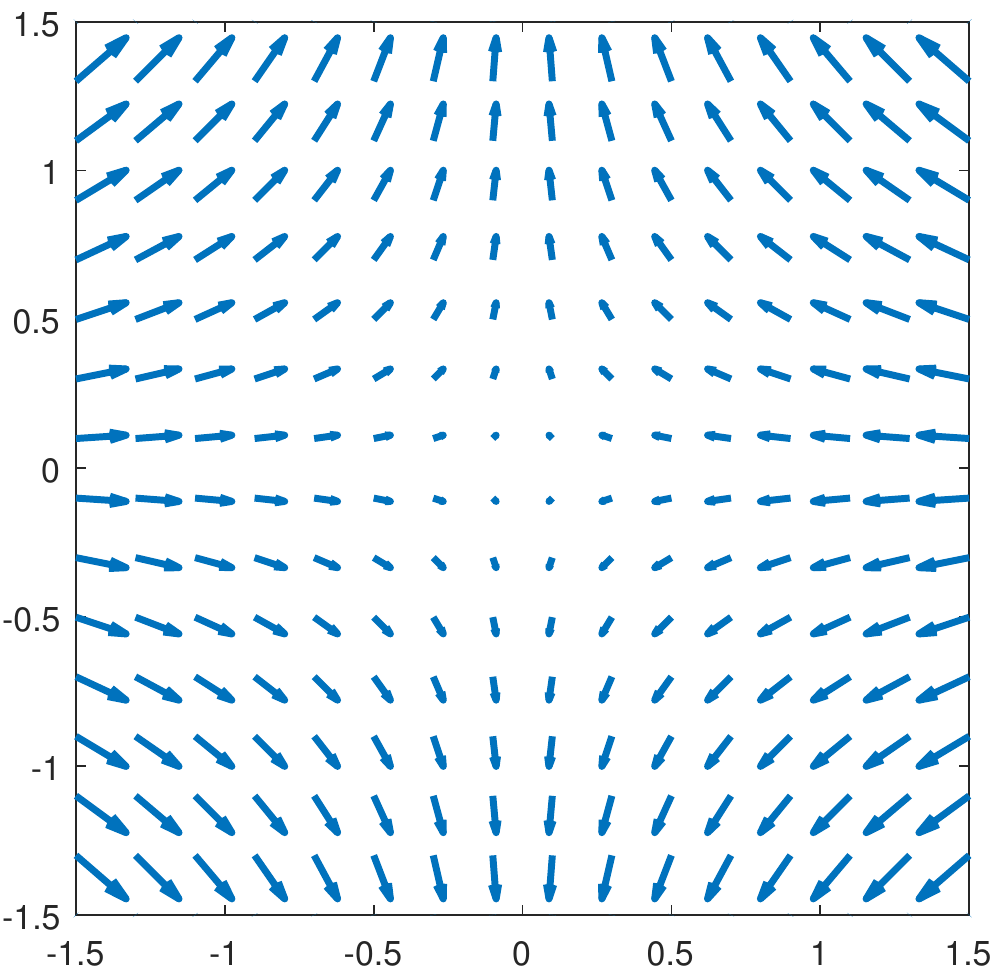}
\caption{Top-bottom example: Velocity field of the surrounding fluid}
\label{fig:vek_feld2}
\end{figure}

In this section we consider an example with surrounding fluid.
We assume that the velocity of the fluid is constant in time and not changed by the particles.
We choose $u(x,y)=(-x,y)^\top$ (see Figure~\ref{fig:vek_feld2}).
Moreover, we set $V_1=V_2=0$ and $A=B=0$ in this case.

\begin{figure}[h!]
\includegraphics[width=0.3\textwidth]{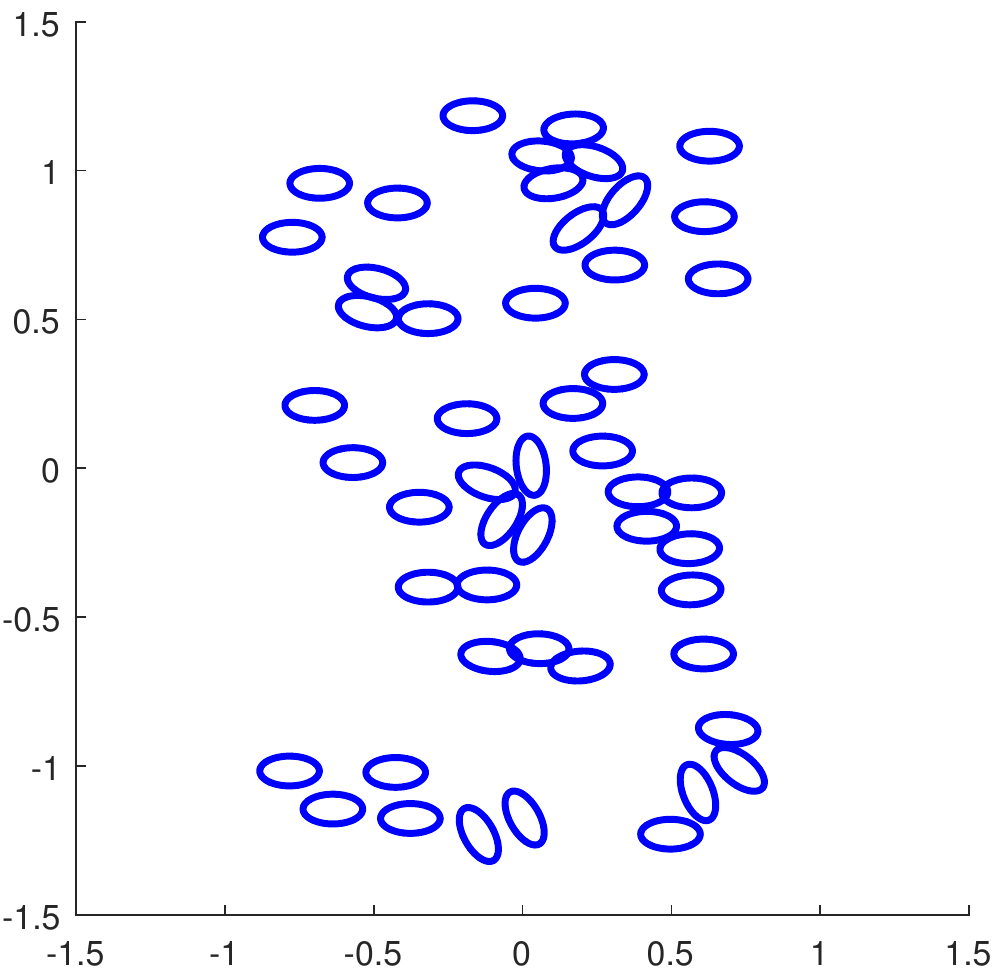}
\includegraphics[width=0.3\textwidth]{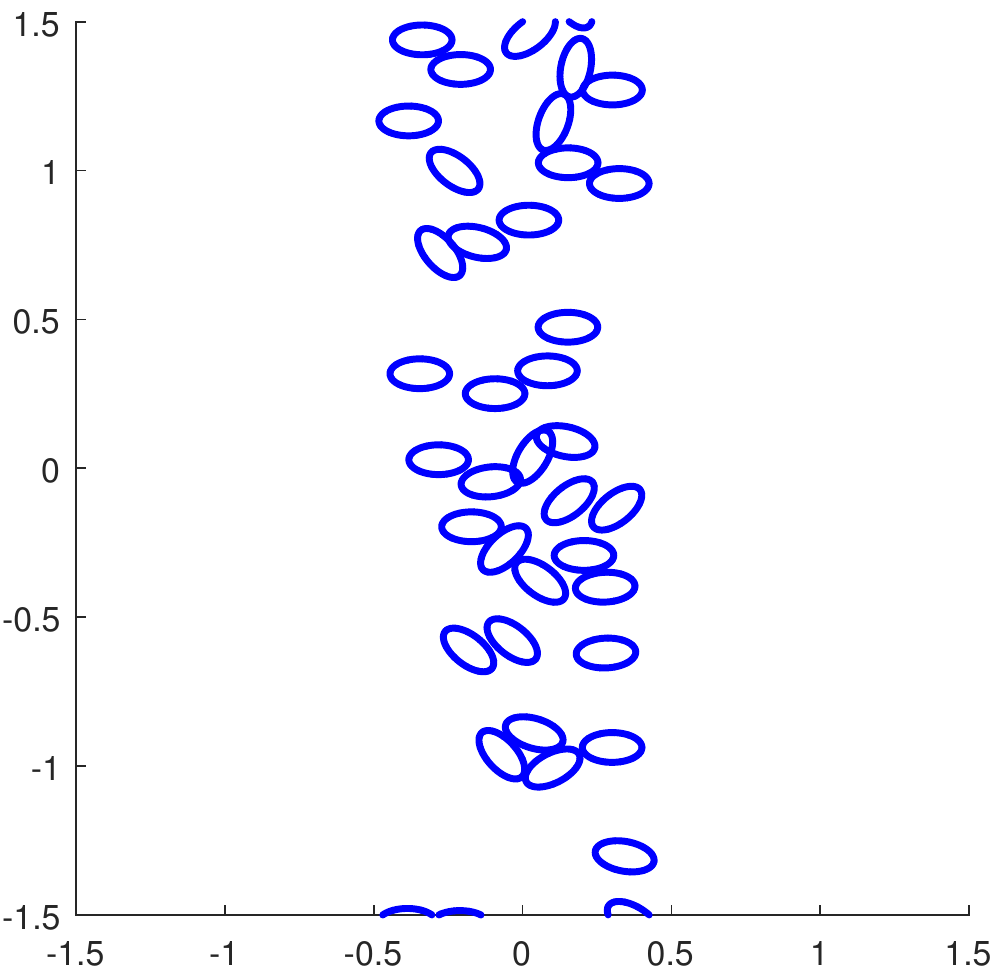}
\includegraphics[width=0.3\textwidth]{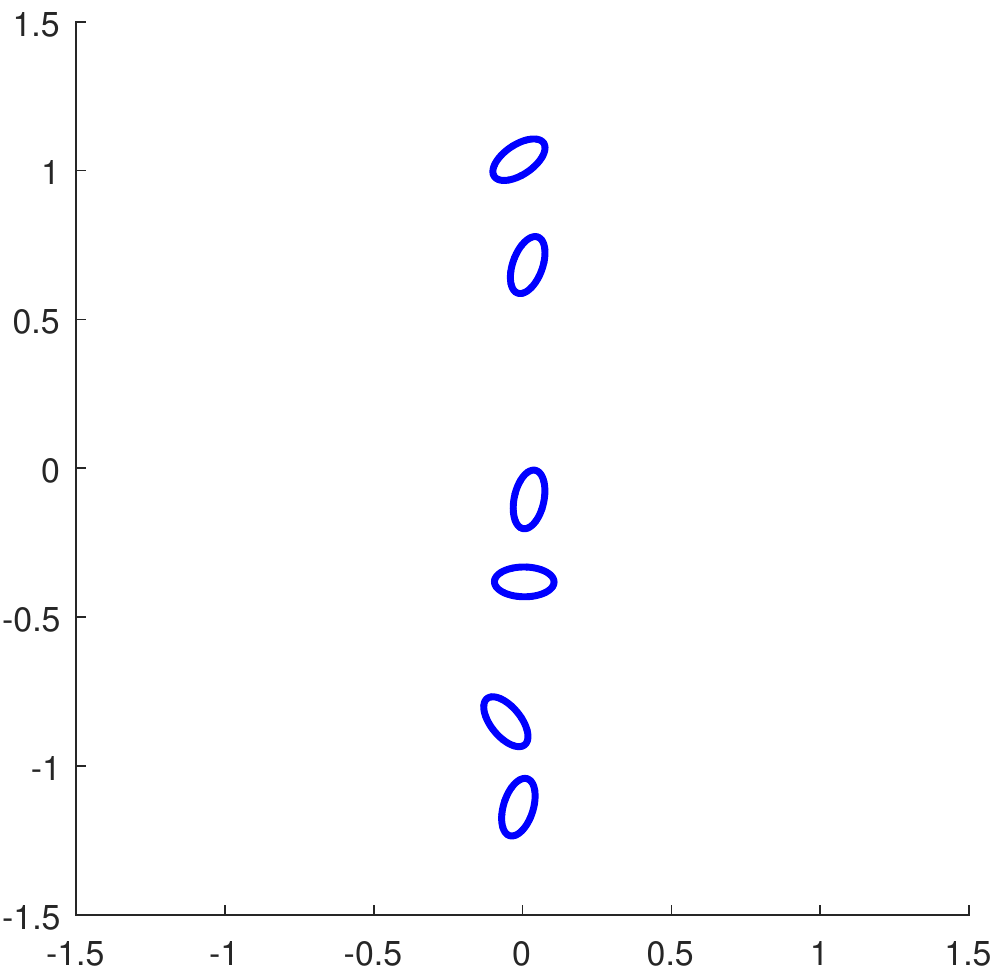}
\caption{Illustration of the microscopic evolution in the top-bottom example with interaction and without stochastic forces at times T$\approx$0.75 (left) and T$\approx$1.5 (middle) and T=5 (right)}
\label{fig:bsp1_mitill}
\end{figure}

\begin{figure}[h!]
\includegraphics[width=0.3\textwidth]{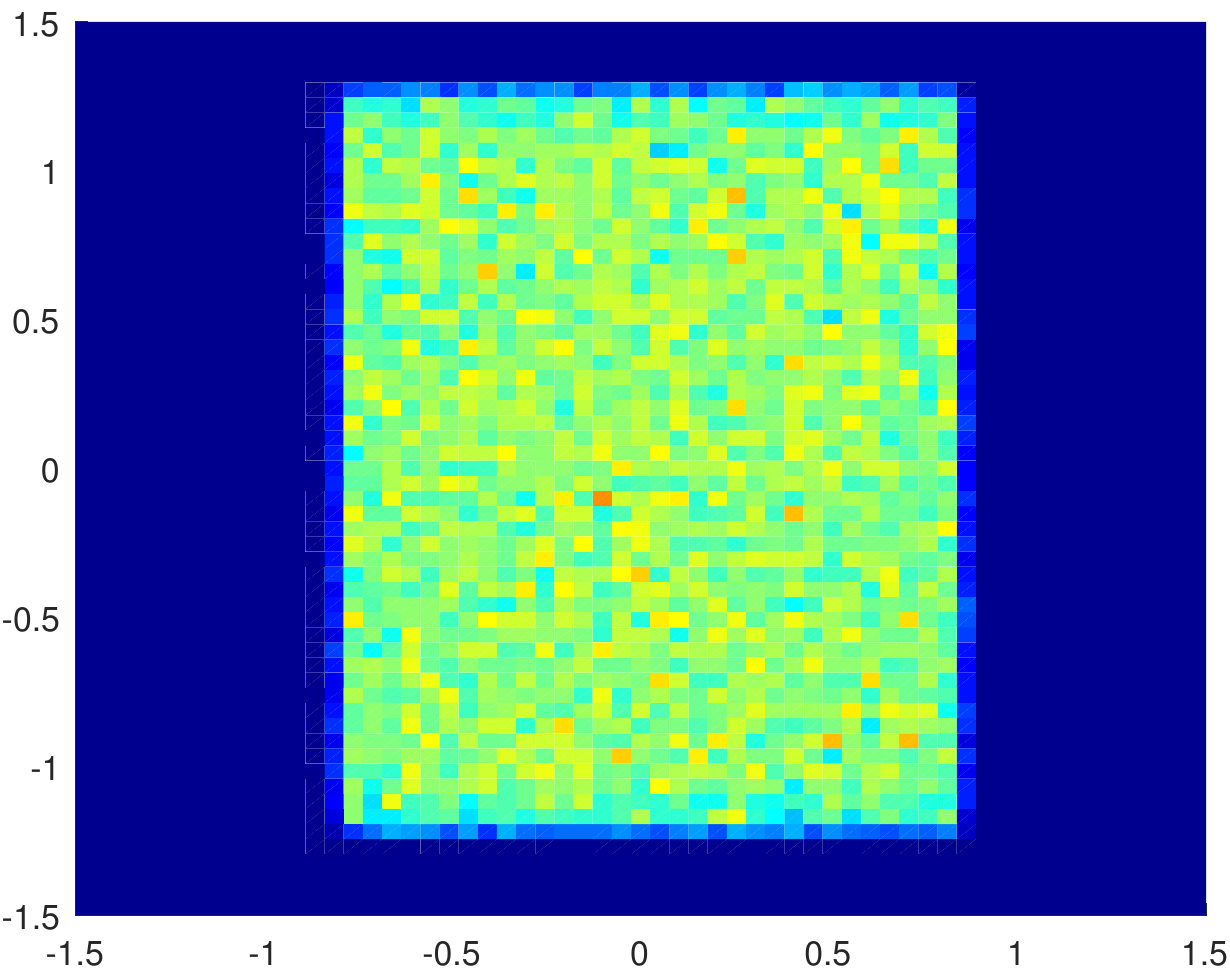}
\includegraphics[width=0.3\textwidth]{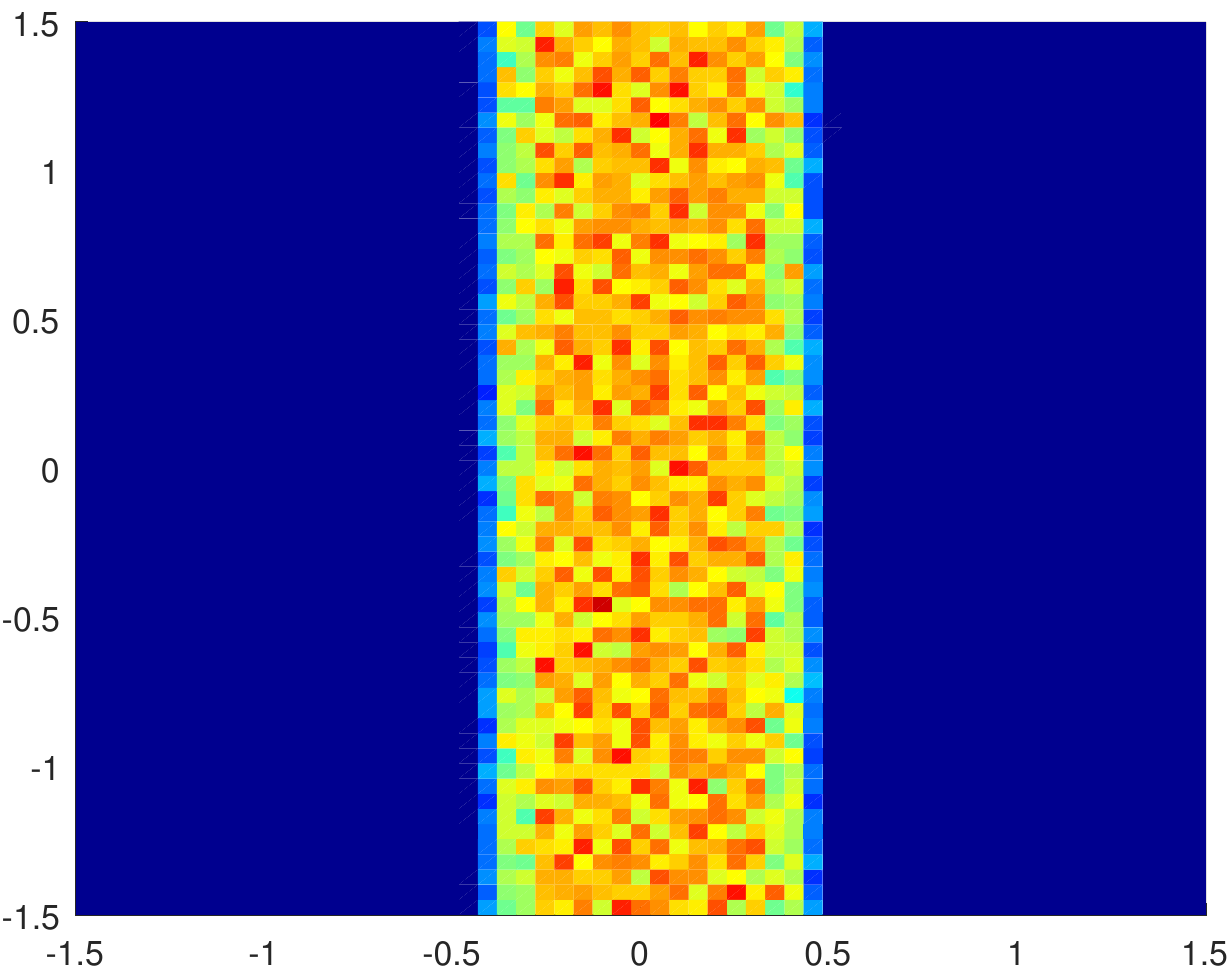}
\includegraphics[width=0.3\textwidth]{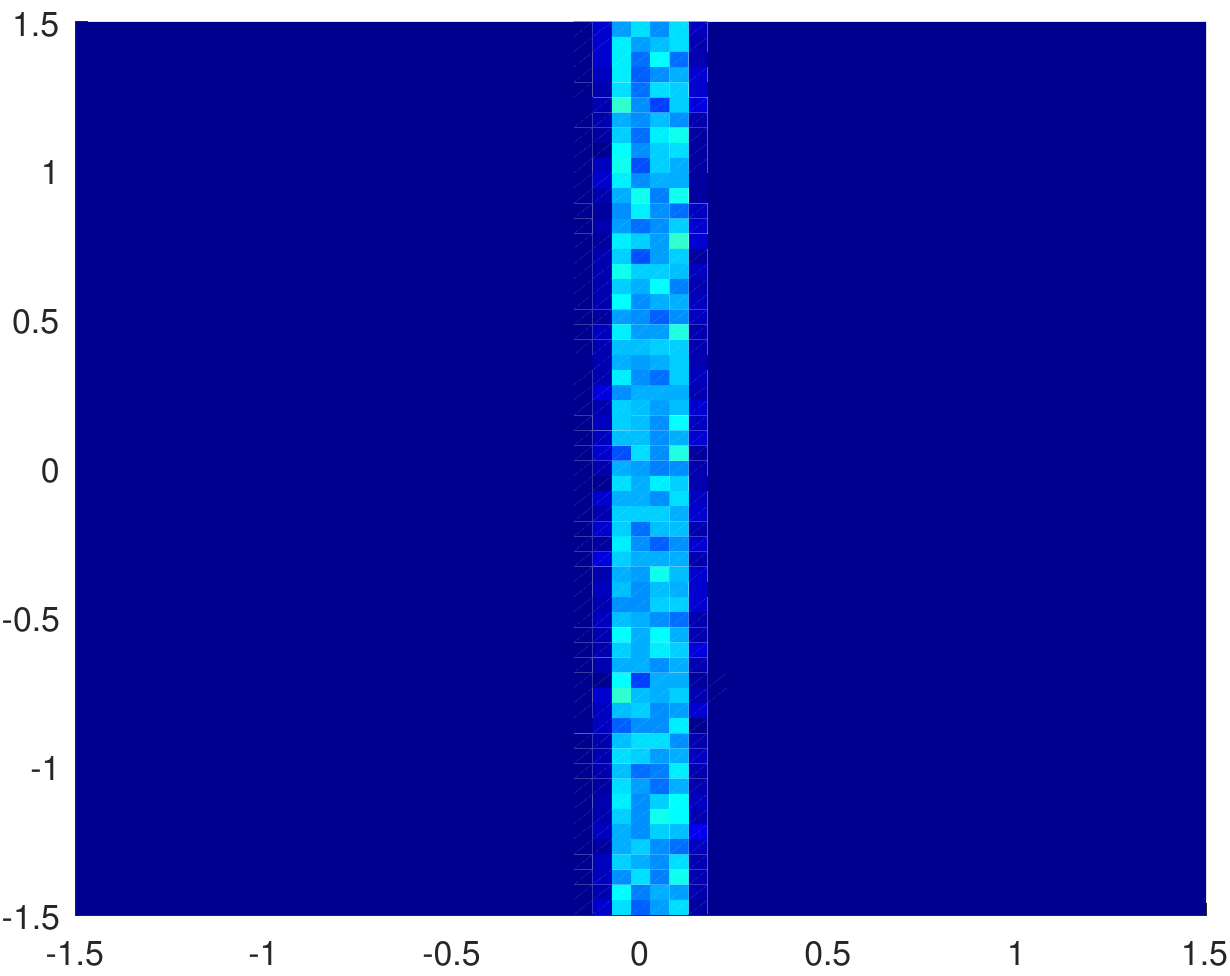}\\
\includegraphics[width=0.3\textwidth]{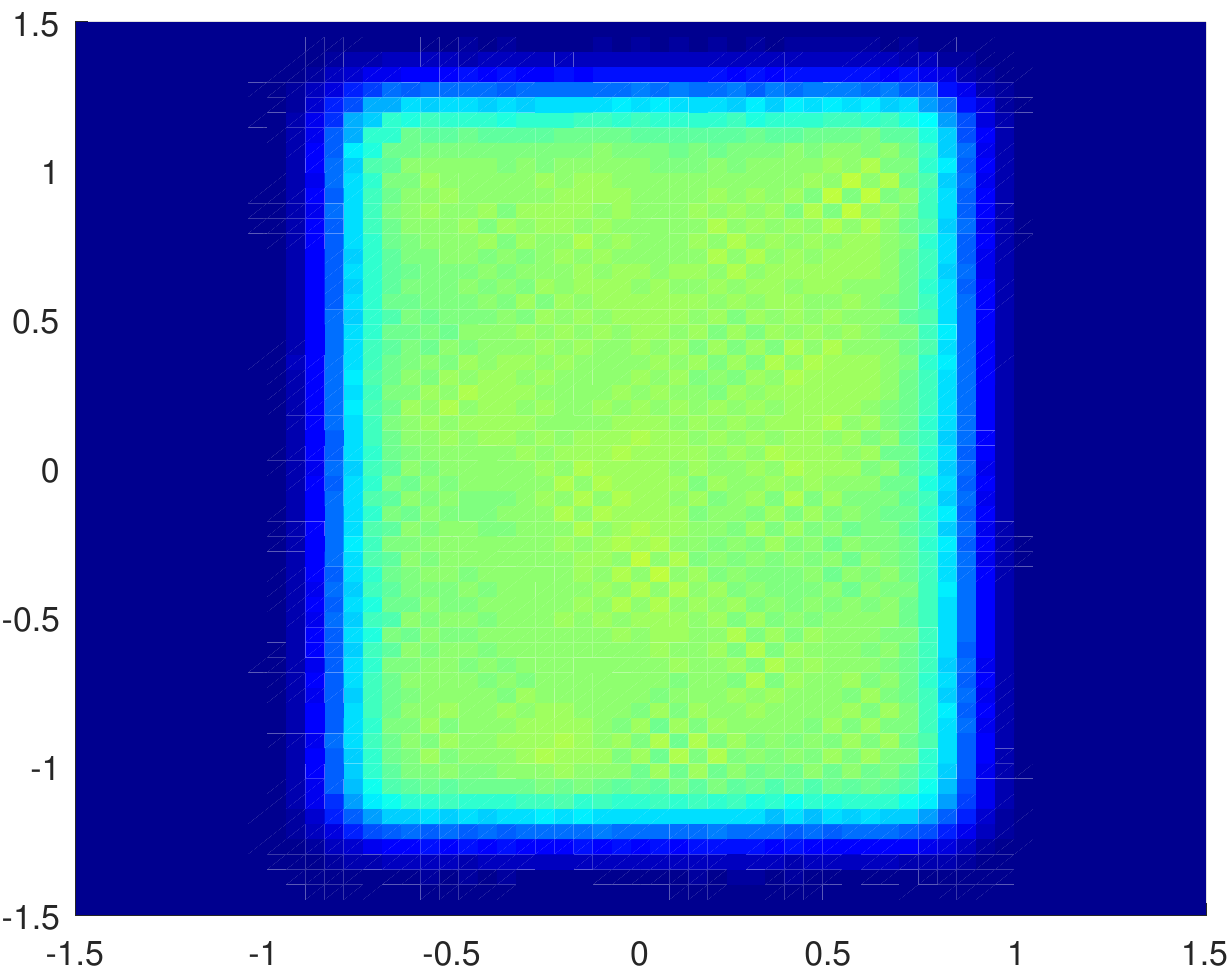}
\includegraphics[width=0.3\textwidth]{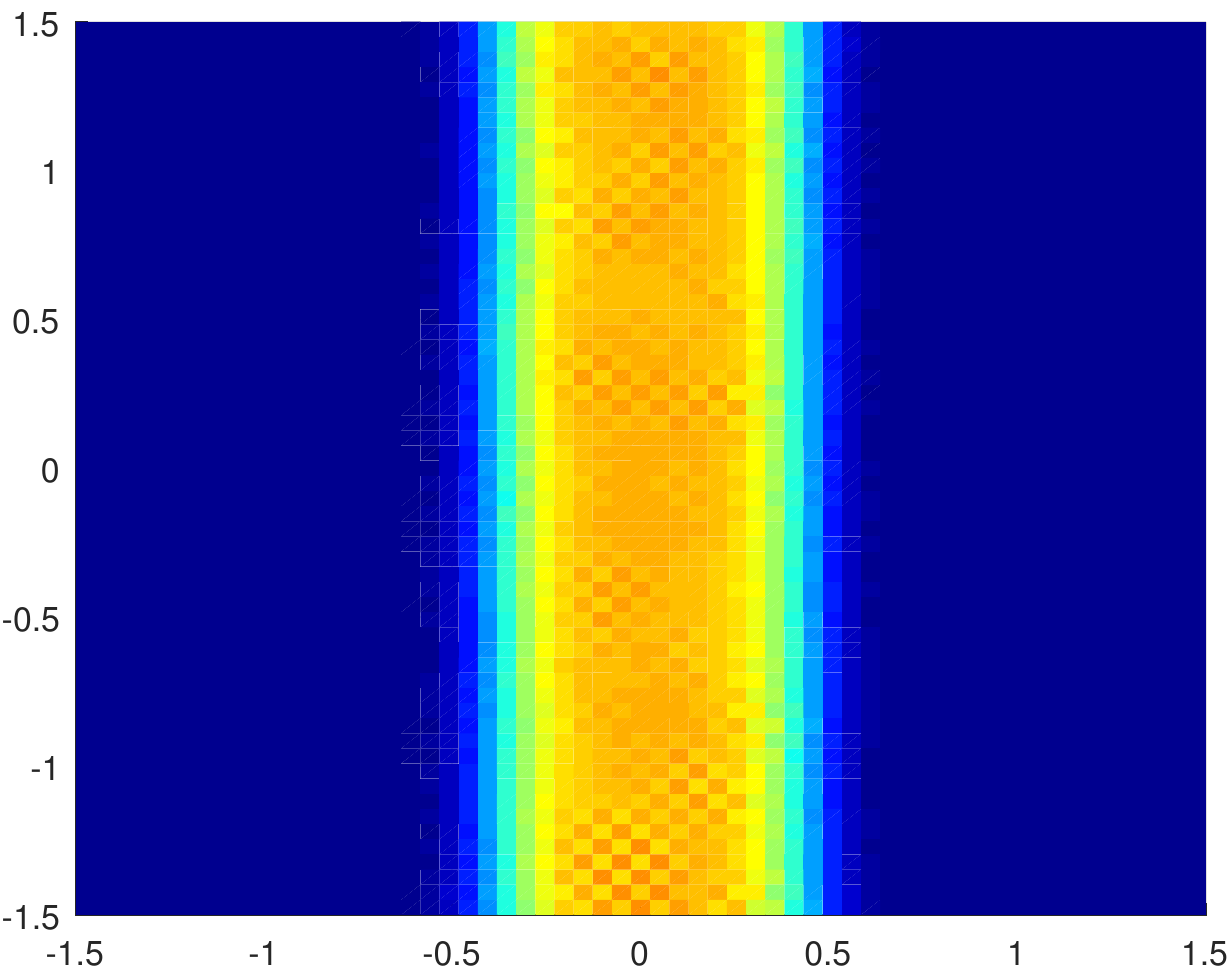}
\includegraphics[width=0.3\textwidth]{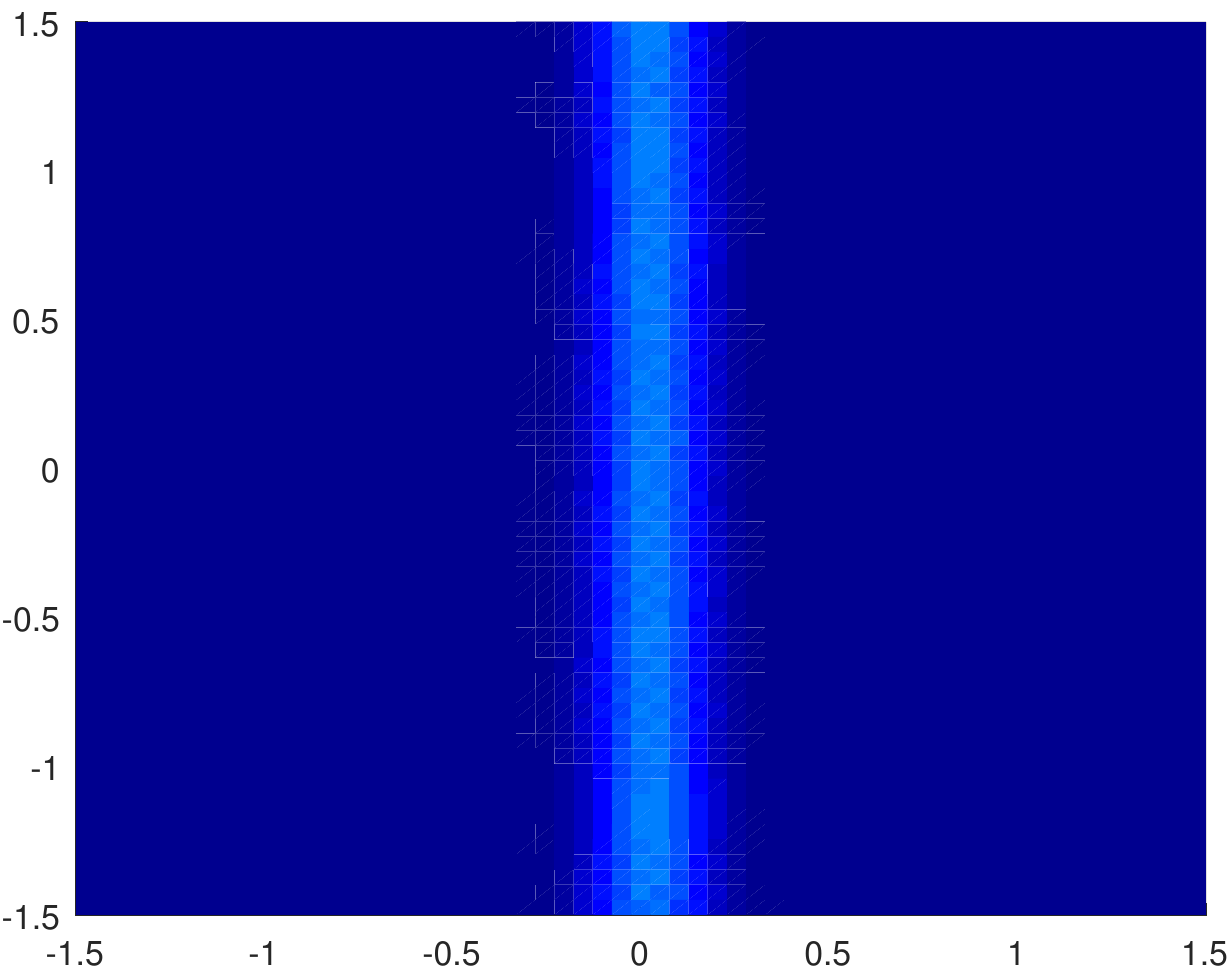}\\
\includegraphics[width=0.3\textwidth]{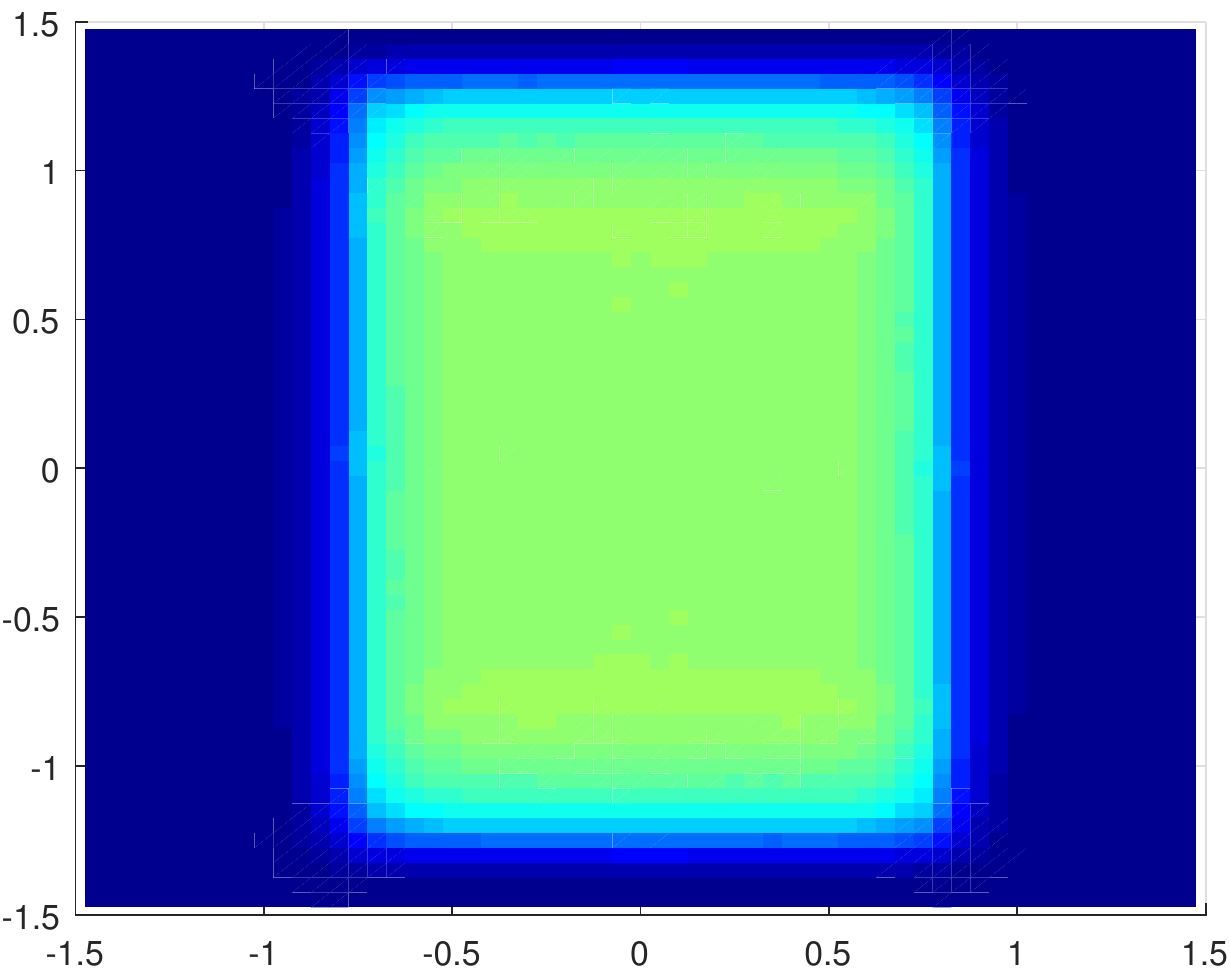}
\includegraphics[width=0.3\textwidth]{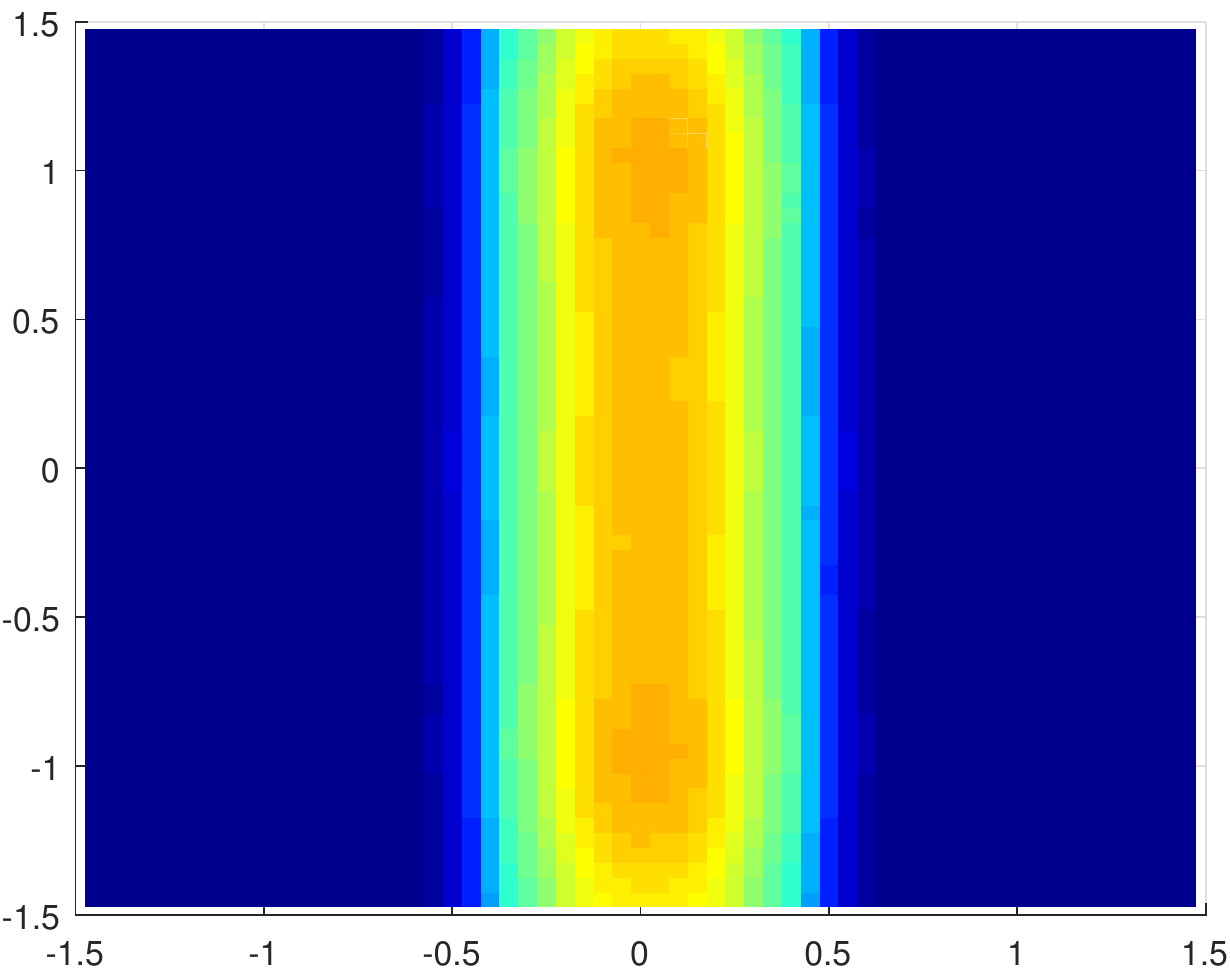}
\includegraphics[width=0.3\textwidth]{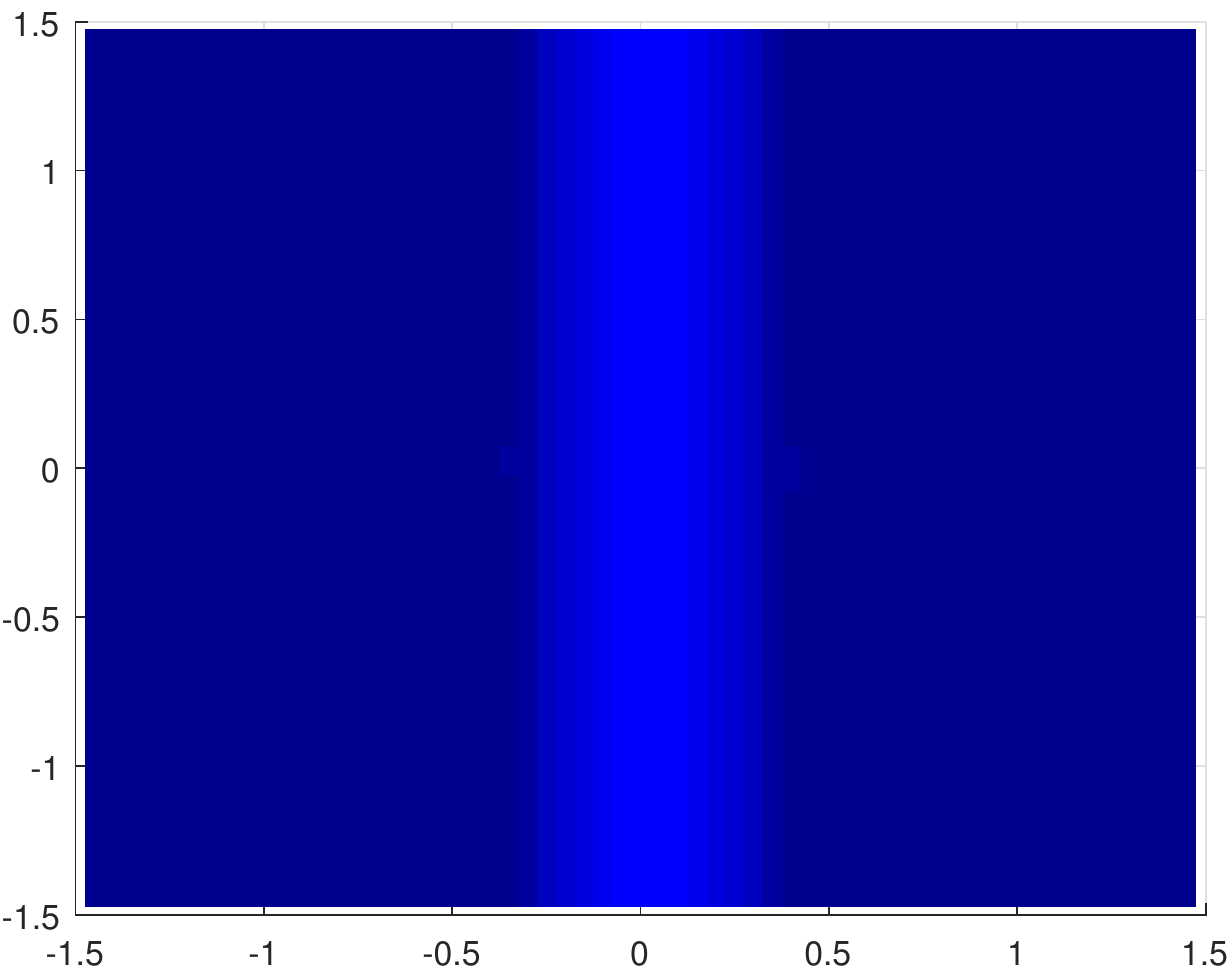}\\
\includegraphics[width=0.3\textwidth]{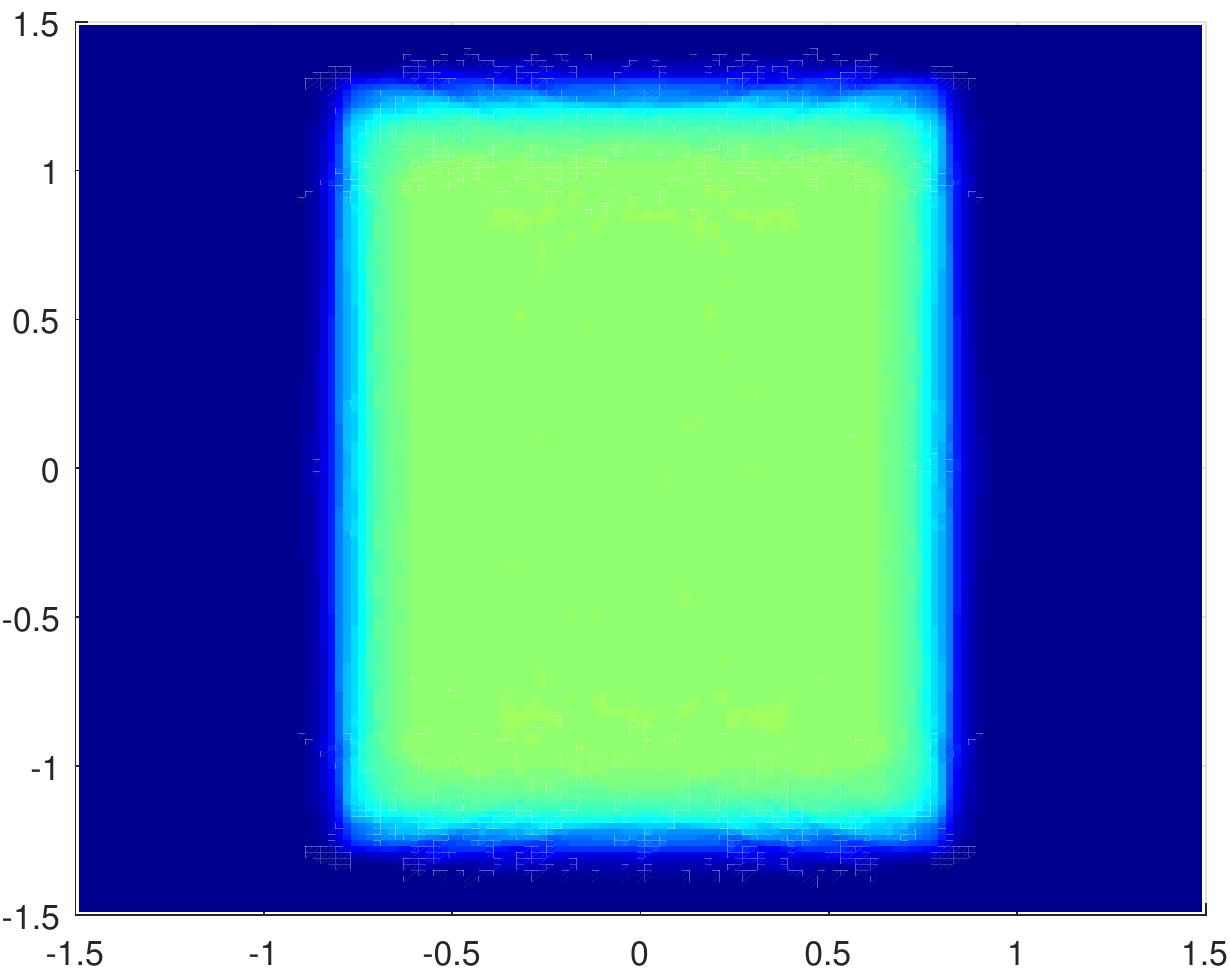}
\includegraphics[width=0.3\textwidth]{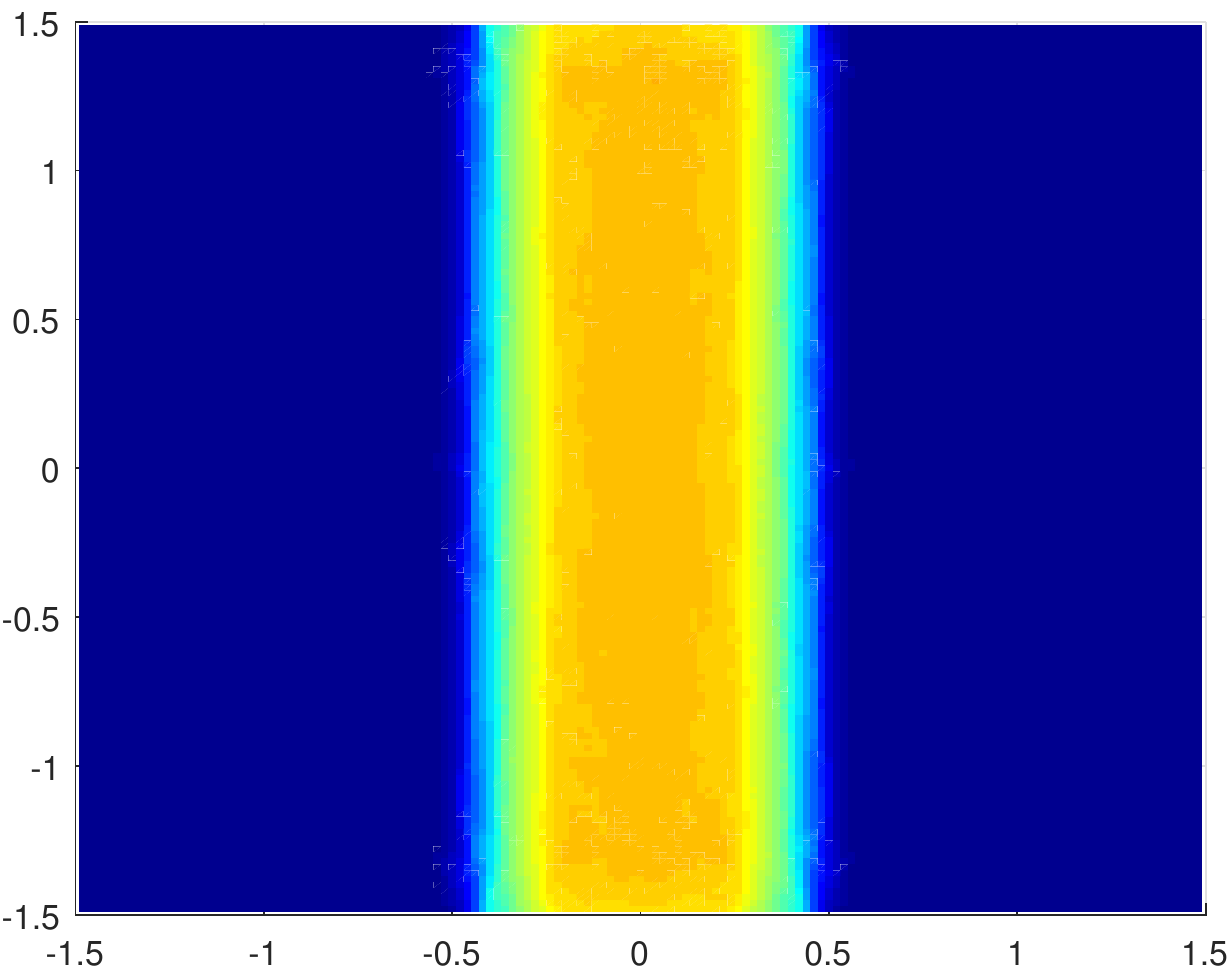}
\includegraphics[width=0.3\textwidth]{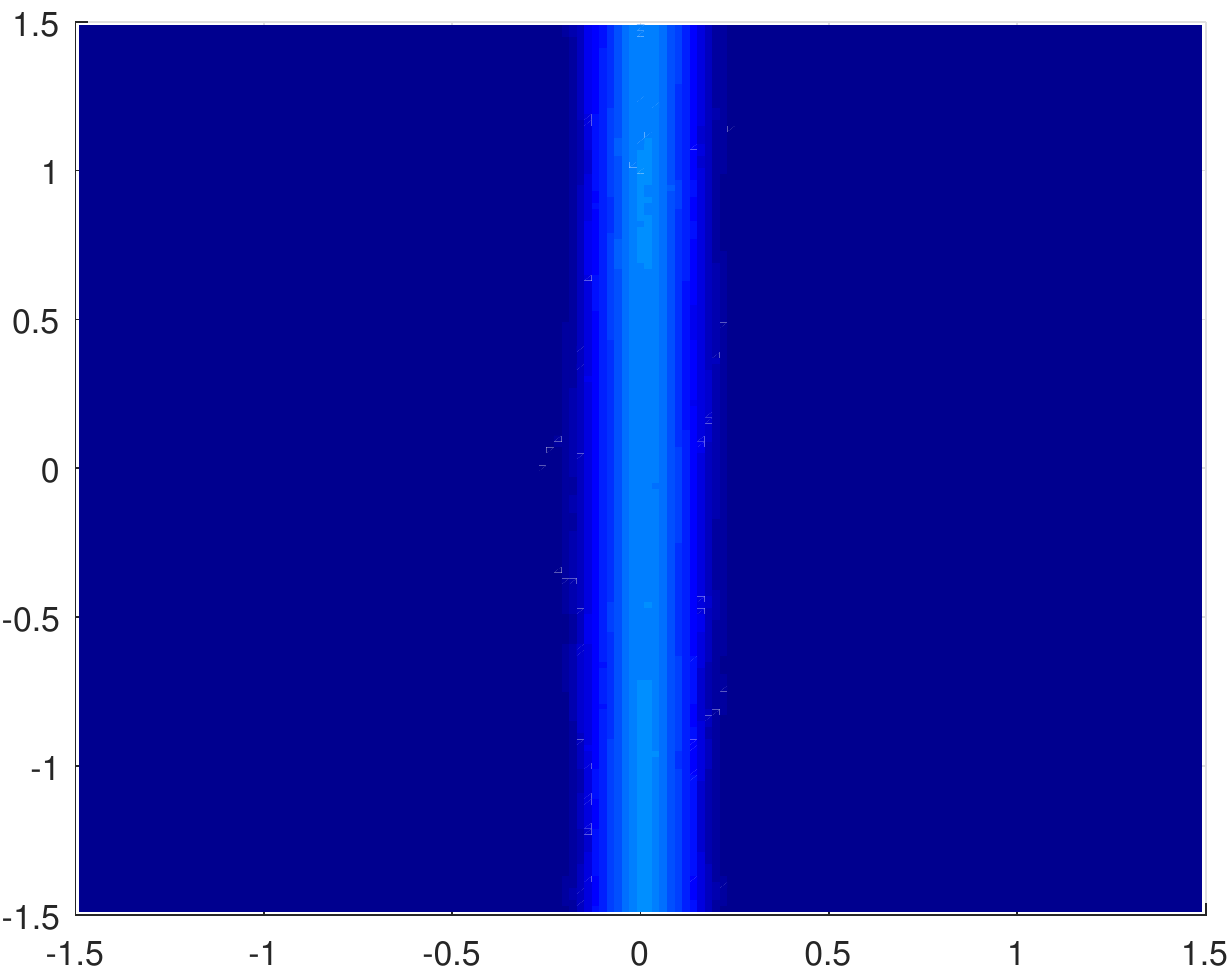}\\
\includegraphics[width=0.3\textwidth]{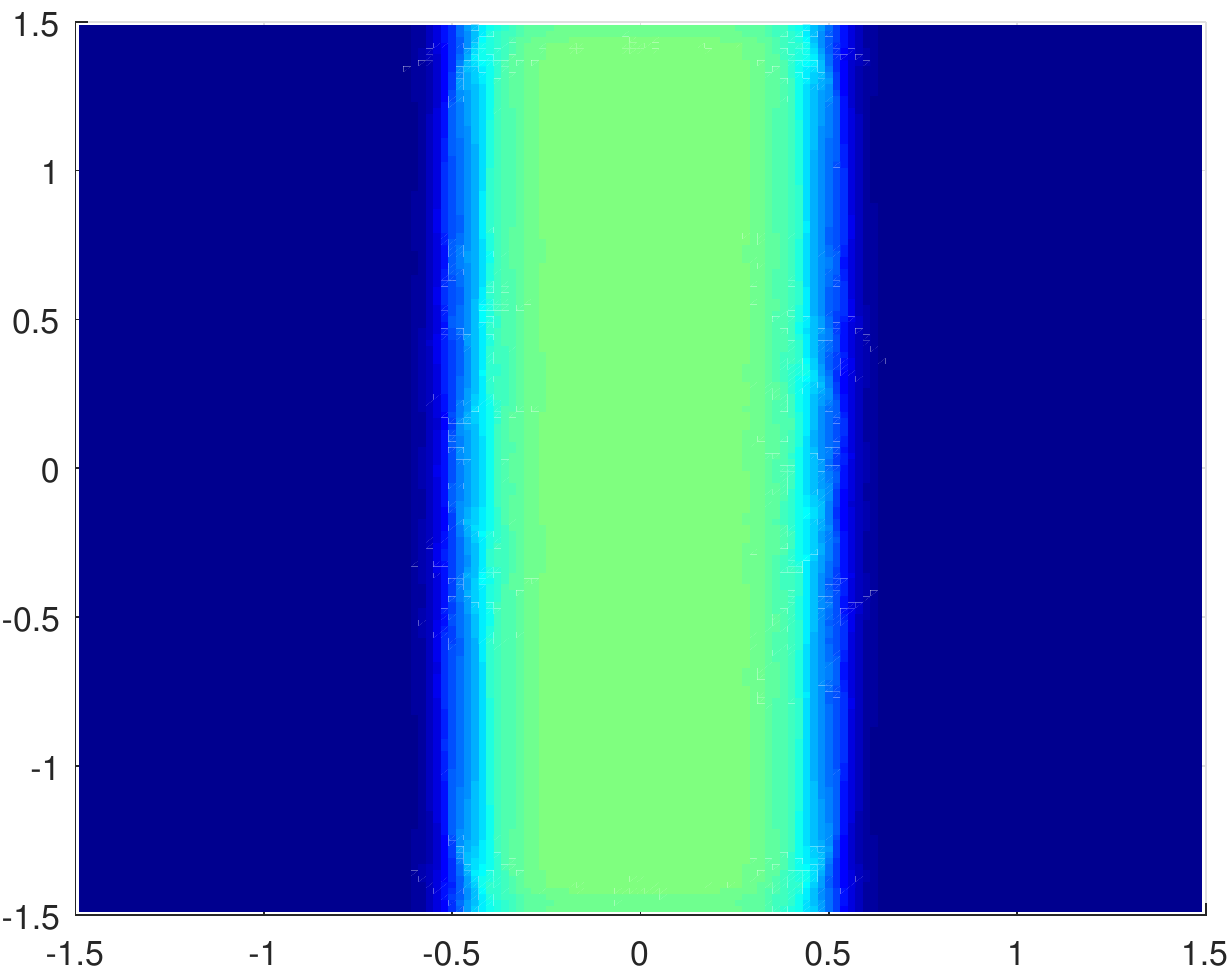}
\includegraphics[width=0.3\textwidth]{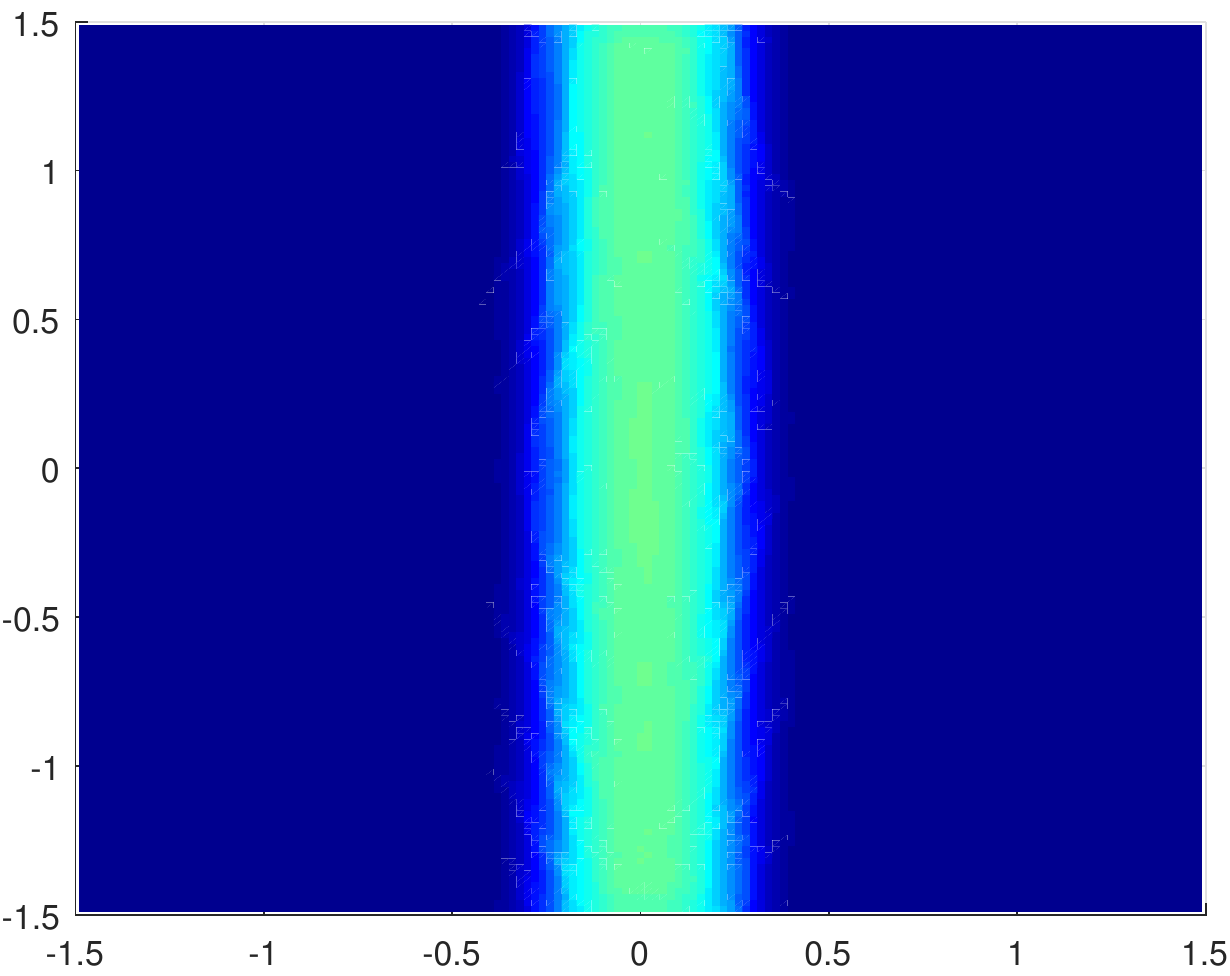}
\includegraphics[width=0.3\textwidth]{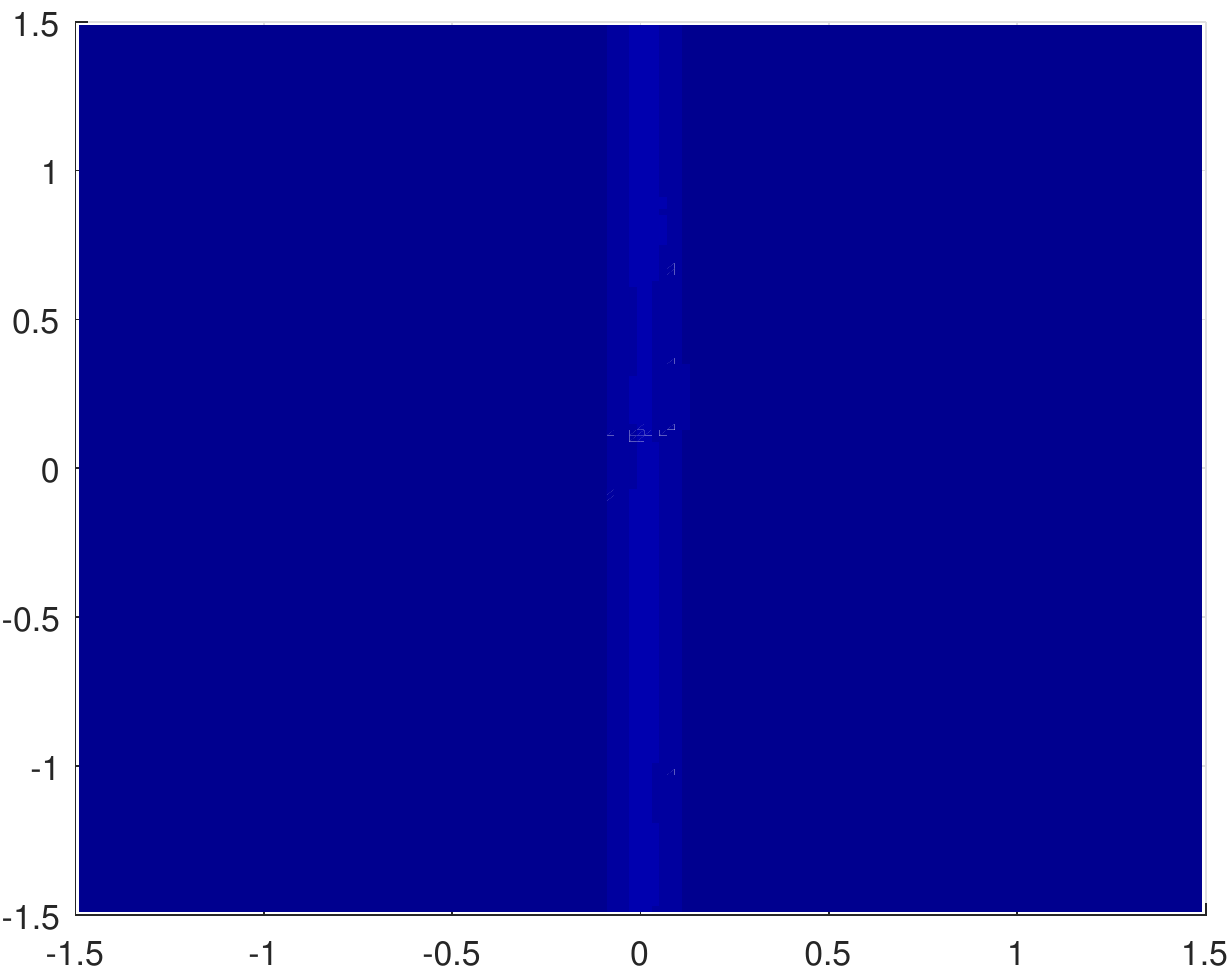}\\
	\centering\begin{tikzpicture}
		\begin{axis}[
		    hide axis,
		    scale only axis,
		    height=0pt,
		    width=0pt,
		    colormap/jet,
		    colorbar horizontal,
		    point meta min=0,
		    point meta max=0.5,
		    colorbar style={
		        width=10cm,
		        xtick={0,0.1,0.2,...,0.5}
		    }]
		    \addplot [draw=none] coordinates {(0,0)};
		\end{axis}
		\end{tikzpicture}
\caption{Top-bottom example with interaction and without stochastic at times T$\approx$0.75 (left) and T$\approx$1.5 (middle) and T=5 (right), for (from top to bottom)  microscopic histogram, smoothed microscopic histogram, mono-kinetic $q$, $\rho$-equations and diffusive $\rho$-equation}
\label{fig:bsp1_mit}
\end{figure}
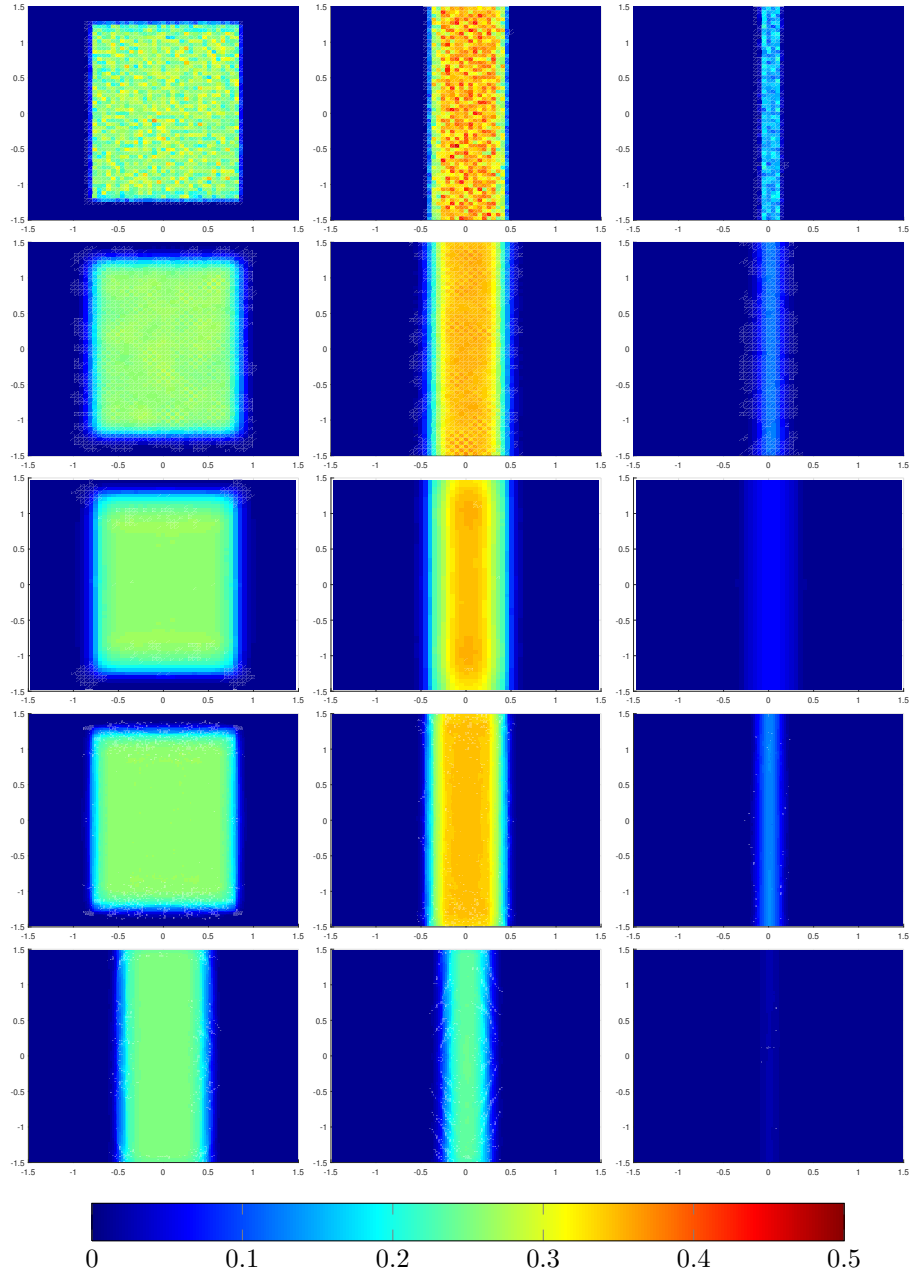

We compare the solutions of the microscopic equation \eqref{eq:micro}, the $q$-equation with mono-kinetic closure \eqref{eq:macro_q_mono}, the  $\rho$-equation \eqref{eq:macro_rho_mono} and the diffusive limit equation \eqref{eq:difflimit}.  For the PDE systems we use Neumann boundary conditions in space and periodic boundary conditions in the orientation.  
For all models we consider a quadratic domain $\Omega=[-1.5,~1.5]\times[-1.5,~1.5]$ and the parameters
\begin{align*}
L&=0.1, & \gamma&=1, & m&=1, & \epsilon_0&=1,\\ 
D&=0.05, & \bar{\gamma}&=1, &I_c&=0.001.
\end{align*}
To visualize the microscopic results in Figure \ref{fig:bsp1_mitill}, we consider $50$ particles and choose as initial conditions randomly distributed positions of the particles inside the domain $\Omega_0=[-1.0,~1.0]\times[-1.0,~1.0]$. Additionally, we set
\begin{align*}
\theta_0^i=0, \qquad v_0^i=(0,0)^\top, \qquad \omega_0^i=0,\qquad \text{for~} i=1,\ldots,50.
\end{align*}
For the microscopic histograms in Figure \ref{fig:bsp1_mit} we consider the same example with 1000 ellipsoidal particles and 128 Monte-Carlo realizations. For these results we show the corresponding density-histogram and a smoothed version of the histogram.
As for the hydrodynamic $q$-equation with mono-kinetic closure, we take a spatial grid size of $h=0.05$ and $k=\pi/30$ in the angle.
The corresponding initial conditions are
\begin{gather*}
q(0,r,\theta)=\begin{cases}
\frac{1}{4 k}, & \text{if,~} r\in\Omega_0 \text{~and~} \theta \in [0,~ k],\\
0, & \text{else},
\end{cases}\\
v(0,r,\theta)=(0,0)^\top,\qquad \omega(0,r,\theta)=0,\qquad \forall (r,\theta) \in\Omega\times[0,~2\pi].
\end{gather*}
For the hydrodynamic and diffusive  $\rho$-equations we choose the same grid size and 
the corresponding   initial conditions for the relevant quantities.

The results are displayed in Figure \ref{fig:bsp1_mit}. As seen in the figure, the hydrodynamic equations give a good approximation of the microscopic model. On the other hand, the results of the diffusive limit equation deviate strongly.

\begin{figure}[h!]
\includegraphics[width=0.49\textwidth]{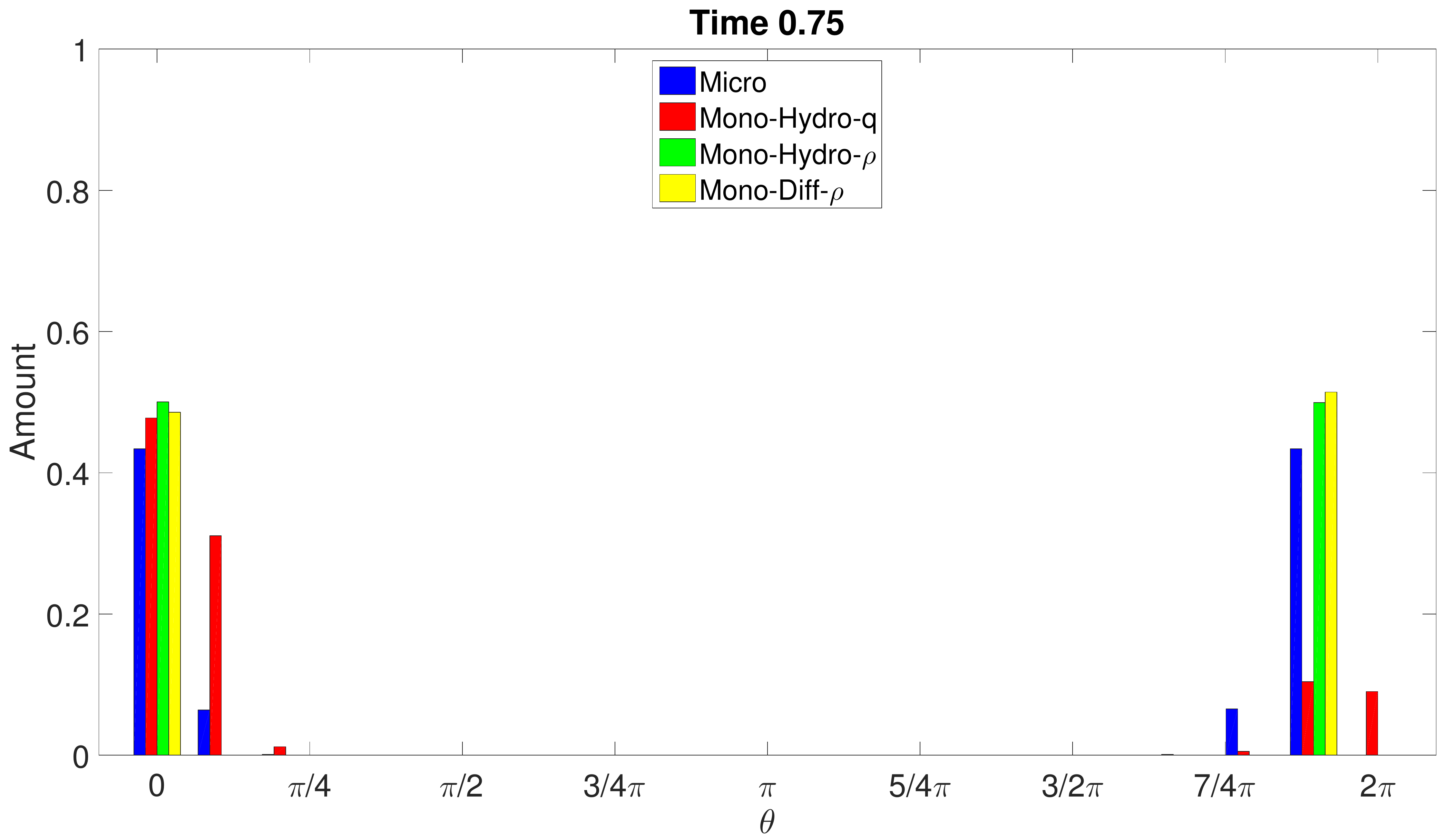}
\includegraphics[width=0.49\textwidth]{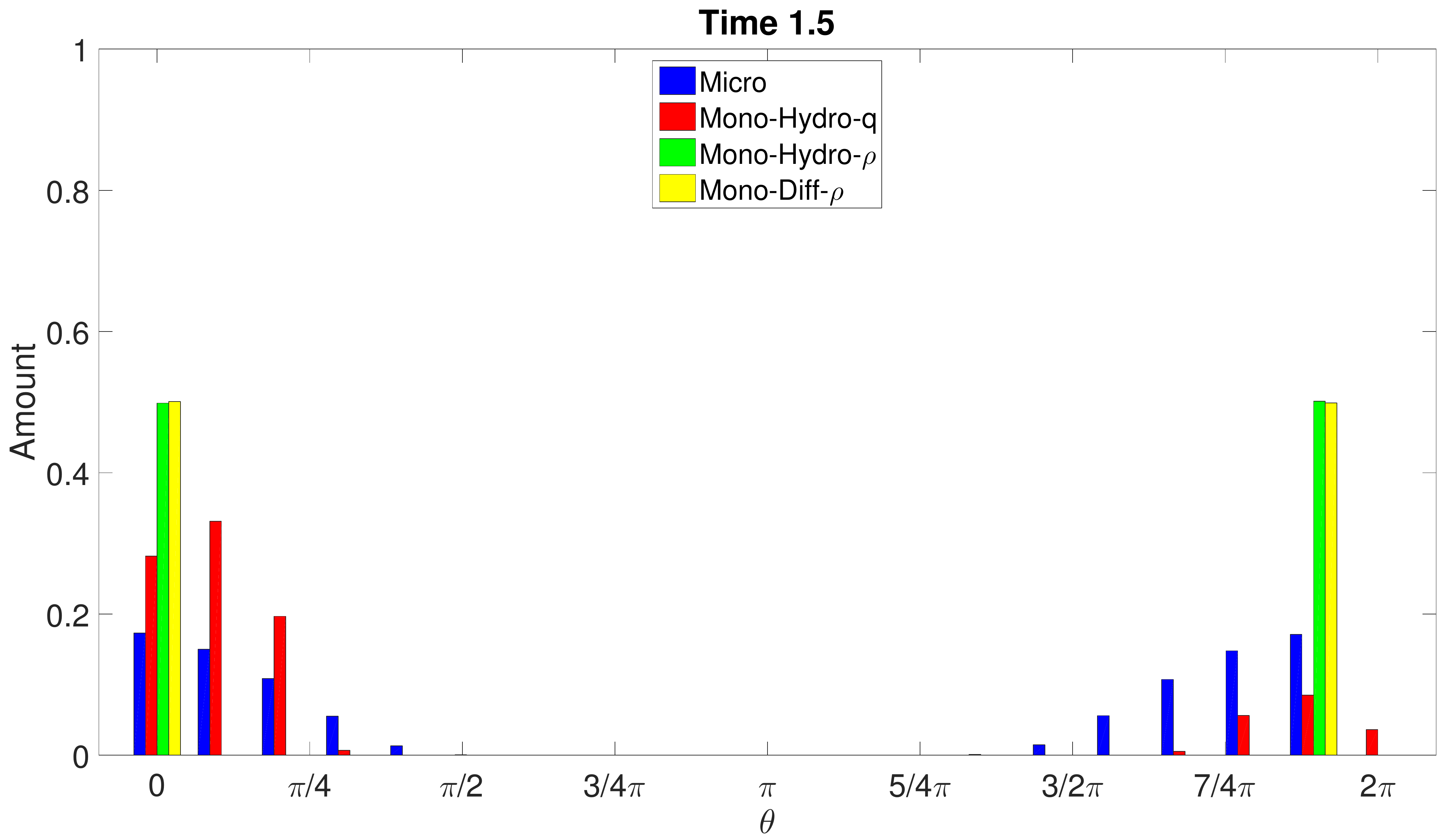}\\
\center \includegraphics[width=0.49\textwidth]{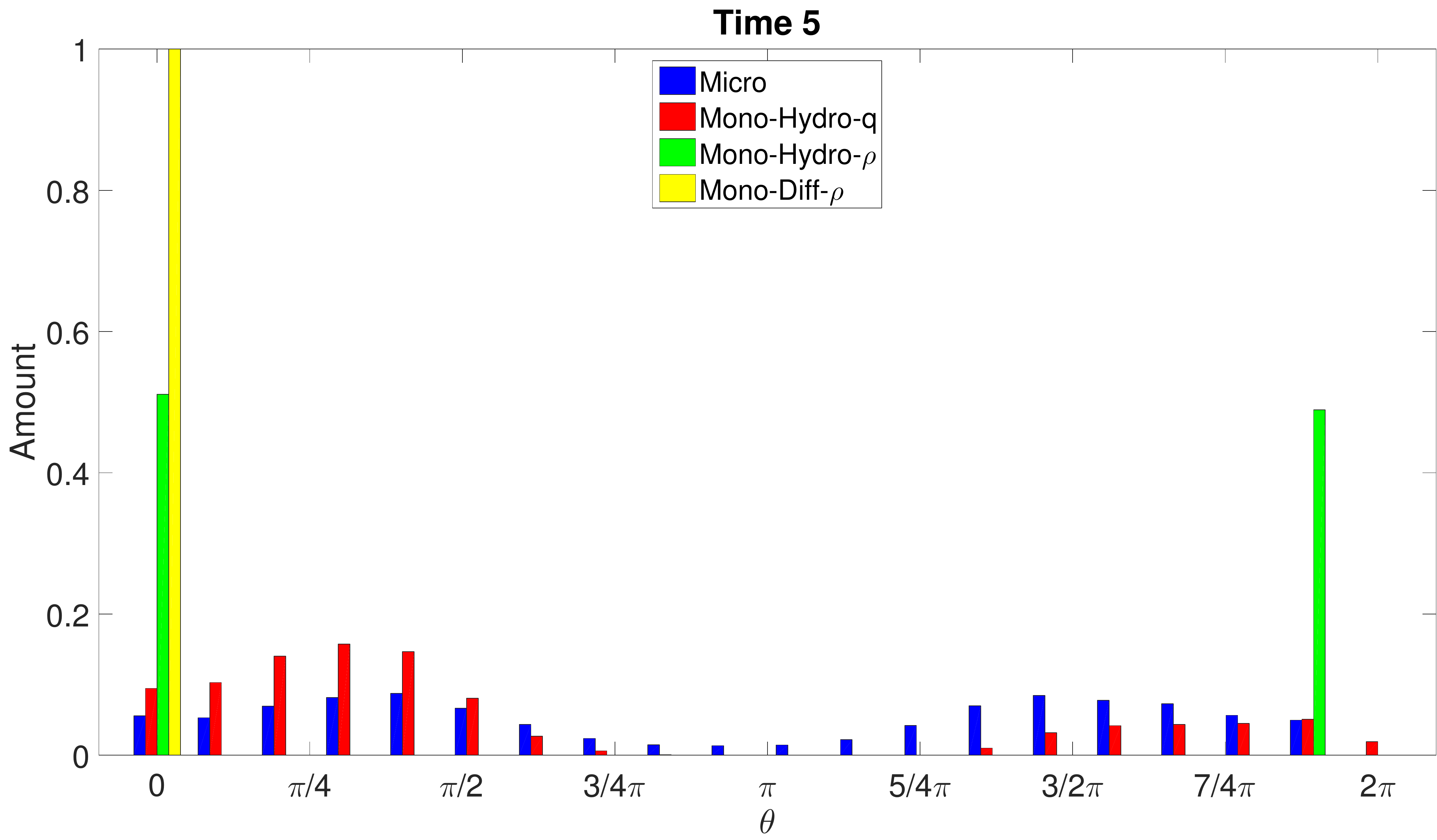}
\textbf{\caption{Angular distribution for the top-bottom example with interaction and without stochastic forces at time  $T\approx0.75$, $T\approx1.5$ and $T=5$ (right) for microscopic, mono-kinetic hydrodynamic $q$, $\rho$-equations and diffusive $\rho$-equation}
\label{fig:bsp1_winkel}}
\end{figure}

Considering the orientation angles, we can observe that, due to interaction, the angle changes towards a stationary equilibrium (see Figure \ref{fig:bsp1_winkel}). 
Again, the behavior of the microscopic model is well captured by the $q$-equation. The $\rho$-equations give a result, which is more concentrated around the angle $\theta=0$.

The results show that in the case without noise the mono-kinetic closure gives a good approximation of the microscopic model. 
Comparing the spatial distribution of the different macroscopic approximations with the microscopic model, we can observe that the mono-kinetic hydrodynamic equation for $q$ and $\rho$ give reasonable results. In these pictures the $\rho$ equation gives slightly less diffusive results. This is due to the larger amount of numerical diffusion in the solution of the $q$ equation. The diffusive approximation is not accurate for the parameters considered here. The comparison of the angular distributions shows that the  closures leading to the $\rho$ equations are not as appropriate as the ones leading to the $q$ equation.

Note that the computation time for the diffusive case is much larger than in the hydrodynamic case. In the hydrodynamic case we are able to pre-compute the weights for the interaction potential, which is not possible for the diffusive equations. Therefore, we will not consider the diffusive equations in the following examples.

\begin{remark}
Considering a situation with large parameters $A,B$ the mono-kinetic closure becomes inadequate and thus the derived equations
loose their validity. In such cases, other closures such as the Maxwellian closure, have to be used.
\end{remark}

\subsection{Rotational flow}
\begin{figure}[h!]
	\centering
	\includegraphics[width=0.5\textwidth]{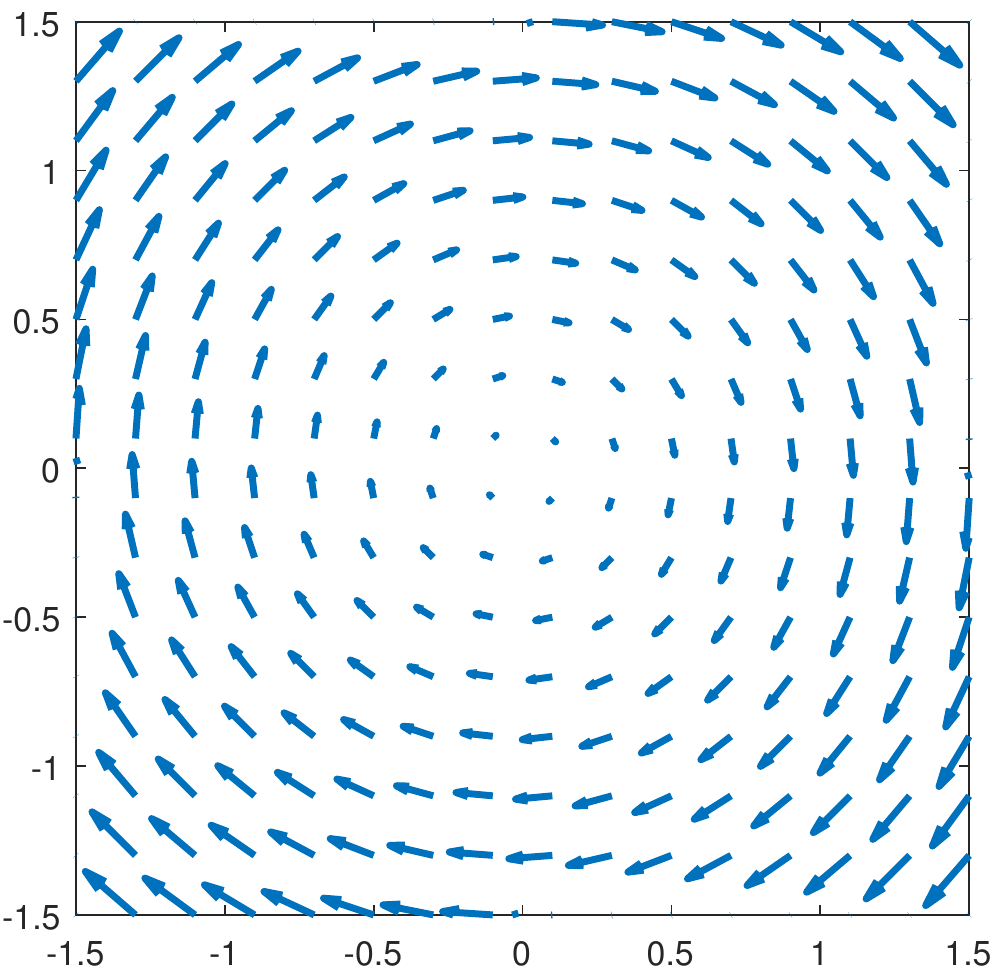}
	\caption{Rotational flow example: Velocity field of the surrounding fluid}
	\label{fig:rotational_vekfeld}
\end{figure}

\begin{figure}[h!]
  	\includegraphics[width=0.3\textwidth]{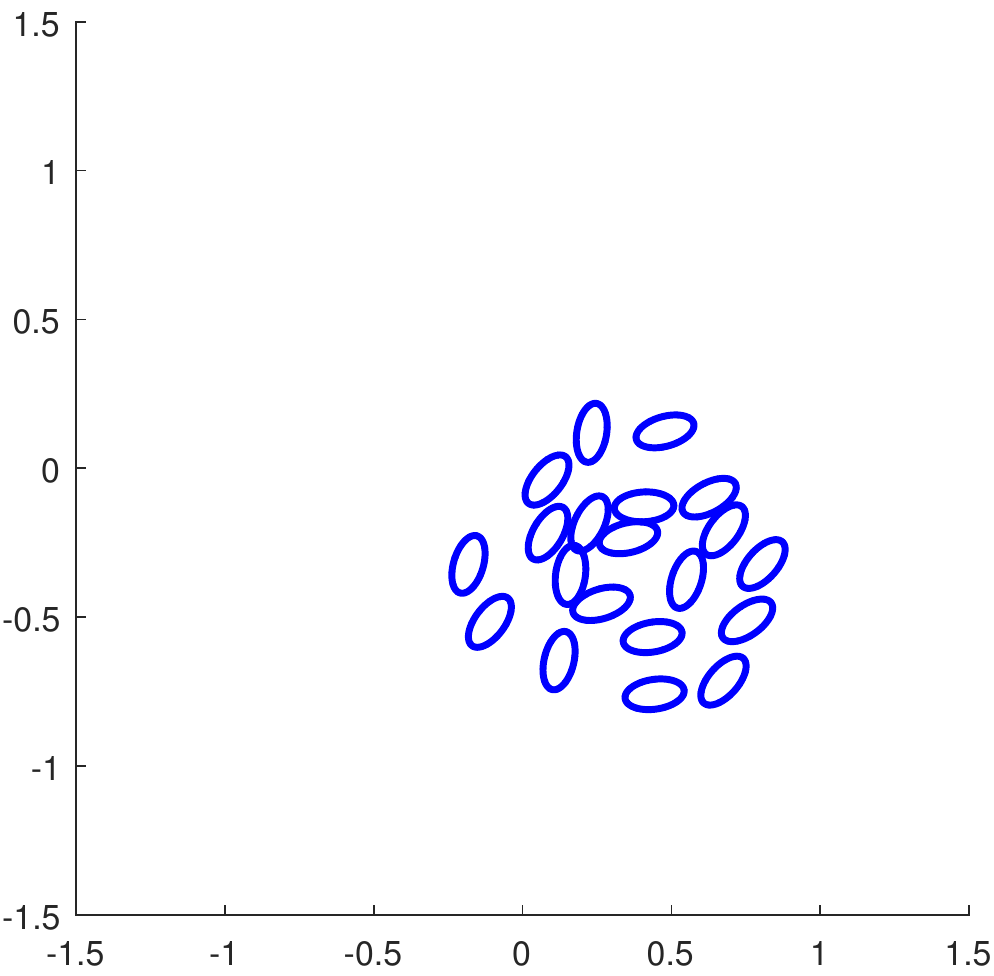}
  	\includegraphics[width=0.3\textwidth]{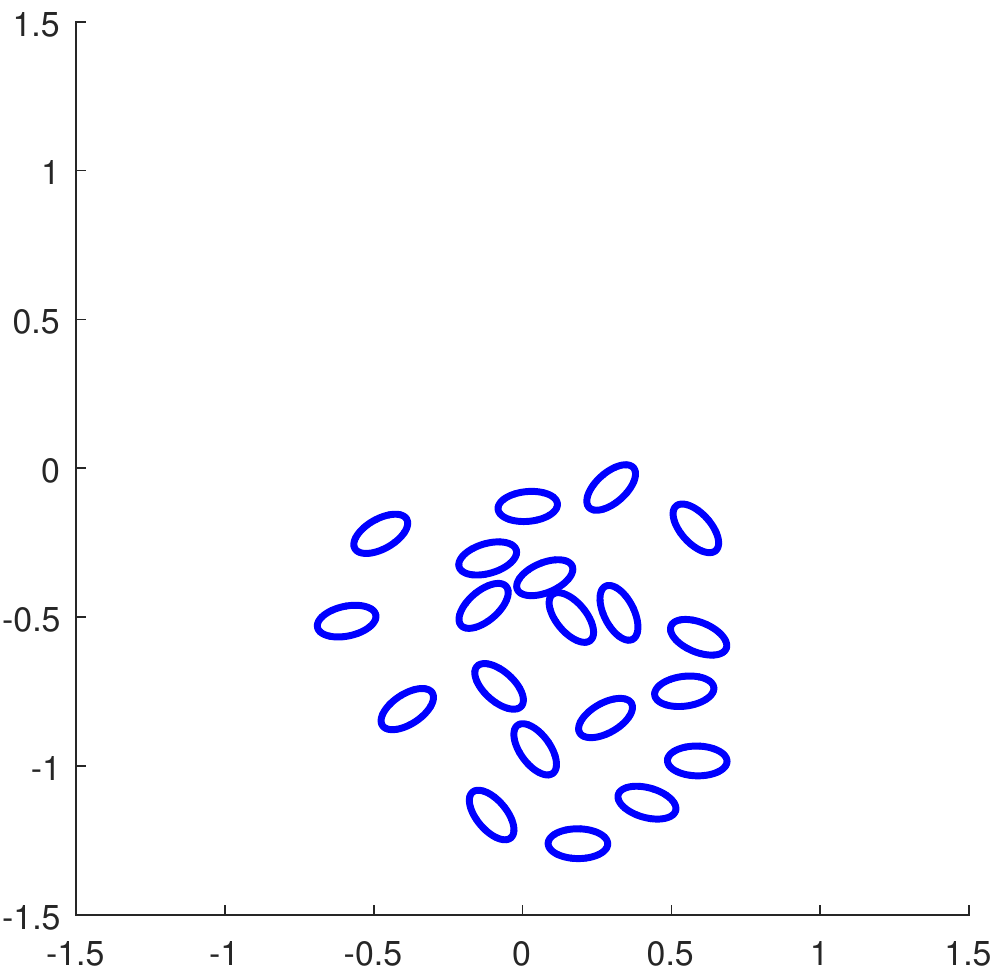}
  	\includegraphics[width=0.3\textwidth]{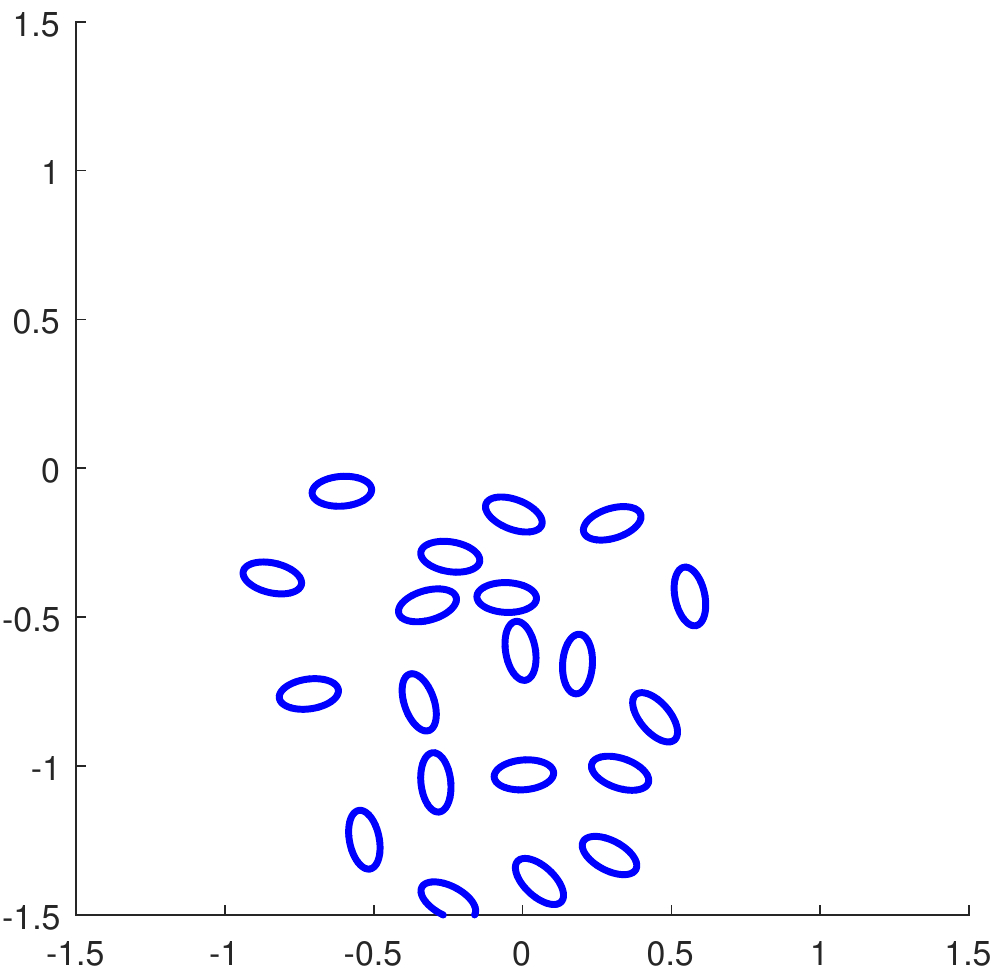}
  	\caption{Illustration of microscopic evolution for rotational flow example with interaction at time T=1.5 (left) and T=2.5 (middle) and  T=3 (right).}
  	\label{fig:bsp1_rotill}
  \end{figure}

\begin{figure}[h!]
	\includegraphics[width=0.3\textwidth]{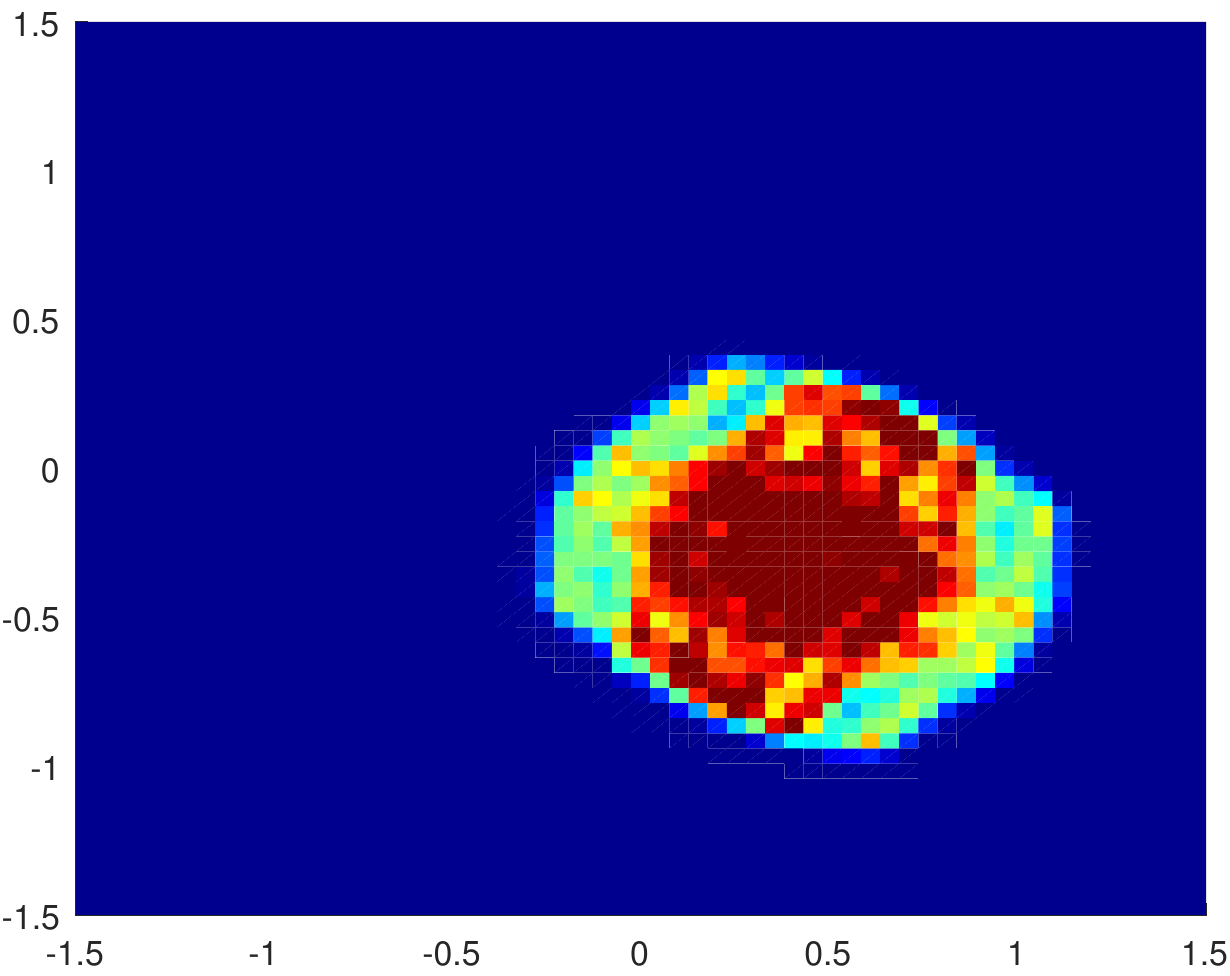}
	\includegraphics[width=0.3\textwidth]{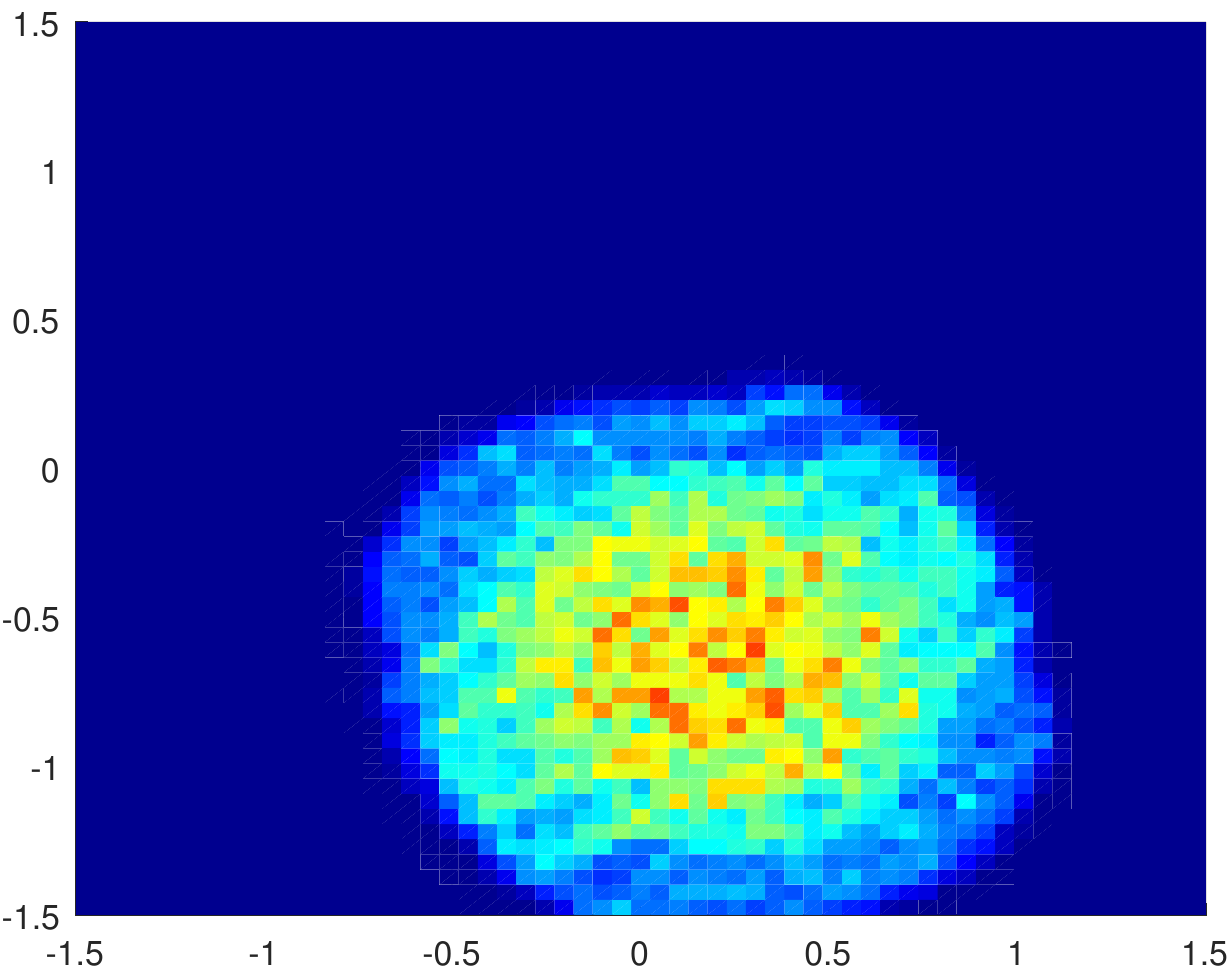}
	\includegraphics[width=0.3\textwidth]{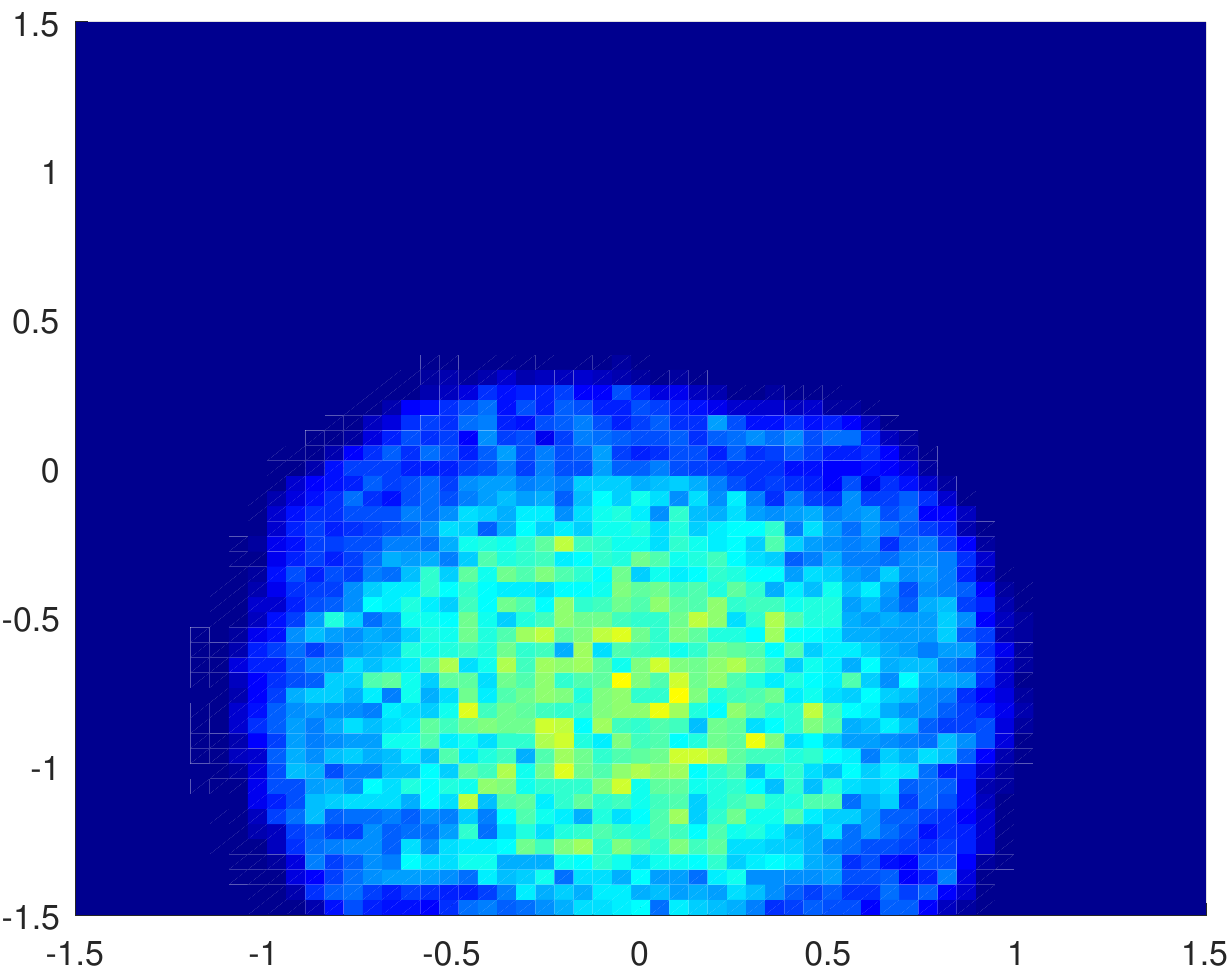}\\
	\includegraphics[width=0.3\textwidth]{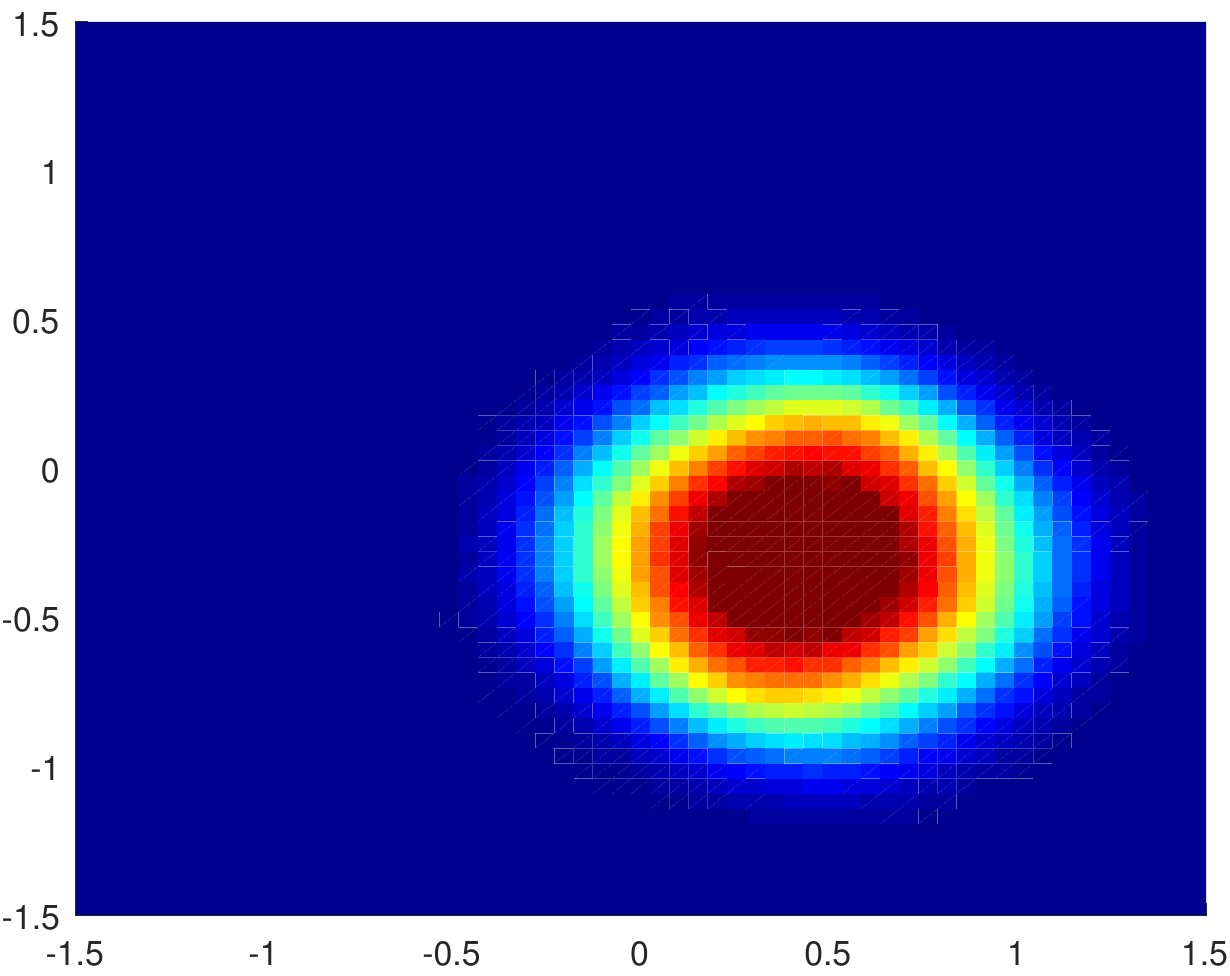}
	\includegraphics[width=0.3\textwidth]{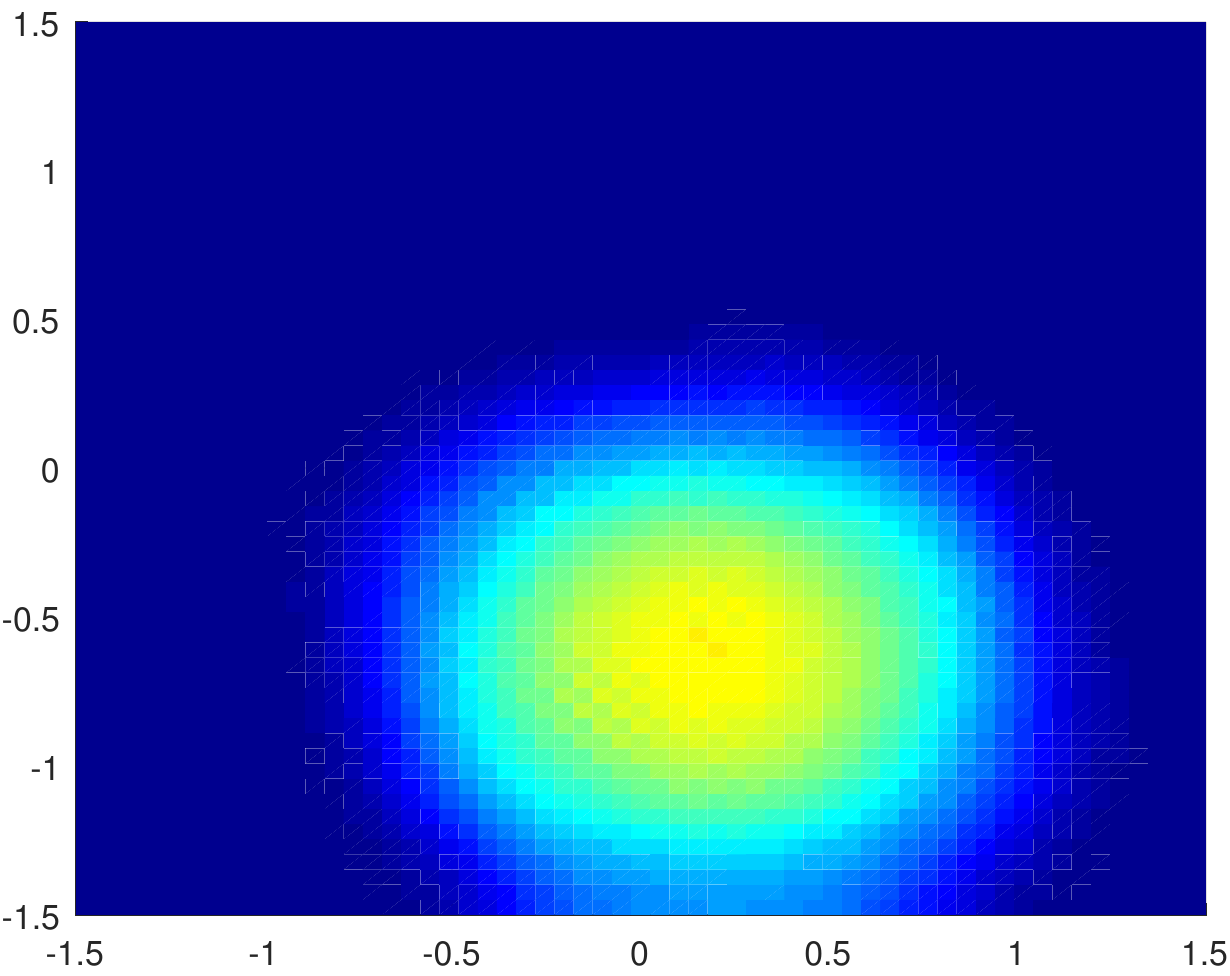}
	\includegraphics[width=0.3\textwidth]{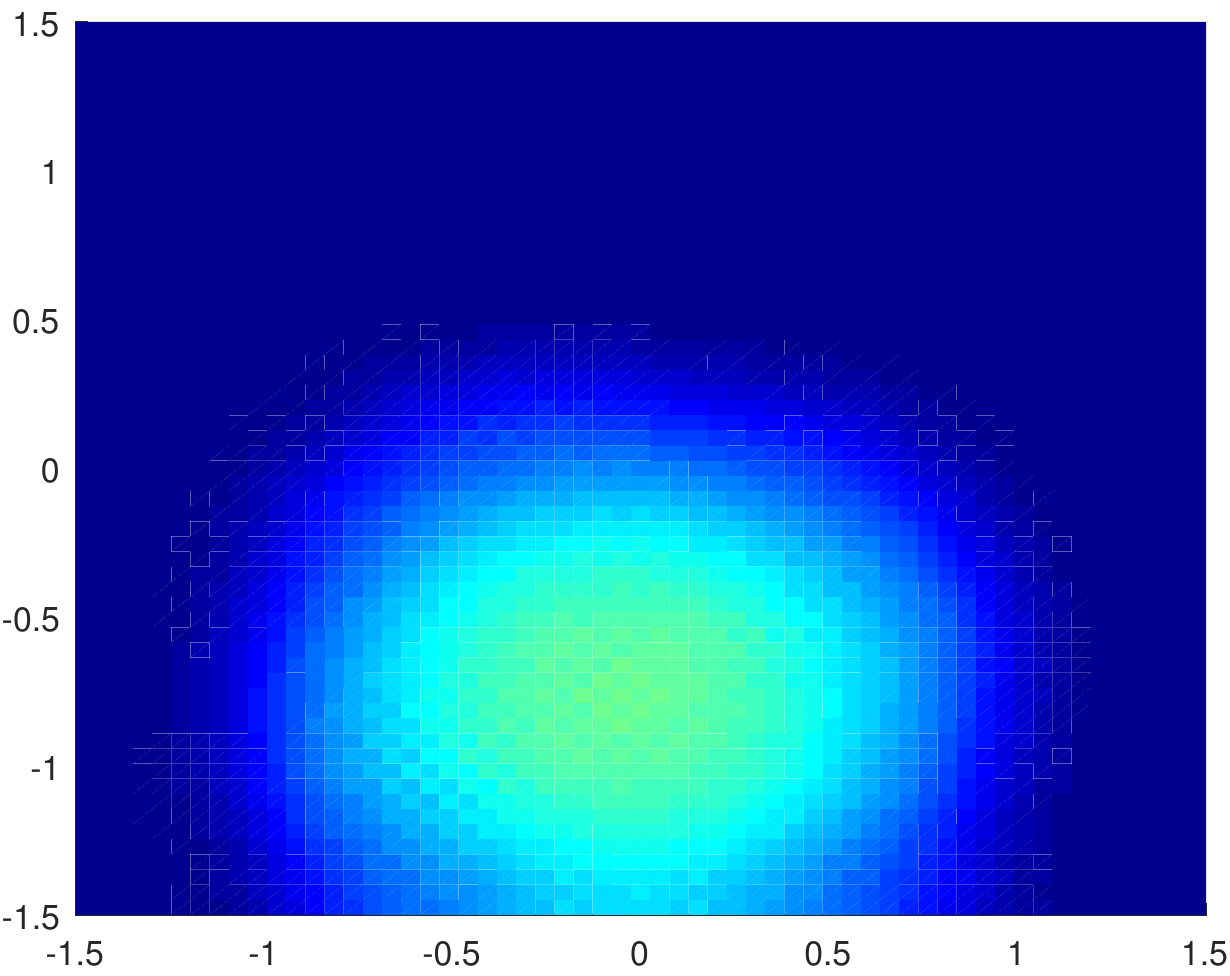}\\
	\includegraphics[width=0.3\textwidth]{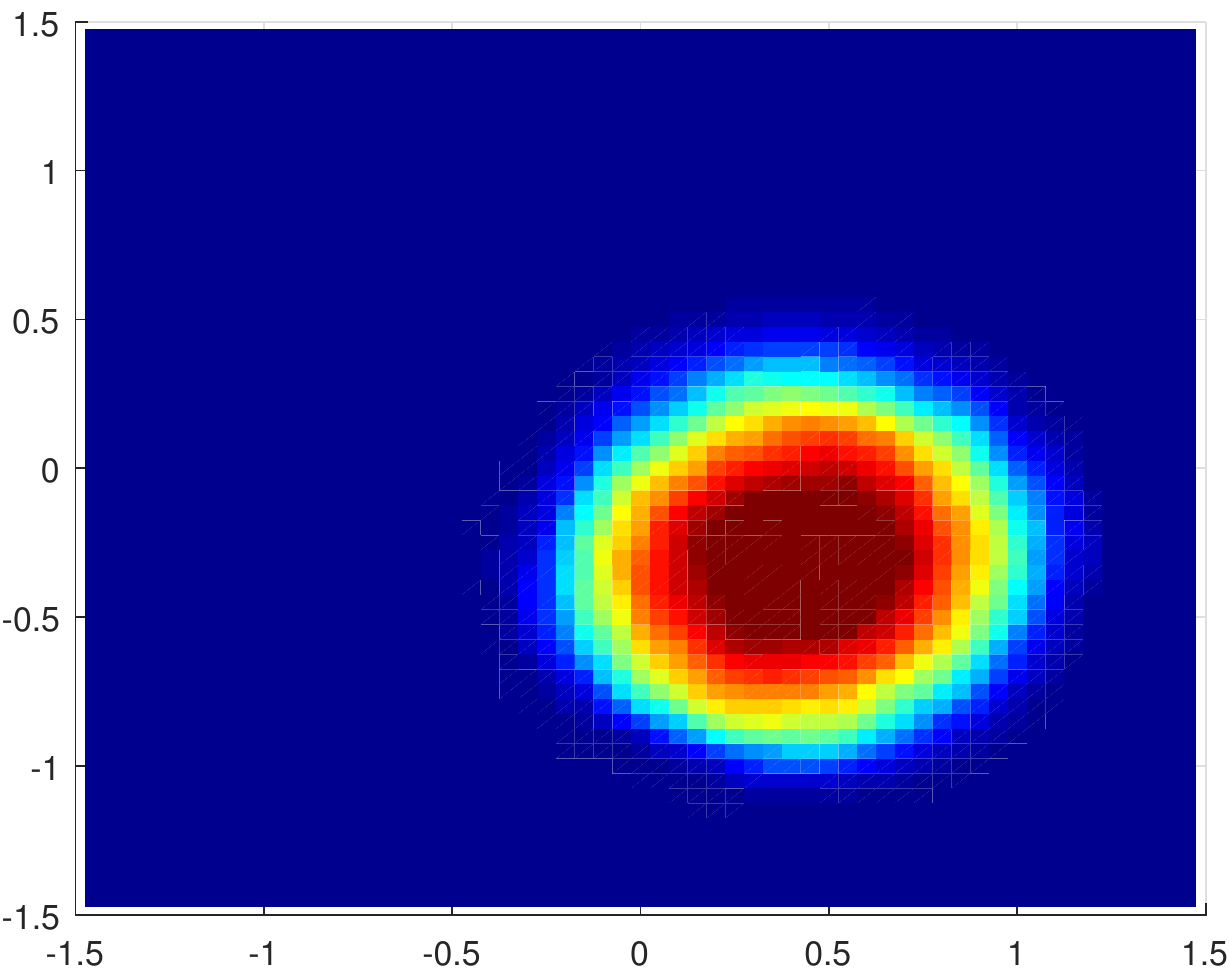}
	\includegraphics[width=0.3\textwidth]{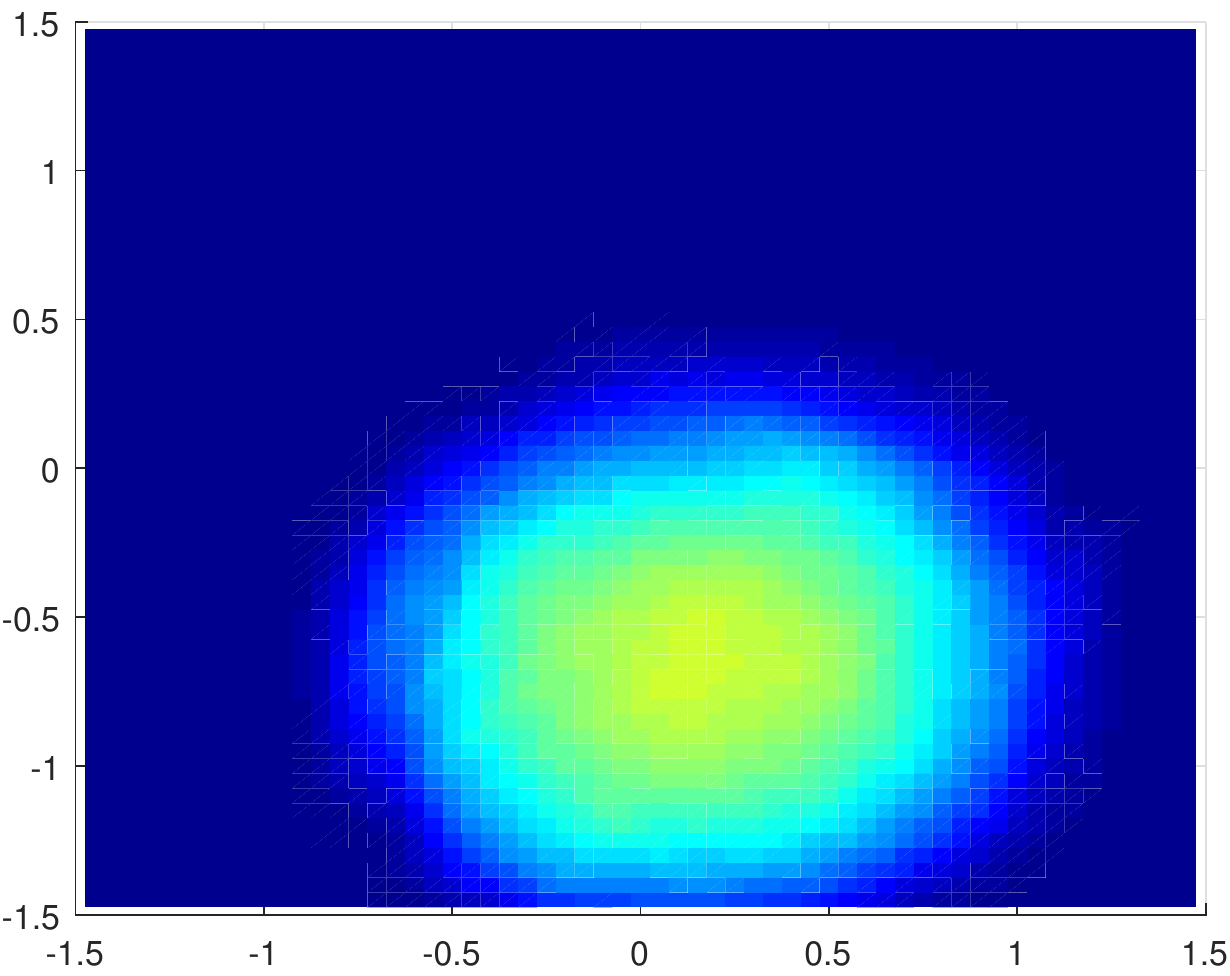}
	\includegraphics[width=0.3\textwidth]{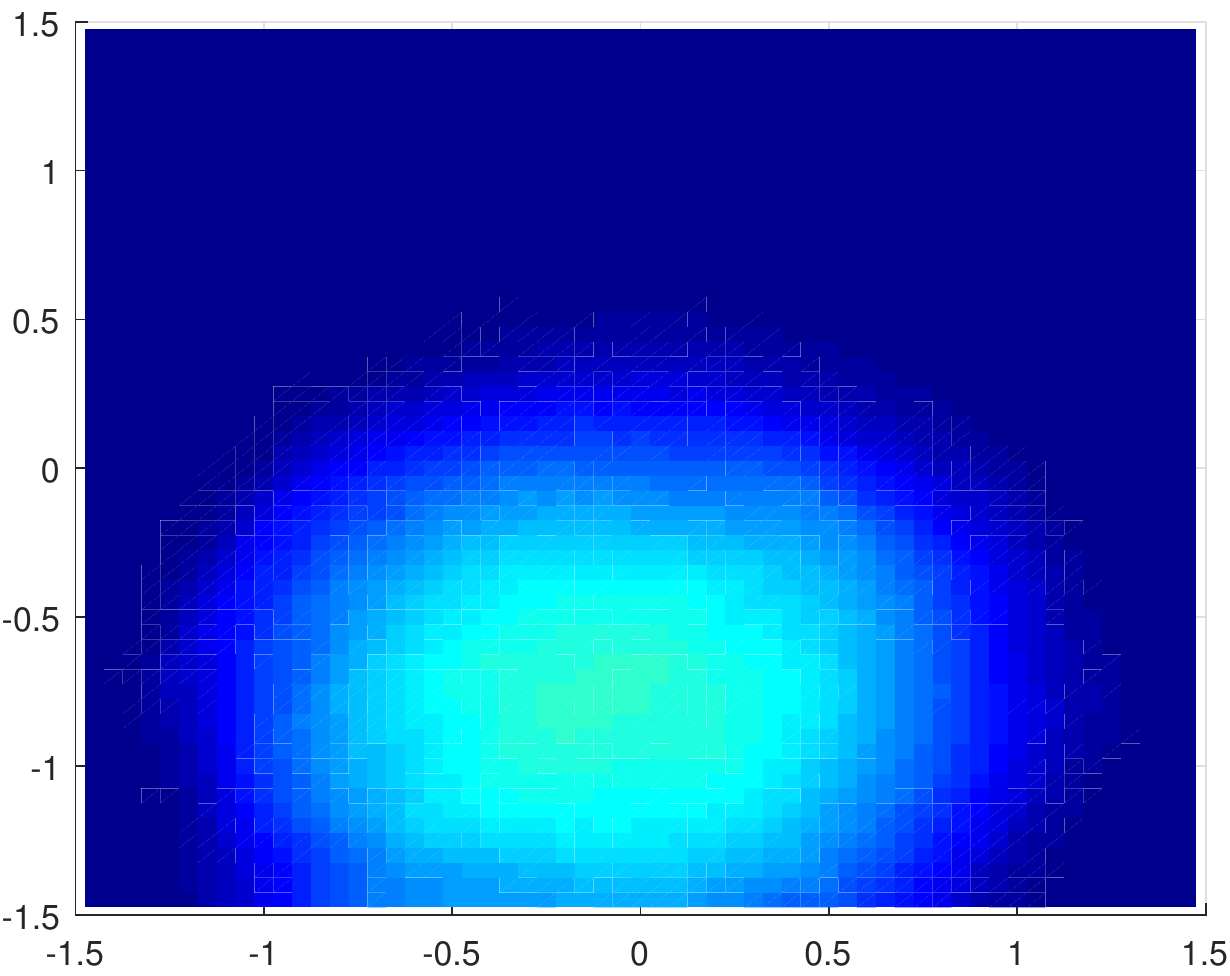}\\
	\includegraphics[width=0.3\textwidth]{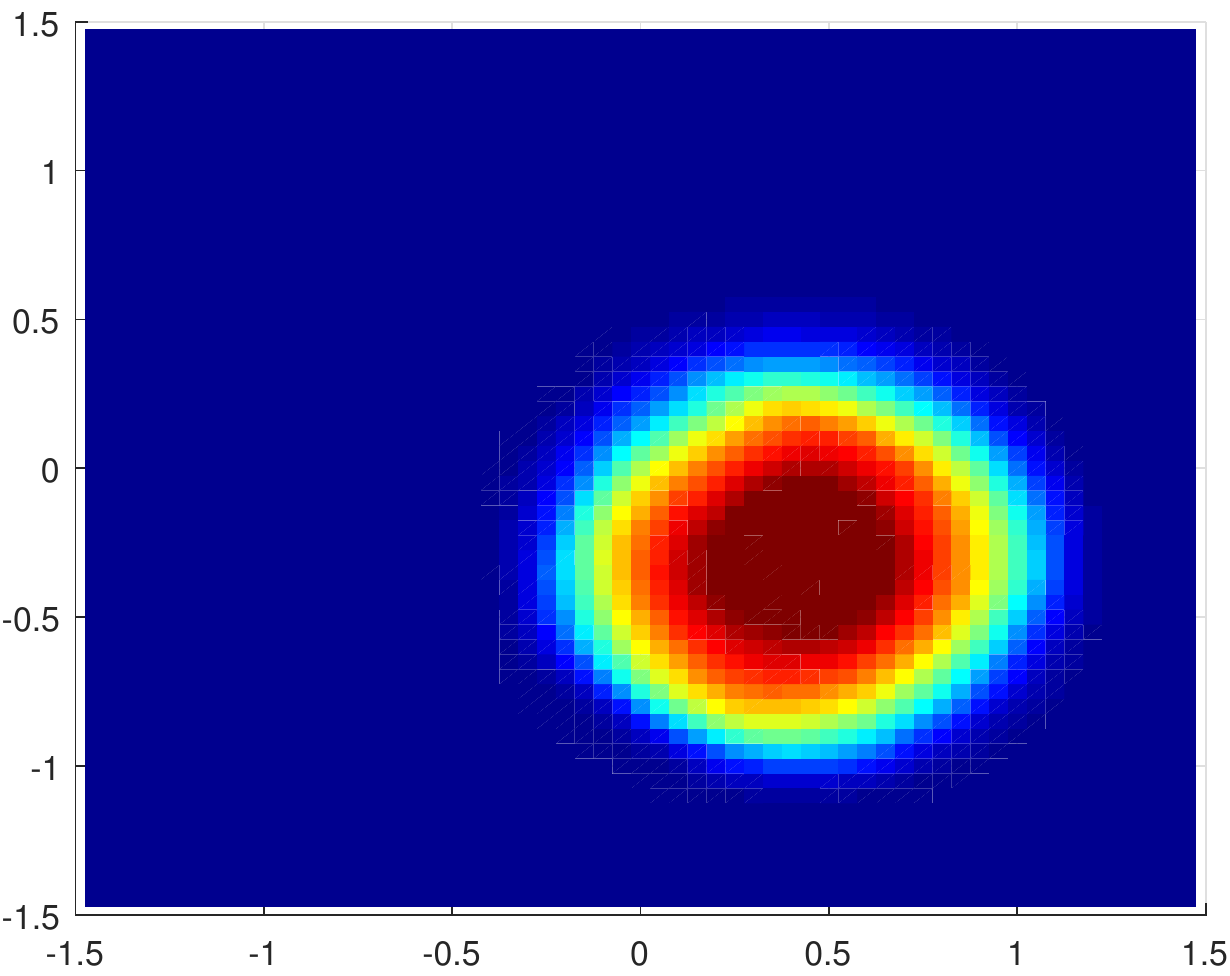}
	\includegraphics[width=0.3\textwidth]{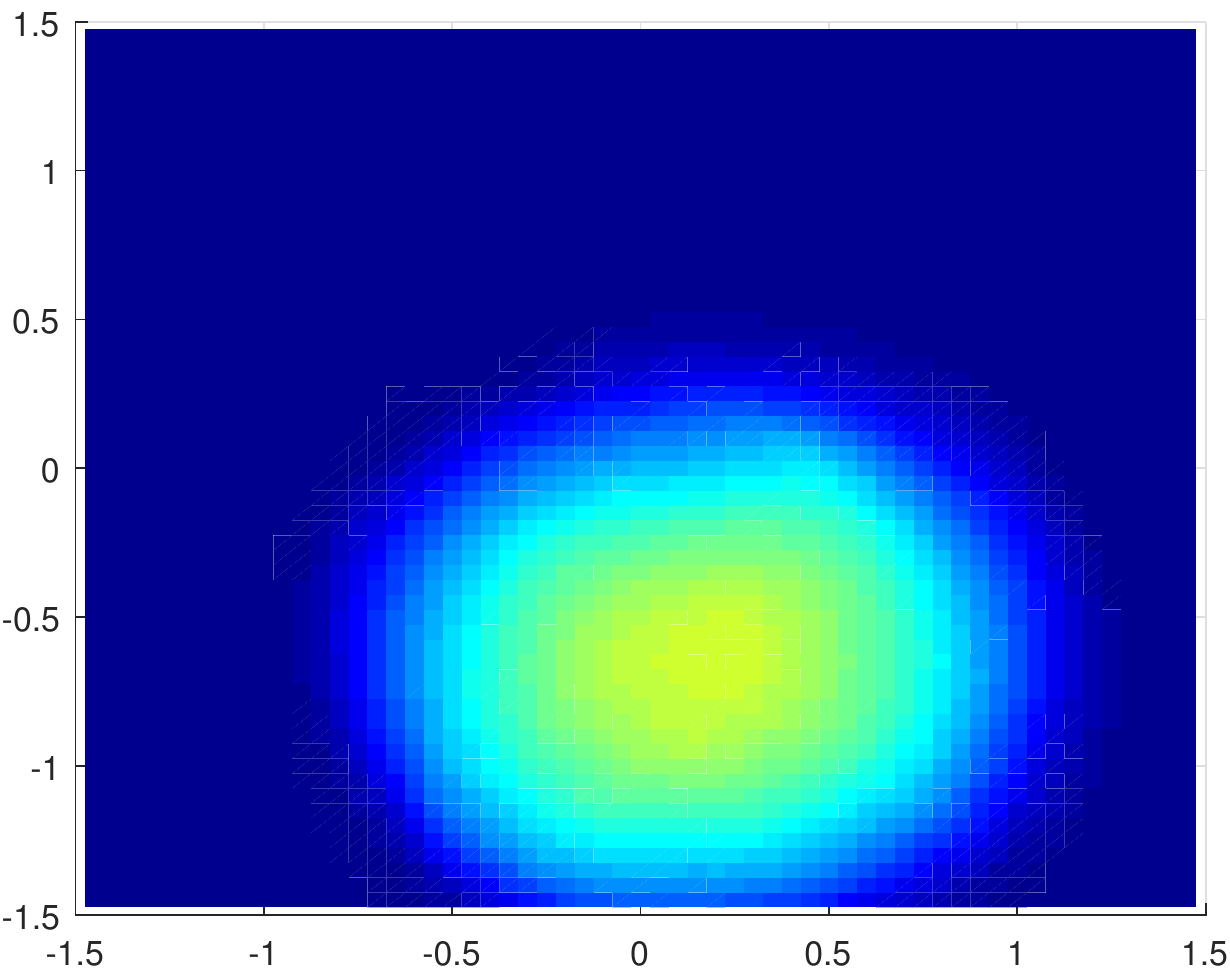}
	\includegraphics[width=0.3\textwidth]{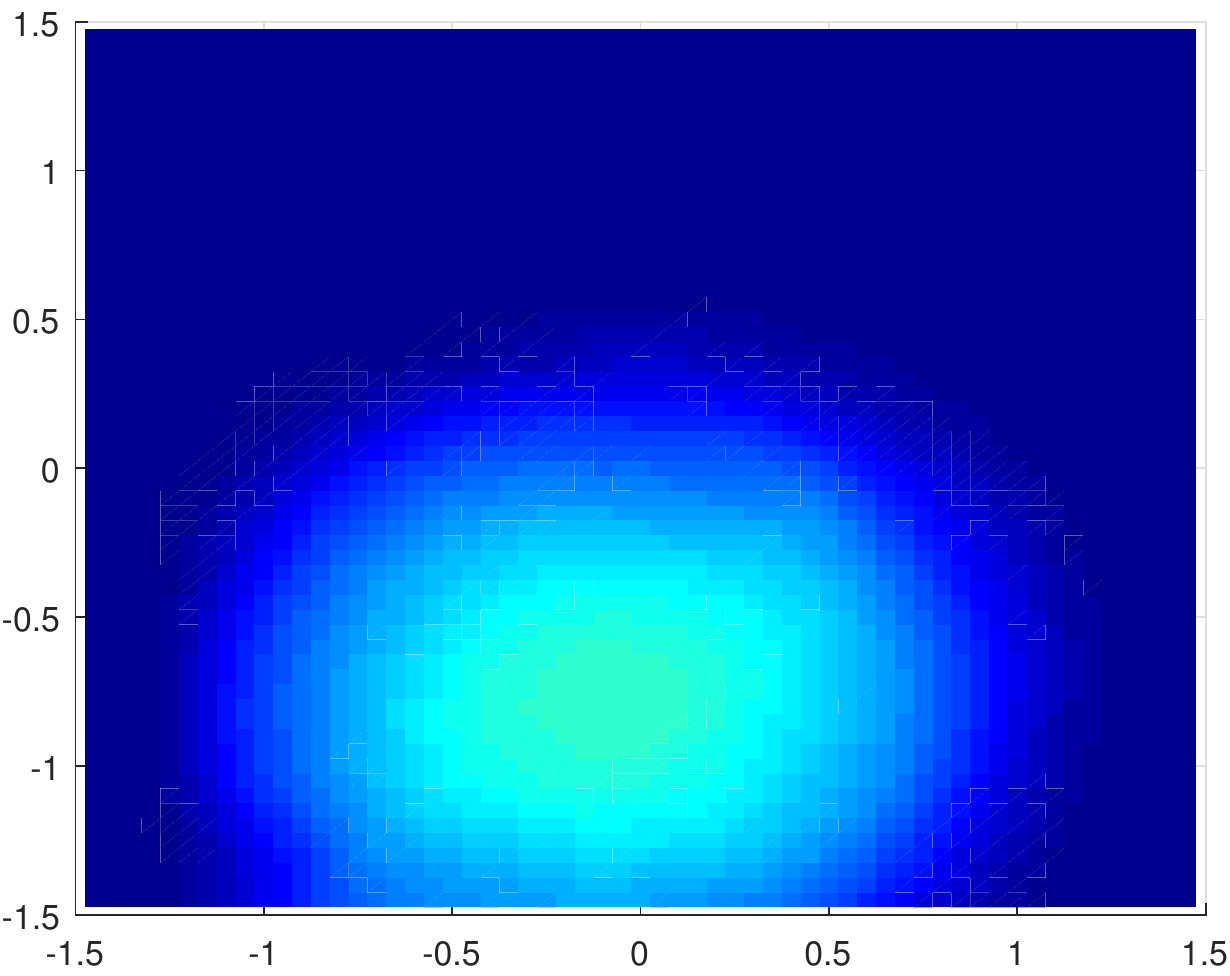}\\
		\centering\begin{tikzpicture}
			\begin{axis}[
			    hide axis,
			    scale only axis,
			    height=0pt,
			    width=0pt,
			    colormap/jet,
			    colorbar horizontal,
			    point meta min=0,
			    point meta max=1,
			    colorbar style={
			        width=10cm,
			        xtick={0,0.1,0.2,...,1}
			    }]
			    \addplot [draw=none] coordinates {(0,0)};
			\end{axis}
			\end{tikzpicture}
	\caption{Microscopic and hydrodynamic rotational flow example with interaction at time T=1.5 (left) and T=2.5 (middle) and T=3 (right), for (from top to bottom) microscopic histogram, smoothed microscopic histogram, mono-kinetic $q$ and mono-kinetic $\rho$-equations}
	\label{fig:bsp1_rot}
\end{figure}

In this example, we analyze the flow of particles in a rotational moving fluid. 
In this scenario we consider the microscopic system, the hydrodynamic $q$-equation \eqref{eq:macro_q_mono} and hydrodynamic $\rho$-equation 
\eqref{eq:macro_rho_mono} with mono-kinetic closure. The quadratic domain is $\Omega=[-1.5,~1.5]\times[-1.5,~1.5]$ and the parameters are chosen to be
\begin{align*}
L&=0.1, & \gamma&=1, & m&=1, & \epsilon_0&=1,\\ 
D&=0.05, & \bar{\gamma}&=1, & I_c&=0.001.
\end{align*}
As before, we neglect the stochastic forces, i.e. $A,B\equiv 0$. We assume that the velocity of the fluid is constant in time and not affected by the particles.
As velocity field of the fluid we choose $u(x,y)=(y,-x)^\top$ (see Figure \ref{fig:rotational_vekfeld}). The particles are expected to move clockwise around the origin. 
To obtain Figure \ref{fig:bsp1_rotill} we choose as initial conditions $20$ randomly distributed  particles inside the domain $\Omega_0=[0.2,~0.7]\times[-0.25,~0.25]$. 
Further values are
\[
 \theta_0^i=\pi/2,\qquad v_0^i=(0,0)^\top, \qquad \omega_0^i =0,\qquad \text{for~} i=1,\ldots,20.
\]
For Figure \ref{fig:bsp1_rot} we consider 1000 particles and 128 Monte-Carlo realizations for the microscopic simulation. Again we show the density histogram and smoothed version of the histogram.
For the hydrodynamic equation we use a spatial grid size of $h=0.02$ and $k=\pi/30$ in the angle coordinate.
The initial conditions are
\begin{gather*}
q(0,r,\theta)=\begin{cases}
\frac{4}{k}, & \text{if,~} r\in\Omega_0 \text{~and~} \theta \in [\pi/2,\pi/2 + k],\\
0, & \text{else},
\end{cases}\\
v(0,r,\theta)=(0,0)^\top,\qquad
\omega(0,r,\theta)=0,\qquad \forall (r,\theta) \in\Omega\times[0,~2\pi].
\end{gather*}
For the hydrodynamic $\rho$-equations we choose the same grid size and 
the corresponding initial conditions for the relevant quantities.

In the microscopic computation we can observe that the particles follow the fluid and rotate around the origin in a clockwise manner, as predicted.

At time $T=1.5$ the ellipsoids are almost completely separated, as shown in Figure \ref{fig:bsp1_rot}.
Due to the rotational flow and their inertia they further diverge during the computation.  
  
\begin{figure}[h!]
	\includegraphics[width=0.5\textwidth]{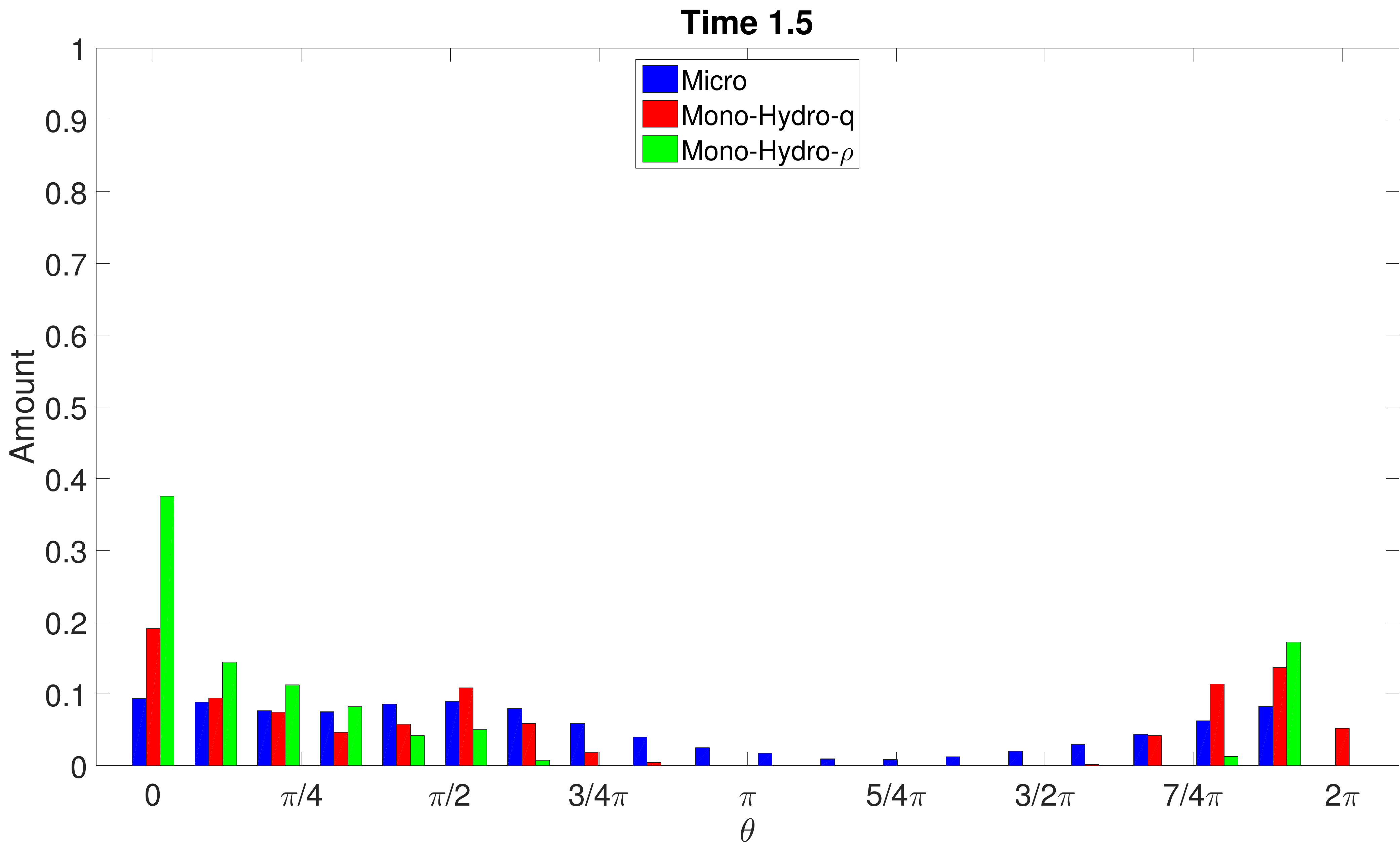}
	\includegraphics[width=0.5\textwidth]{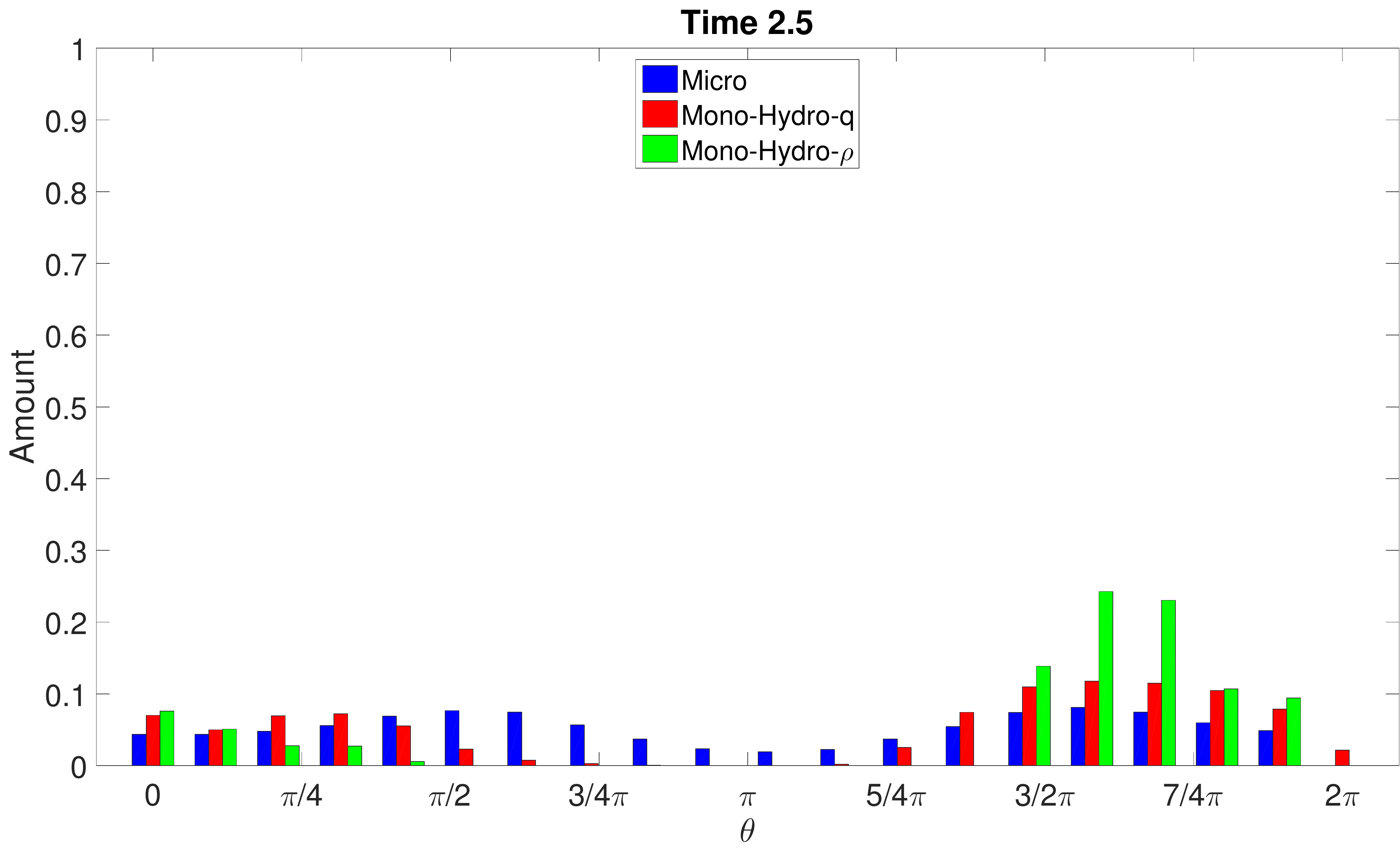}\\
	\center \includegraphics[width=0.5\textwidth]{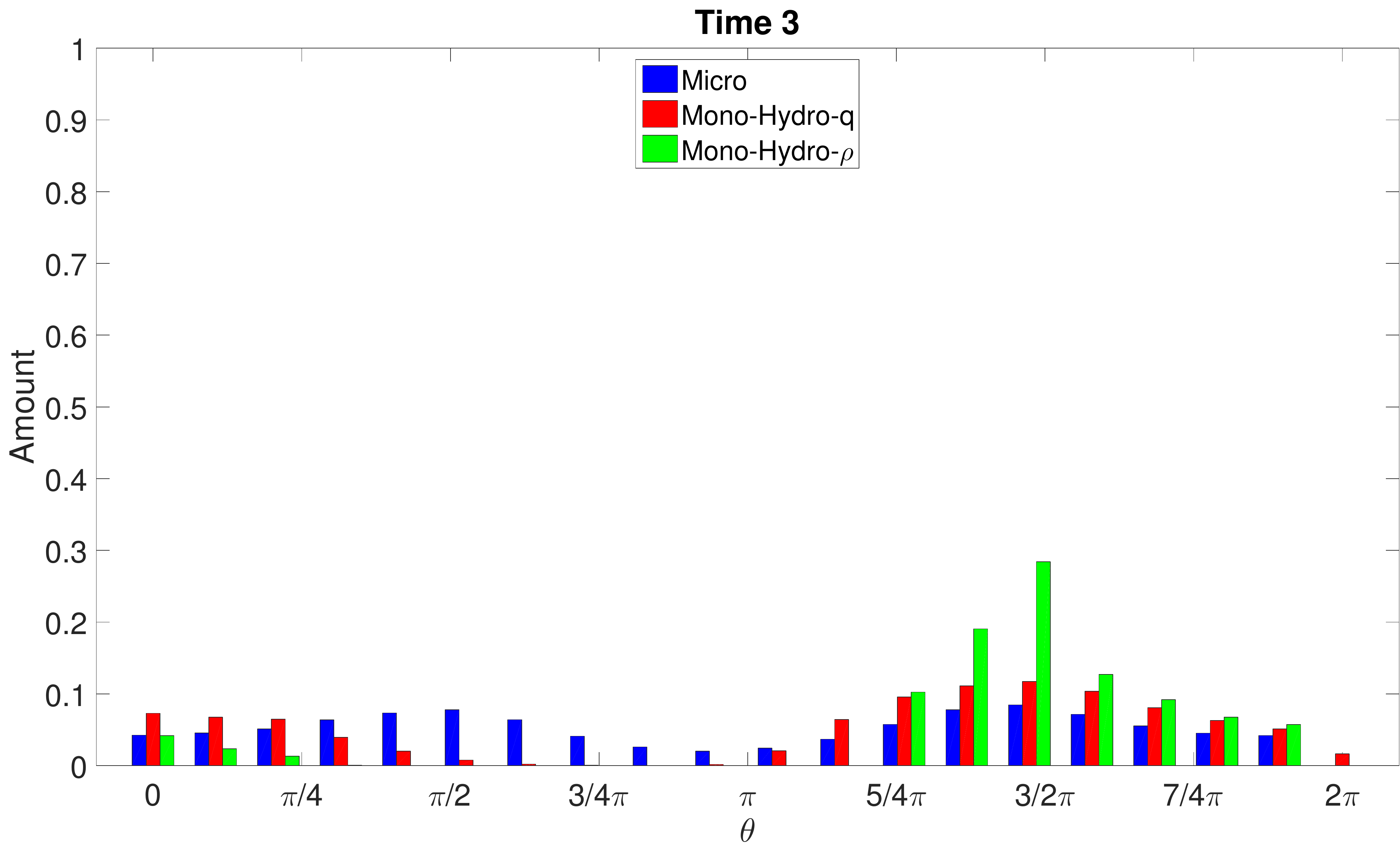}
	\caption{Angular distribution for rotational flow with interaction and without stochastic forces at time  $T\approx1.5$, $T\approx2.5$, $T=3$ for microscopic, mono-kinetic $q$  and  mono-kinetic $\rho$-equations}
	\label{fig:bsp3_winkel}
\end{figure}

 The angular distribution at the different time steps is displayed in Figure \ref{fig:bsp3_winkel}. Here, one also observes a good agreement of the distribution functions.

Comparing the microscopic results with the values of the hydrodynamic approximations, all models show a similar behavior.
In particular, the orientation of the ellipsoids show a very good agreement for all models.
In this case we cannot observe any significant differences between both hydrodynamic models
for the spatial distribution. The angular distribution is again better captured by the $q$ equation.

We note that for the rotational flow problem the results of  the equations without interaction differ strongly from the one with interaction. 
Without repulsion, the orientation of the ellipsoids is simply given by an angular distribution, which is strongly peaked in the direction of the flow.

% % % % %
% Cavitylam
% % % % %
\subsection{Driven cavity}
In this last test case we want to analyze the behavior of the particles in a more complex flow.
We compare the microscopic solutions with those of the hydrodynamic $\rho$-equations with mono-kinetic closure \eqref{eq:macro_rho_mono}.

\begin{figure}[h!]
	\centering
	\includegraphics[width=0.5\textwidth]{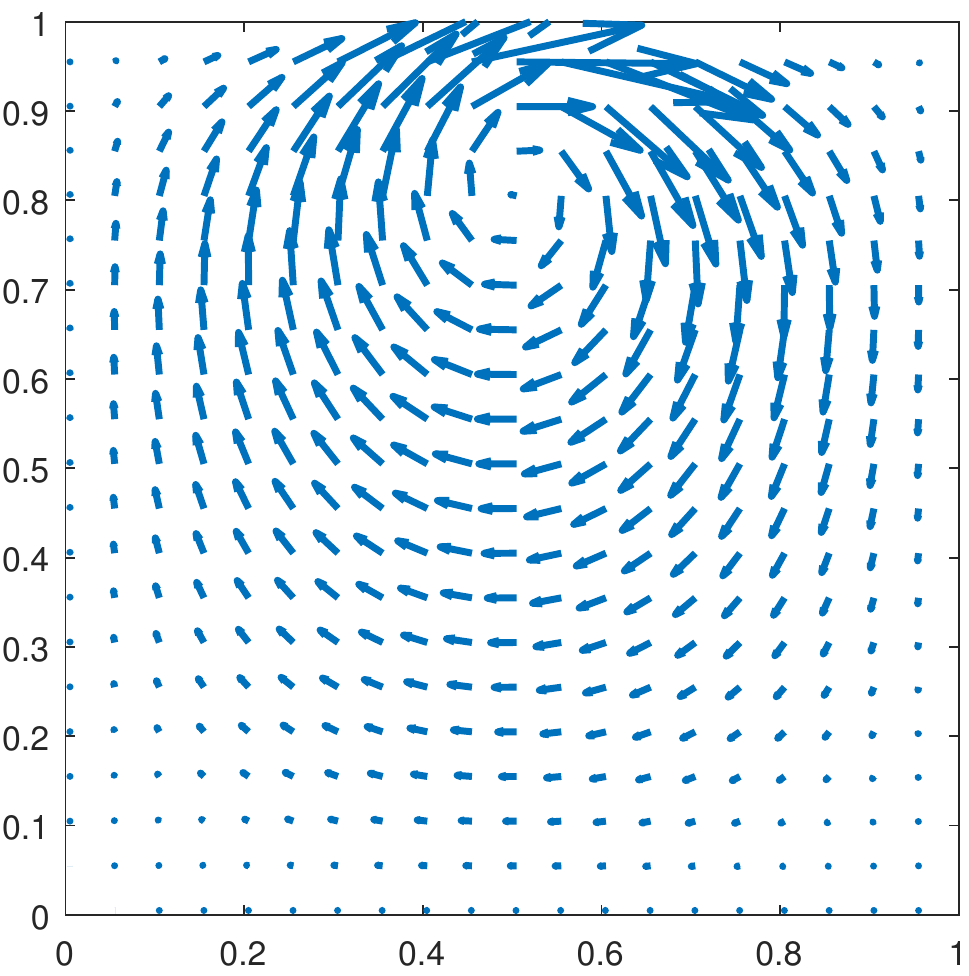}
	\caption{Vector field of the surrounding fluid.}
	\label{fig:vecfeld_cavitylam_gauss}
\end{figure}

In the quadratic domain $\Omega=[0,~1]\times[0,~1]$, the velocity field is given by the numerical solution of the incompressible Navier-Stokes equation \cite{foam}.
As boundary conditions for the fluid equations  we use no slip conditions at the left, right and bottom boundaries 
and the following velocity profile at the top boundary
\begin{align*}
u_1=\exp\left(-\frac{1}{2}\frac{(x-0.5)^2}{0.01}\right),\qquad u_2=0.
\end{align*}
The kinematic viscosity $\nu=0.1$ leads to a maximal Reynolds number of $Re=10$.
The resulting velocity field is sketched in Figure \ref{fig:vecfeld_cavitylam_gauss}.

To include wall boundaries for the ellipsoidal particles in the microscopic case, we insert ghost particles with distance $l/2$ onto the boundaries, where their orientations lie parallel to the respective wall. The interaction potential is the same as before \eqref{eq:interaction}. For the interaction of the boundary particles with the inner particles, we increase the value of $\epsilon_0$ by the factor $10$.
For the PDE models a similiar procedure using ghost cells is used. We extend the domain such that all interaction kernels inside the domain 
can be fully discretized. The density at the cells outside of the computational domain is set to the maximal density of $1/h^2$.
%Their orientation is parallel to the wall. 
The outside cells will only be considered for computing the interaction potential. 
As in the microscopic case, the strength of the boundary potentials will be enlarged by increasing $\epsilon_0$ by the factor $10$. 
In the hyperbolic part of the  splitting algorithm we use reflective boundary conditions.

Furthermore we choose the following parameters
\begin{align*}
L&=0.05, & \gamma&=10, & m&=1, & \epsilon_0&=1,\\ 
D&=0.025, & \bar{\gamma}&=10, & I_c&=0.001,
\end{align*}
and neglect the stochastic forces $A,B\equiv0$ and exterior potentials $V_1,V_2\equiv0$.

For the microscopic case, we choose randomly distributed positions of $1000$ particles inside $\Omega_0=[0.4,~0.6]\times[0.4,~0.6]$ and
\[
\theta_0^i=0, \qquad v_0^i=(0,0)^\top, \qquad \omega_0^i=0,\qquad \text{for~} i=1,\ldots,1000,
\]
for the initial condition, while for the hydrodynamic equation the spatial grid size is $h=0.005$ with initial conditions
\begin{gather*}
\rho(0,r)=\begin{cases}
25, &\text{if,~} r\in\Omega_0,\\
0, &\text{else},
\end{cases}\\
v(0,r)=(0,0)^\top,\qquad \omega(0,r)=0,\qquad \phi(0,r)=0,\qquad \forall r \in\Omega.
\end{gather*}

\begin{figure}[h!]
	\includegraphics[width=0.3\textwidth]{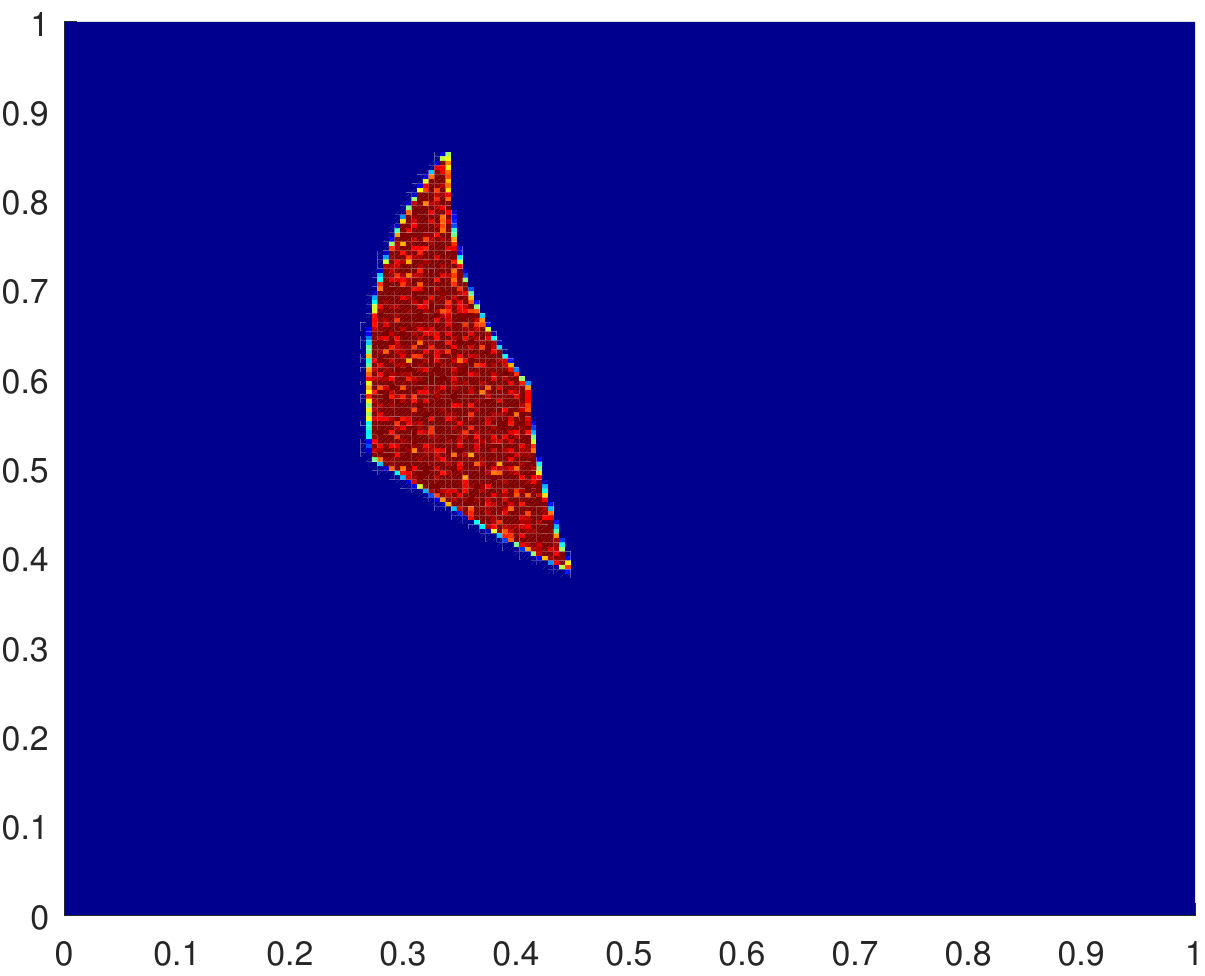}
	\includegraphics[width=0.3\textwidth]{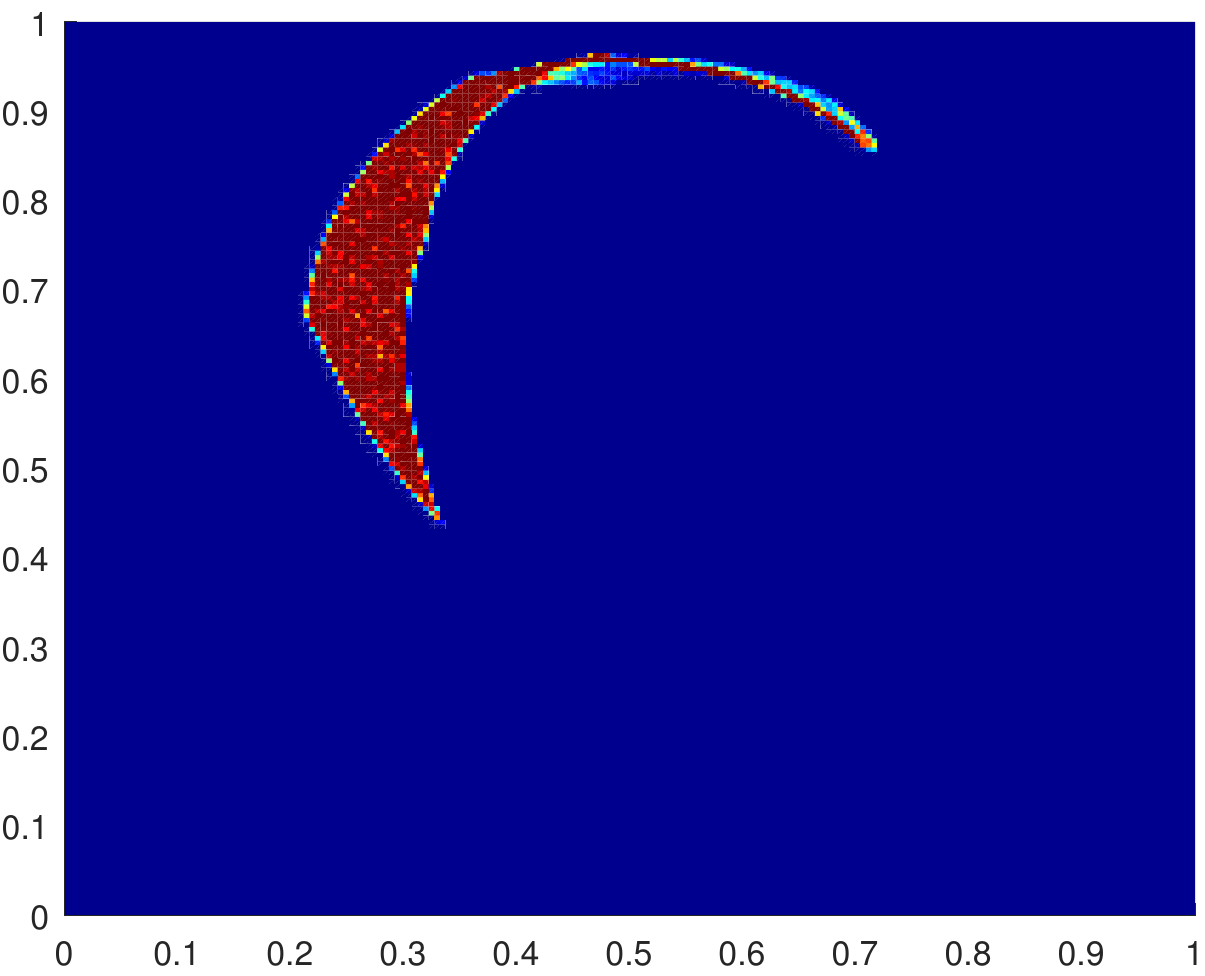}
	\includegraphics[width=0.3\textwidth]{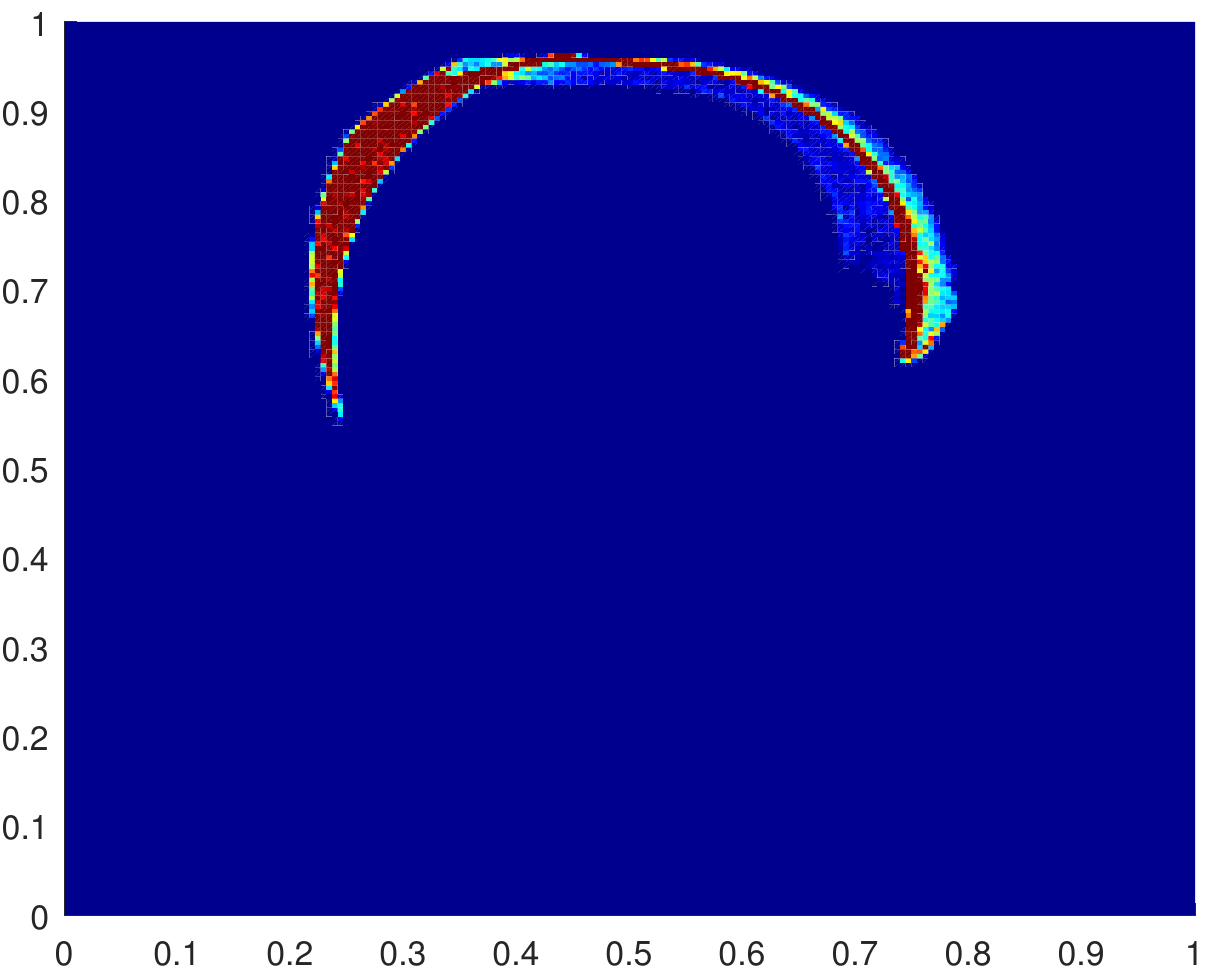}\\
	\includegraphics[width=0.3\textwidth]{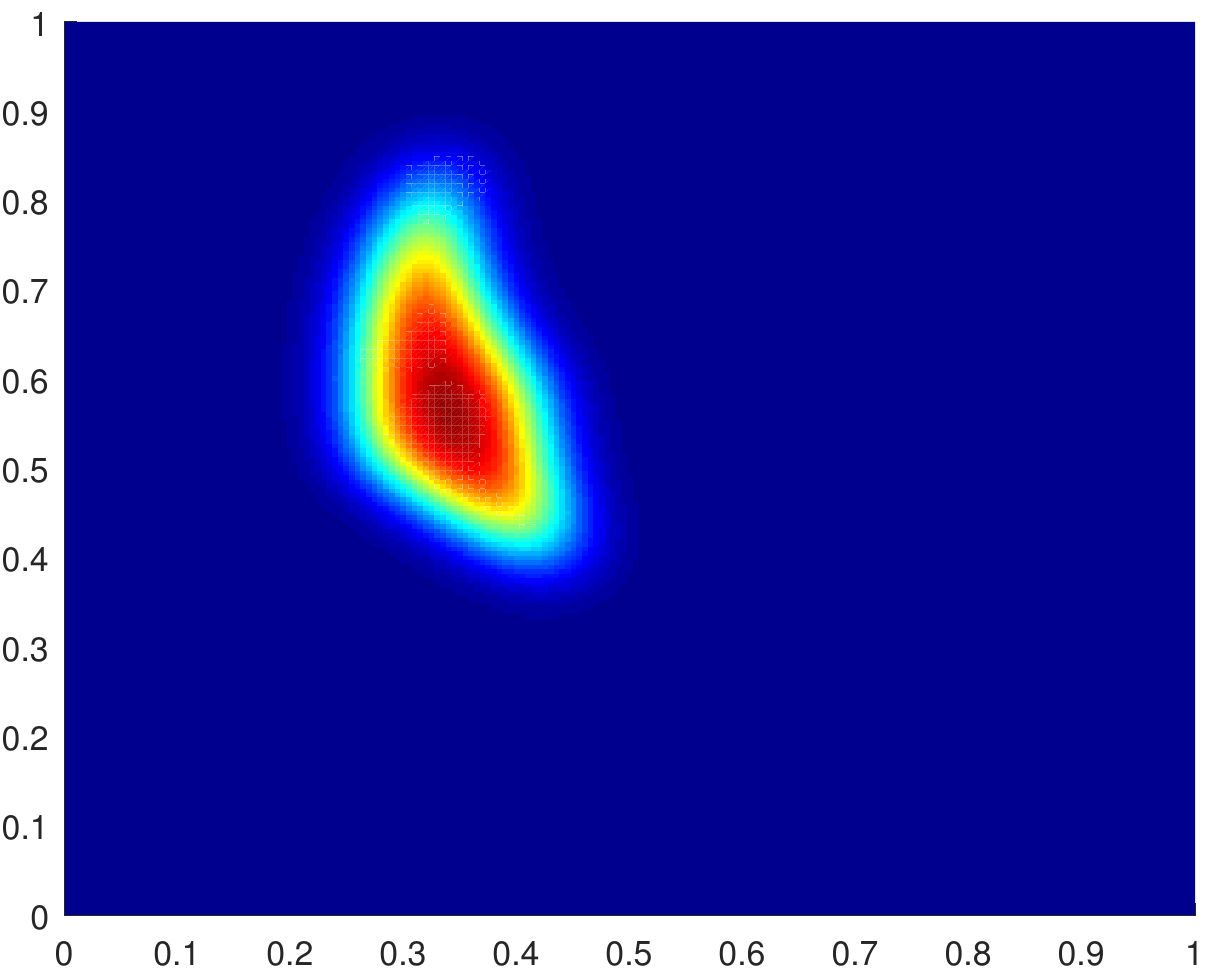}
	\includegraphics[width=0.3\textwidth]{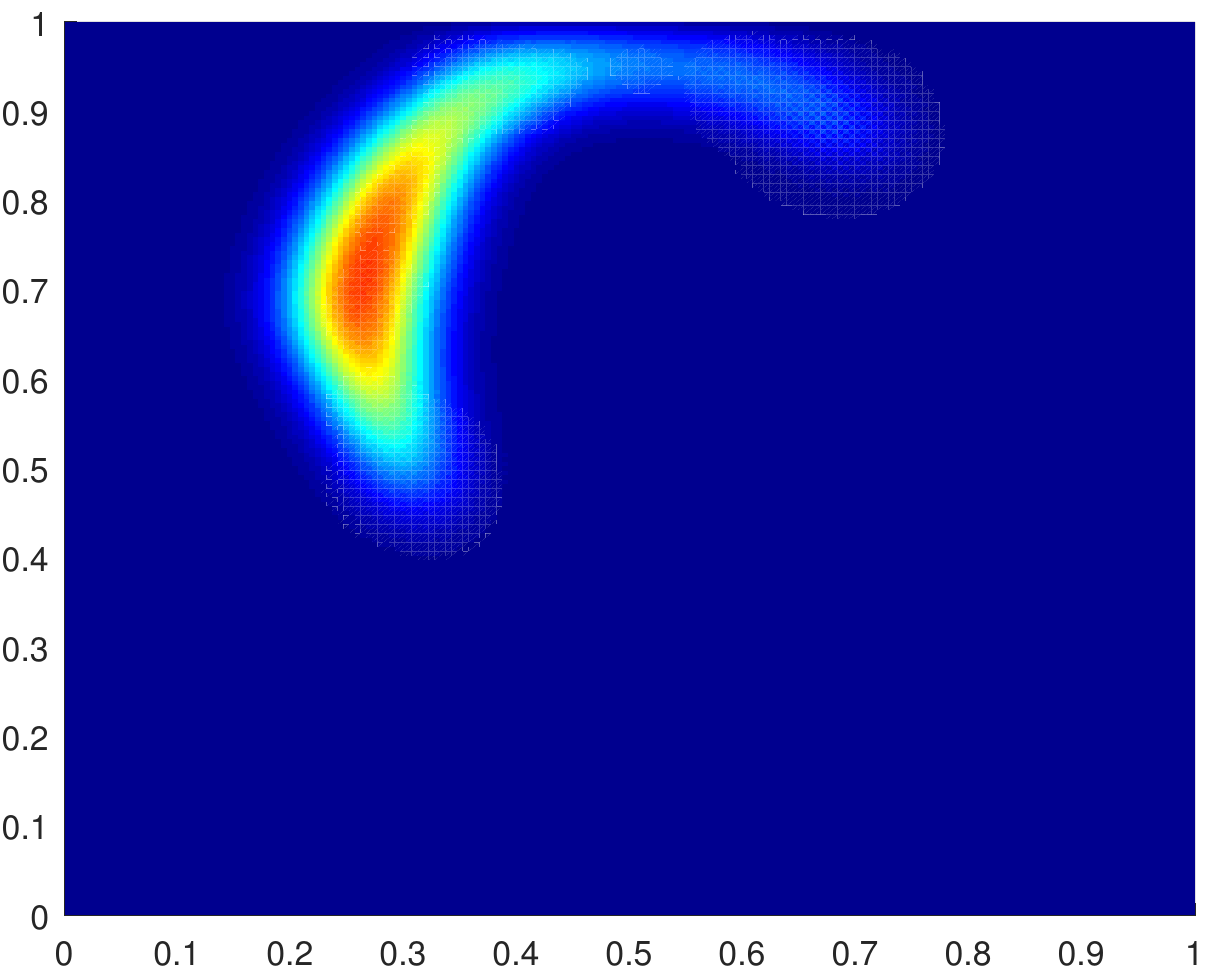}
	\includegraphics[width=0.3\textwidth]{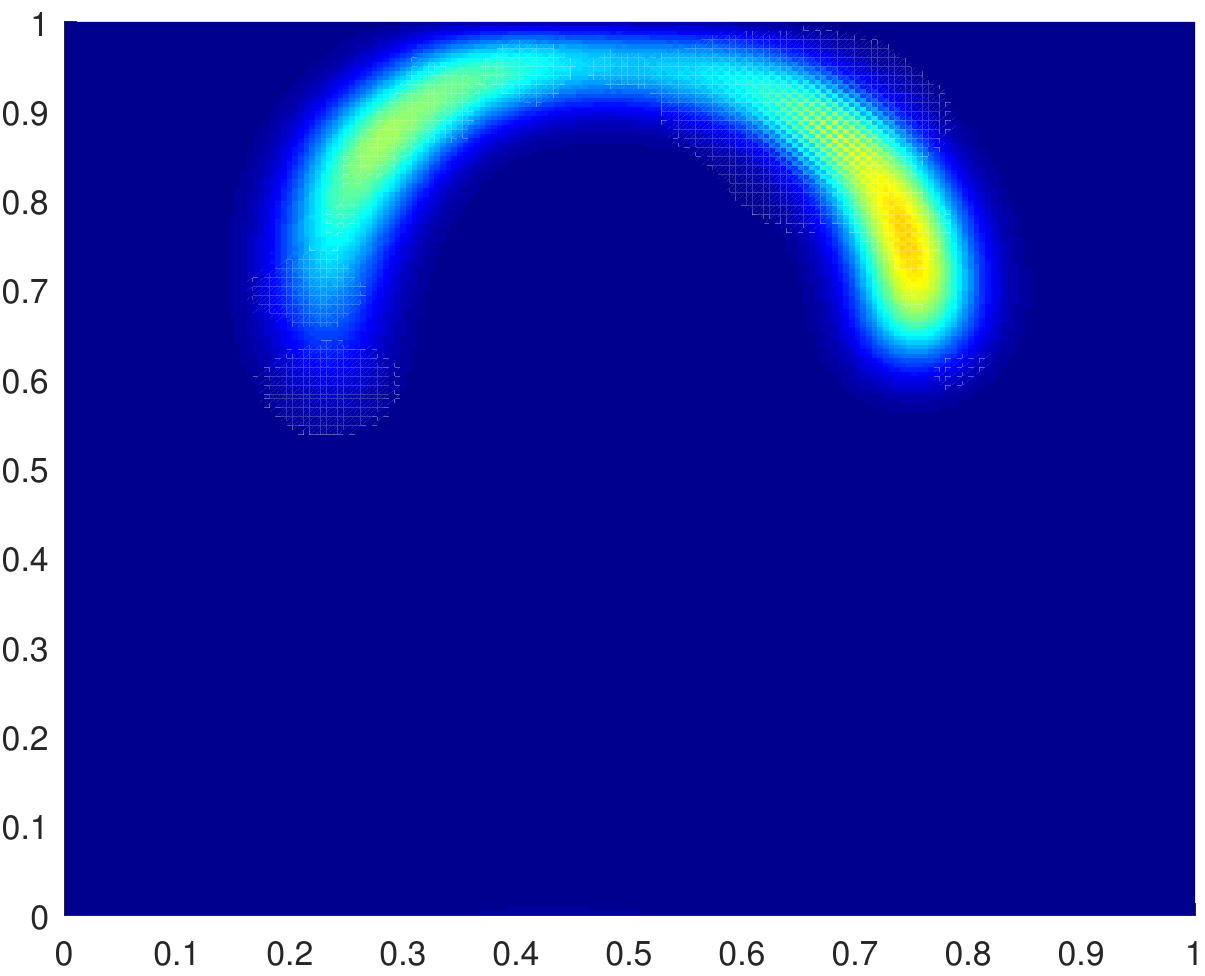}\\
	\includegraphics[width=0.3\textwidth]{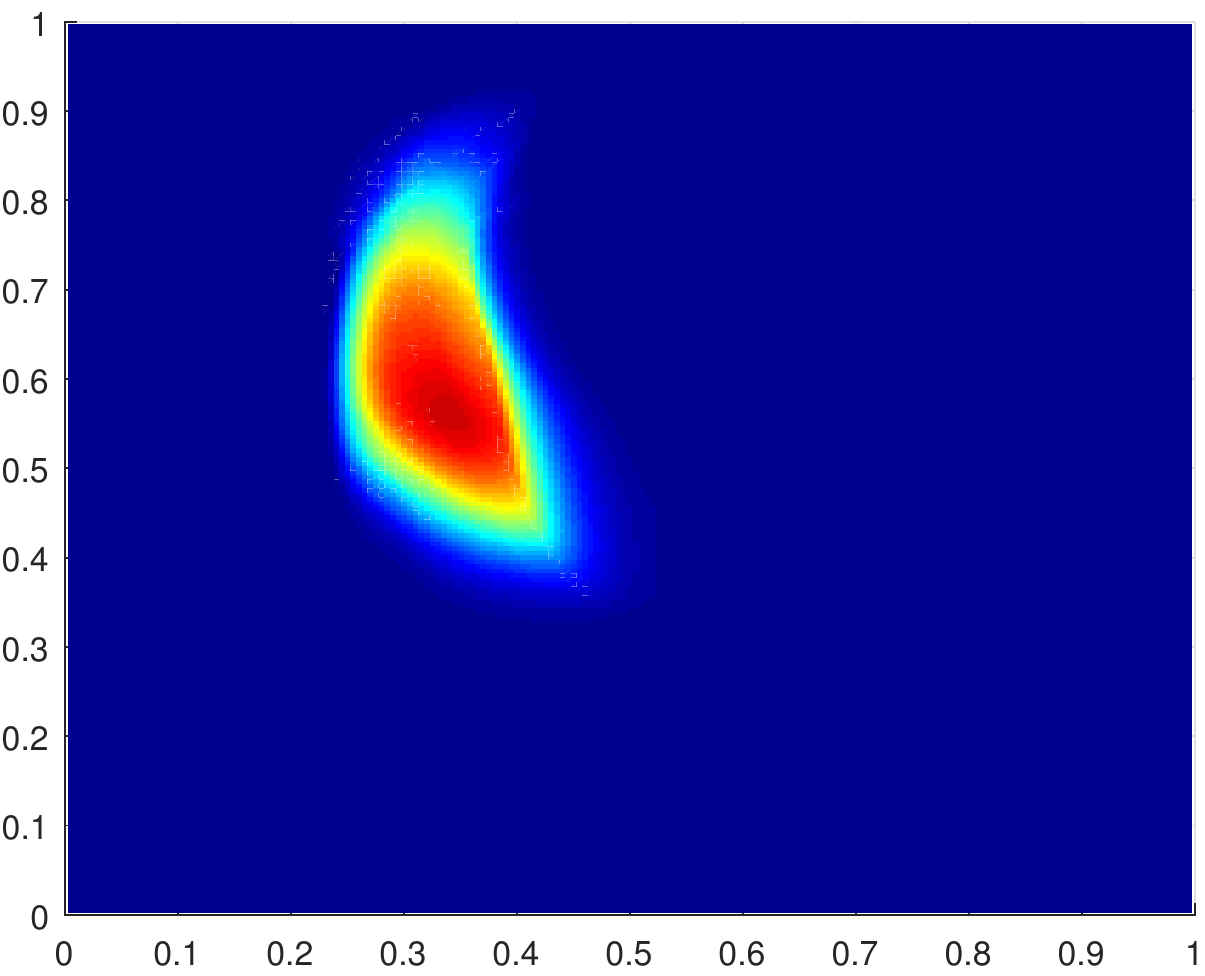}
	\includegraphics[width=0.3\textwidth]{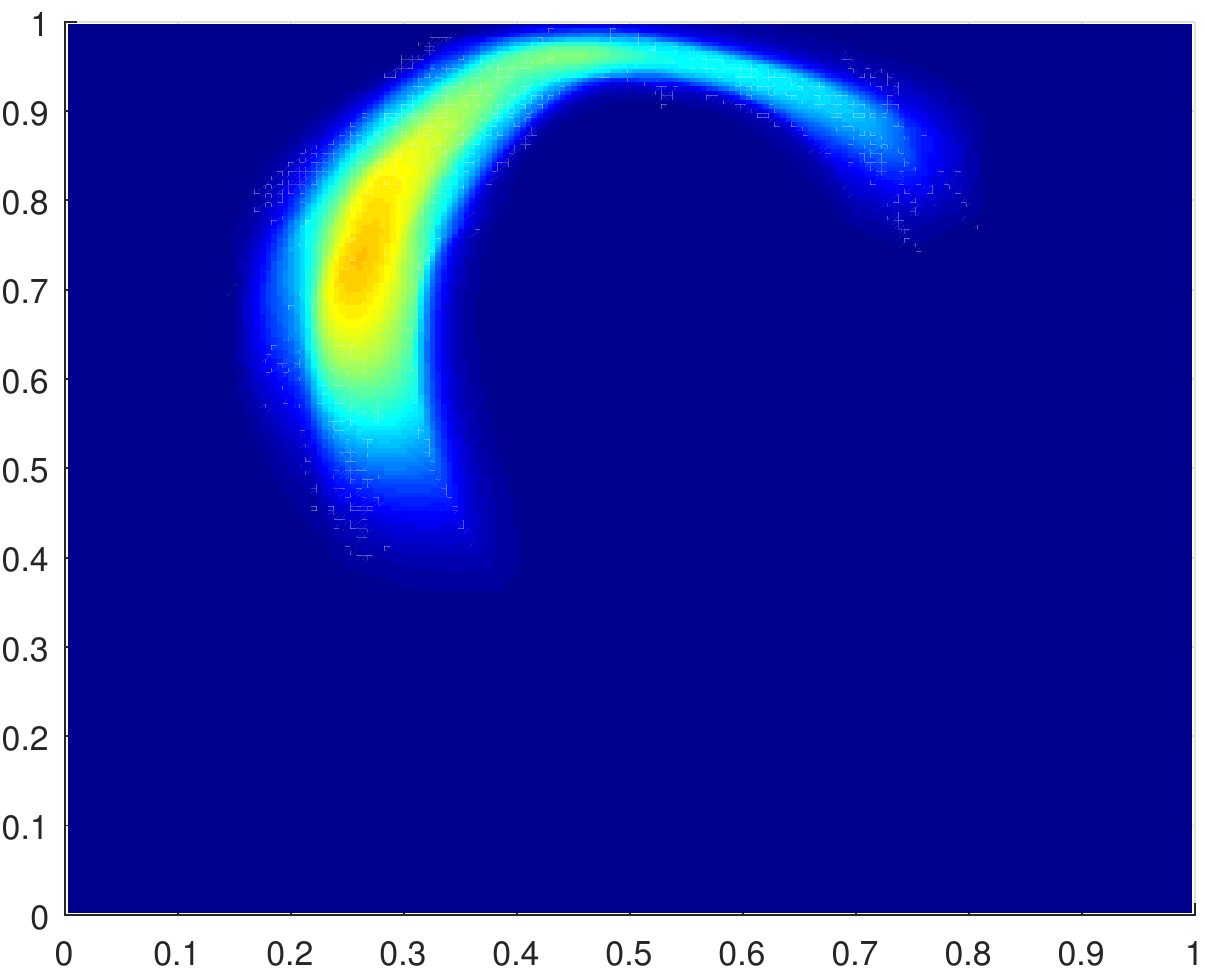}
	\includegraphics[width=0.3\textwidth]{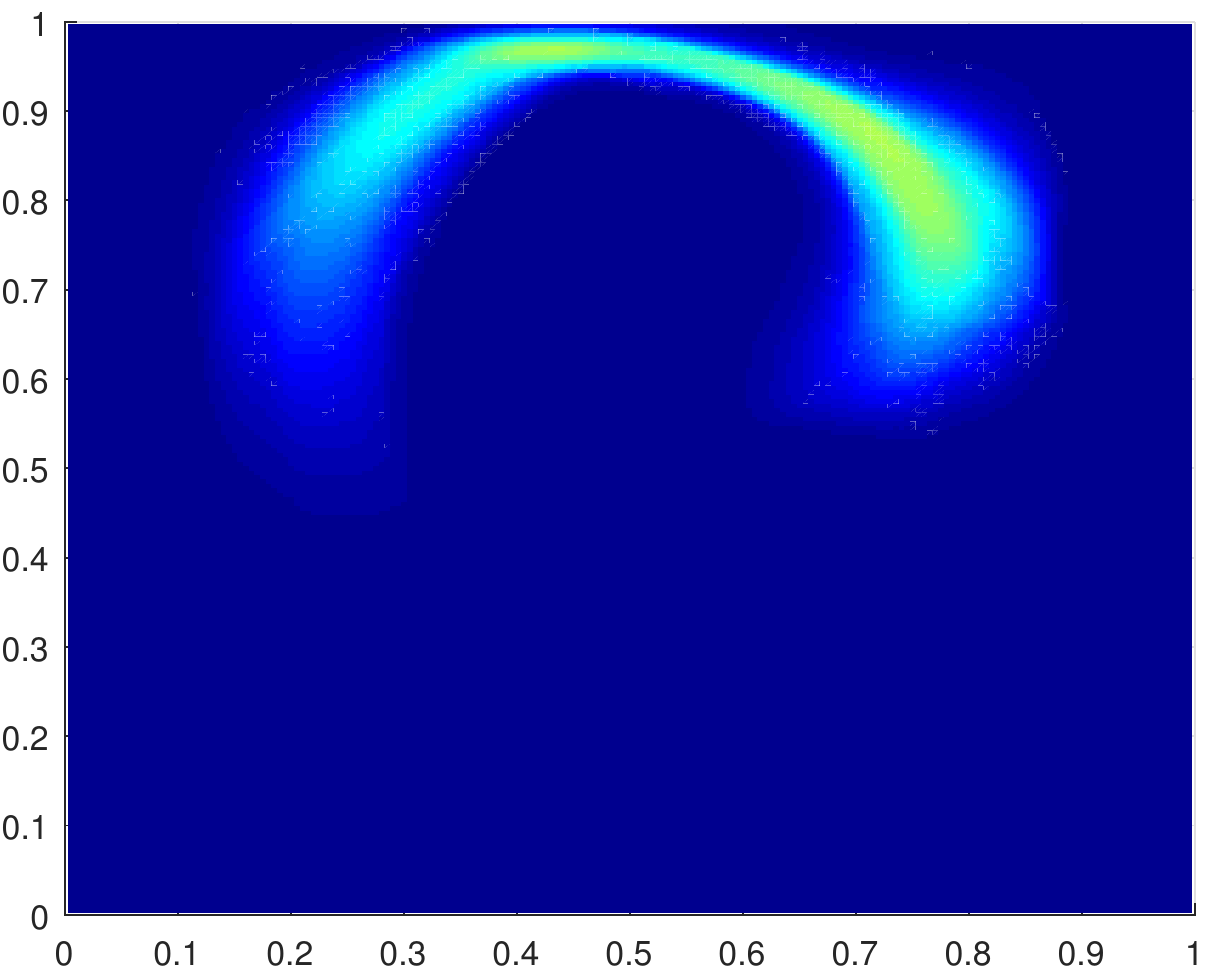}\\
	\centering\begin{tikzpicture}
	\begin{axis}[
	    hide axis,
	    scale only axis,
	    height=0pt,
	    width=0pt,
	    colormap/jet,
	    colorbar horizontal,
	    point meta min=0,
	    point meta max=25,
	    colorbar style={
	        width=10cm,
	        xtick={0,5,10,...,25}
	    }]
	    \addplot [draw=none] coordinates {(0,0)};
	\end{axis}
	\end{tikzpicture}
	\caption{Microscopic result is plotted on top, smoothed microscopic results  are plotted in the middle and the macroscopic result is plotted on the bottom. From left to right the results at time $T=2$, $3.5$, $5$ are shown.}
	\label{fig:cavitylam_ohne}
\end{figure}

To compare the microscopic results with the macroscopic ones, we compute $128$ Monte-Carlo realizations and the corresponding histograms.
For the model without interaction these are shown for different times in the top figure in Figure \ref{fig:cavitylam_ohne}.
In the middle, smoothed versions of the microscopic histograms are plotted. At the bottom, the solution of the macroscopic equation \eqref{eq:macro_rho_mono} without interaction is shown.

\begin{figure}[h!]
	\includegraphics[width=0.45\textwidth]{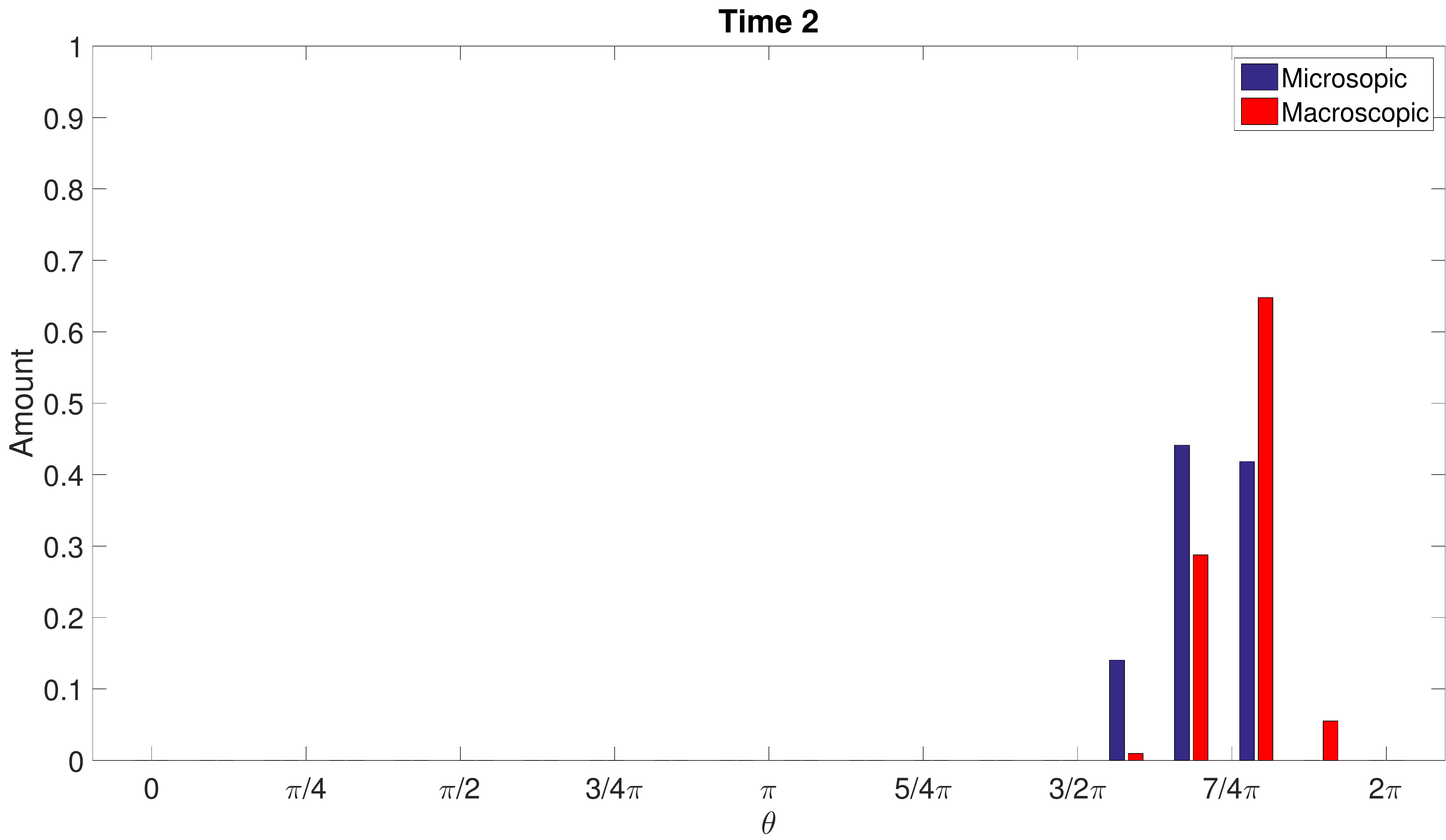}
	\includegraphics[width=0.45\textwidth]{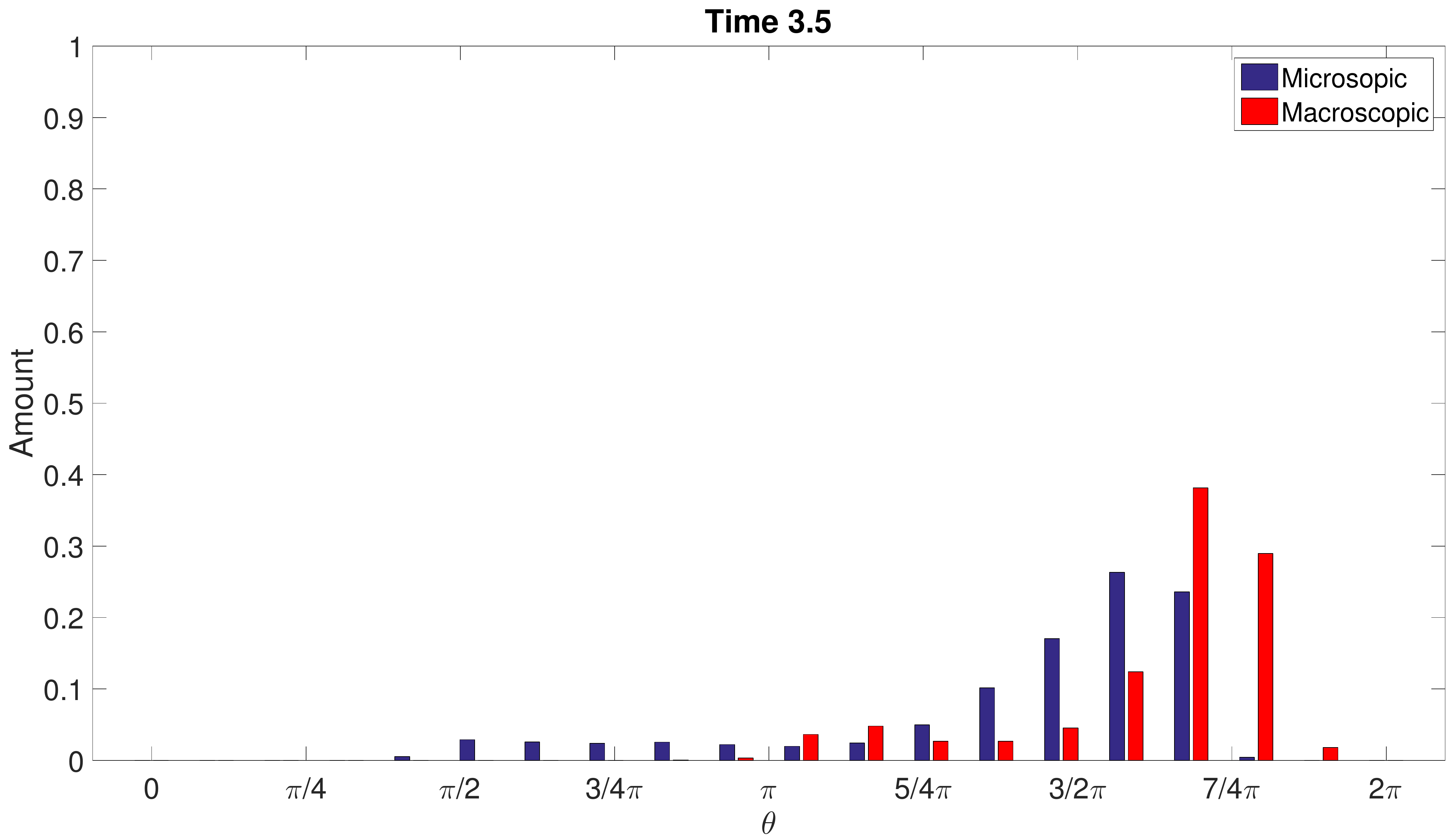}\\
	\center \includegraphics[width=0.45\textwidth]{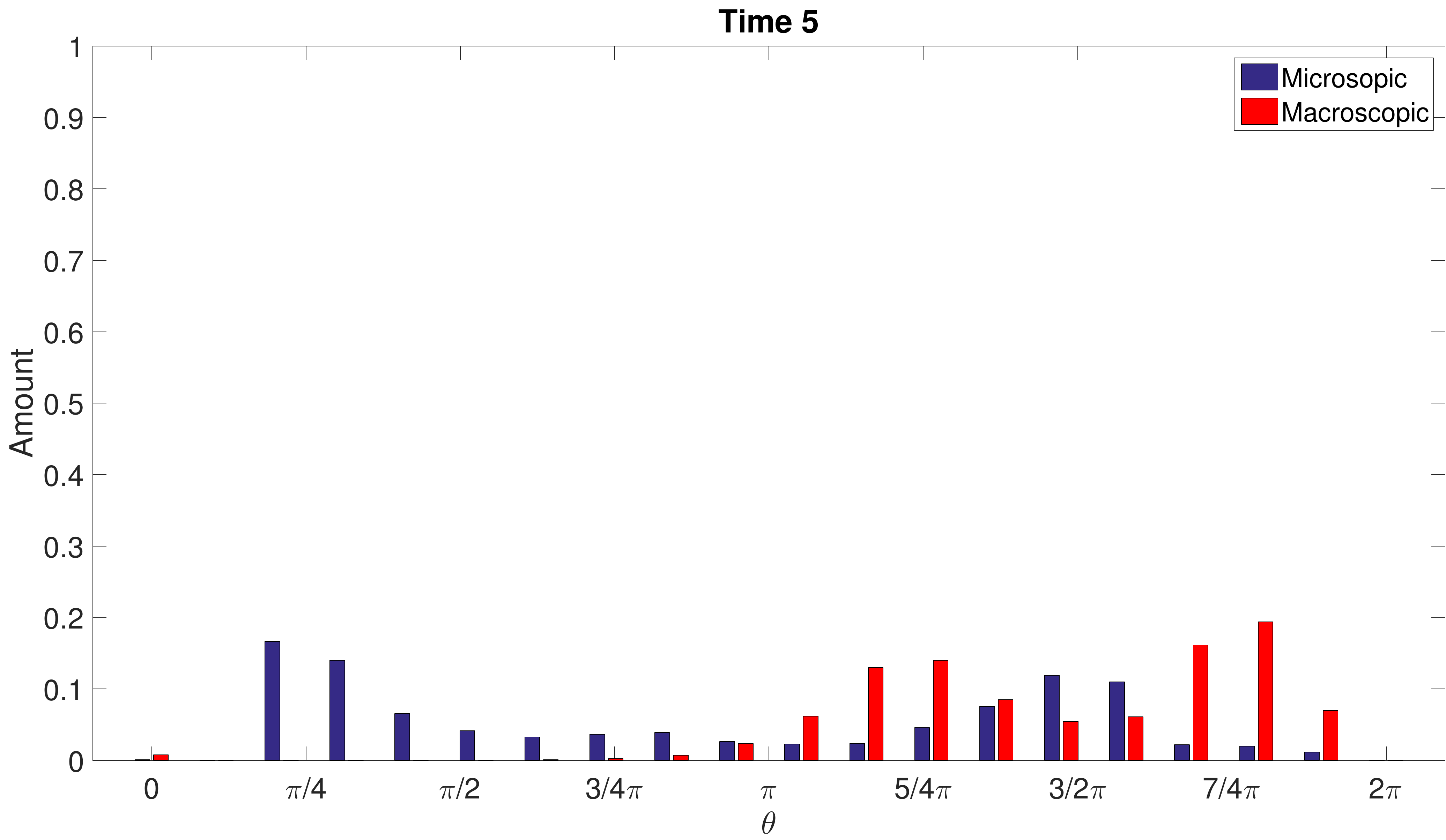}
	\caption{Angular distribution for the microscopic (blue) and the macroscopic results, at time T=2, 3.5, 5.}
	\label{fig:cavitylam_ohne_winkel}
\end{figure}

The initial group of particles follows closely the streamlines of the flow 
and is stretched where speeds increase.
At all time levels we can observe a good agreement of the microscopic and macroscopic simulation. 
This also holds true for the angular distribution, which is shown in Figure \ref{fig:cavitylam_ohne_winkel}.
In both models the ellipsoids are oriented  along  the streamlines of the surrounding fluid.
At $t=5$, slight differences in the angular distributions are observed.
Again, a numerical solution of the $q$ equation would yield a slighlty better approximation of the microscopic angular distribution.

\begin{figure}[h!]
	\includegraphics[width=0.3\textwidth]{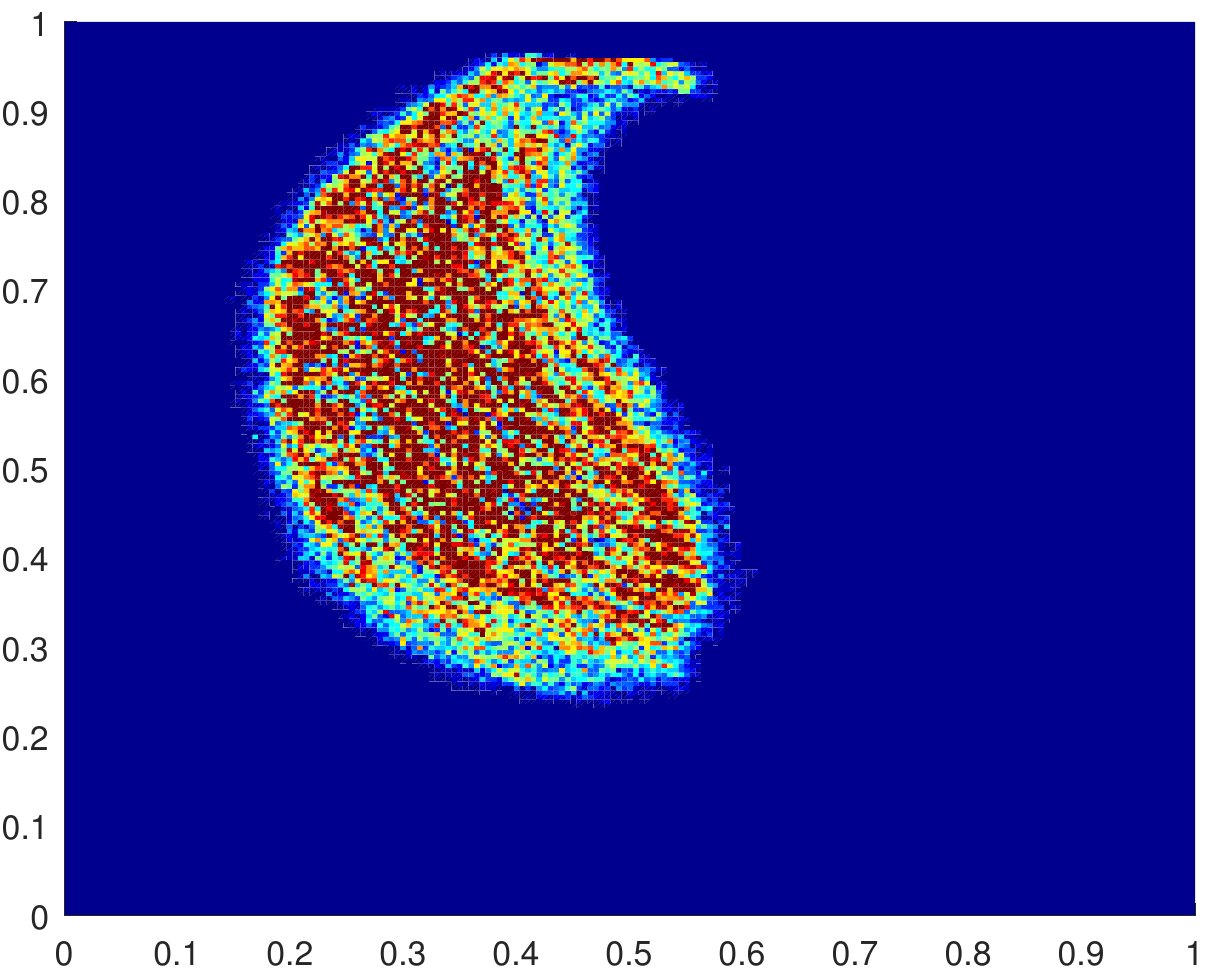}
	\includegraphics[width=0.3\textwidth]{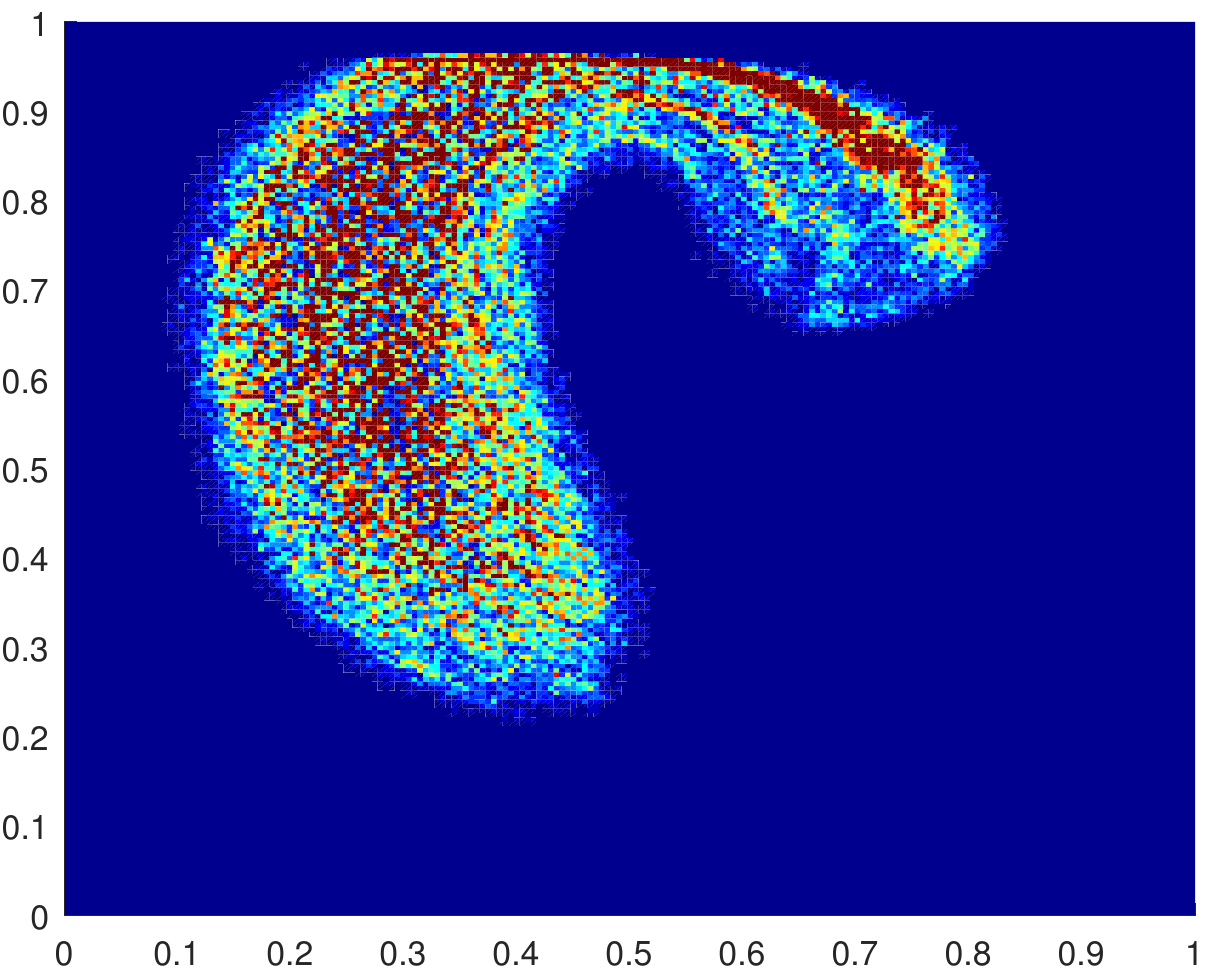}
	\includegraphics[width=0.3\textwidth]{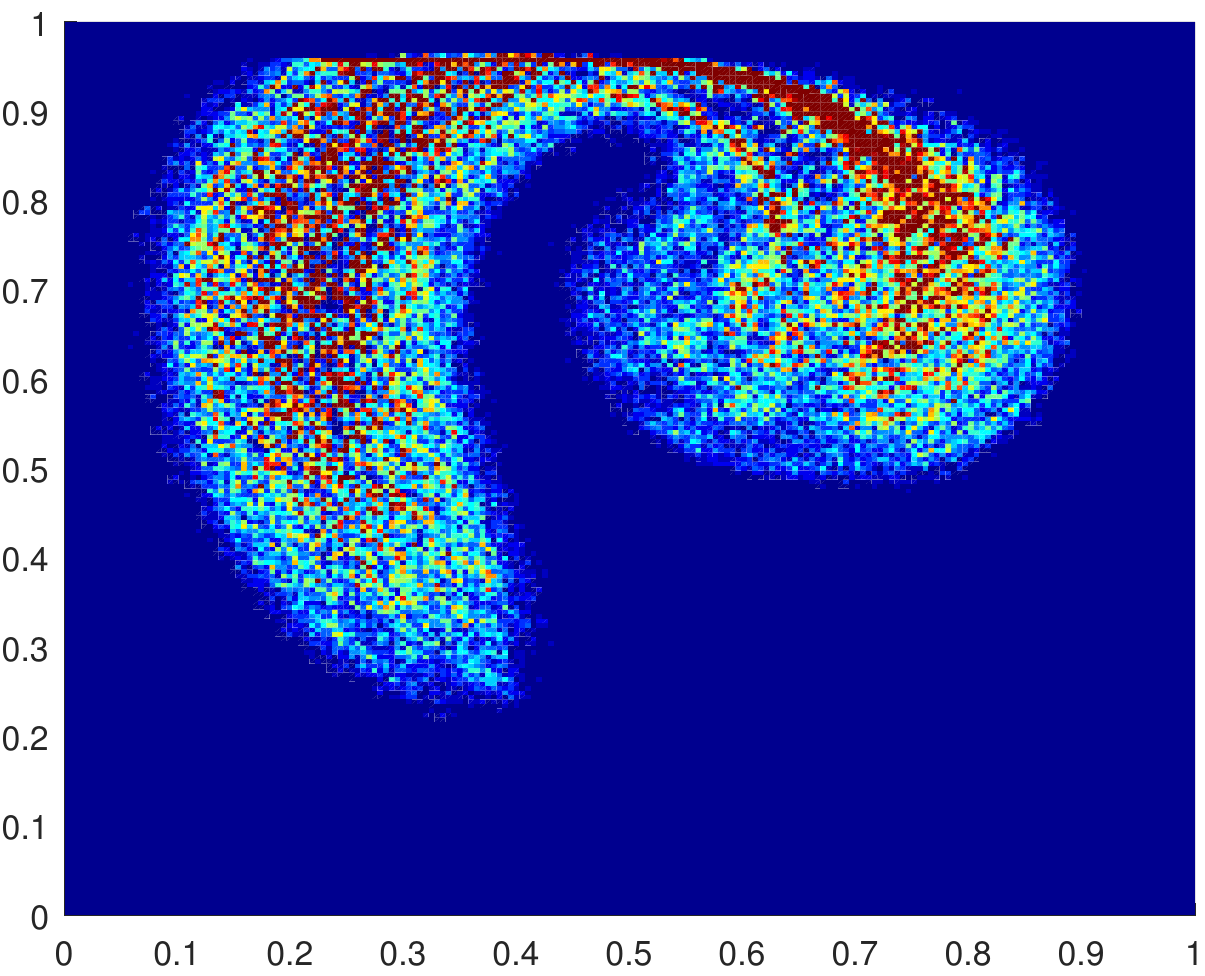}\\
	\includegraphics[width=0.3\textwidth]{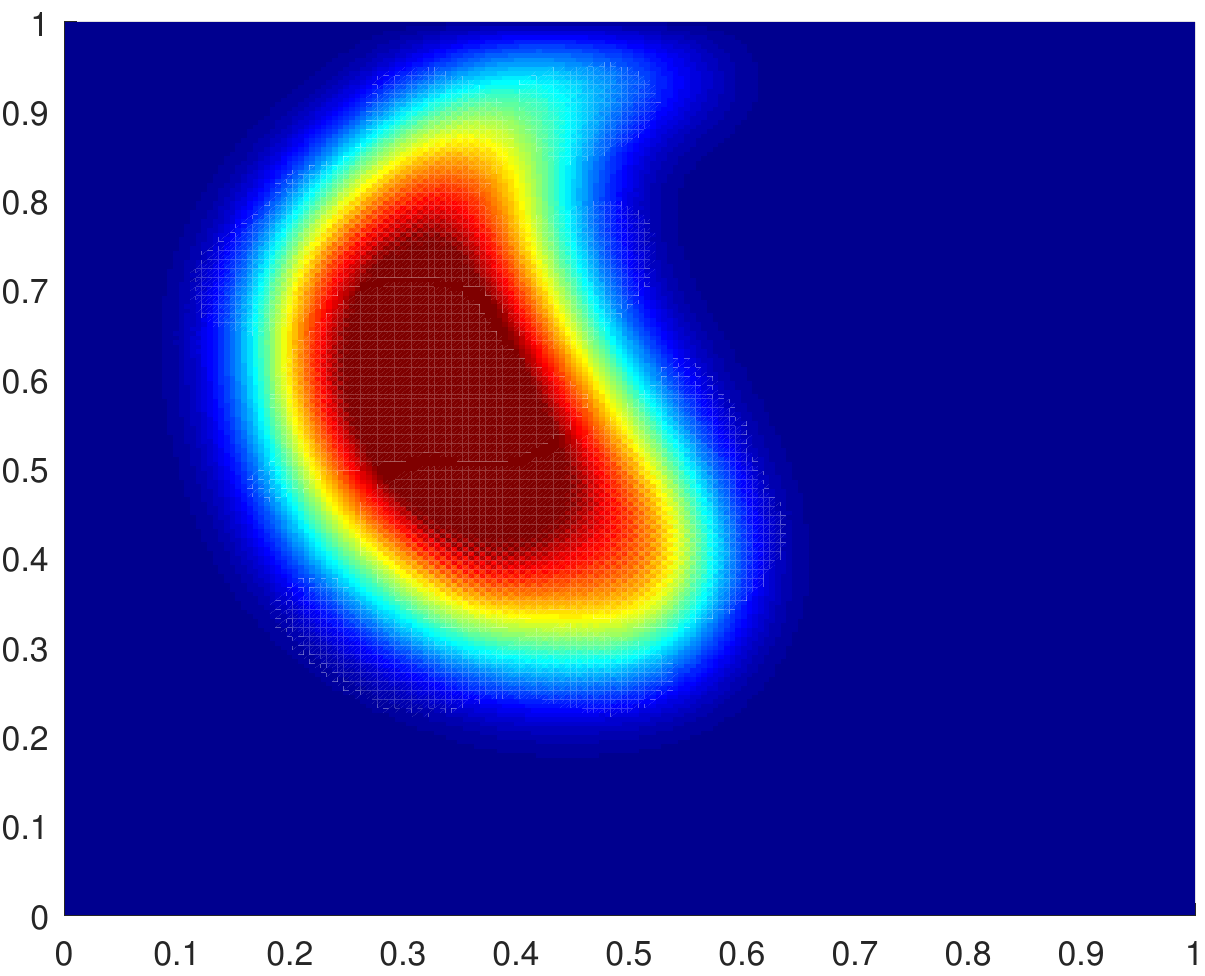}
	\includegraphics[width=0.3\textwidth]{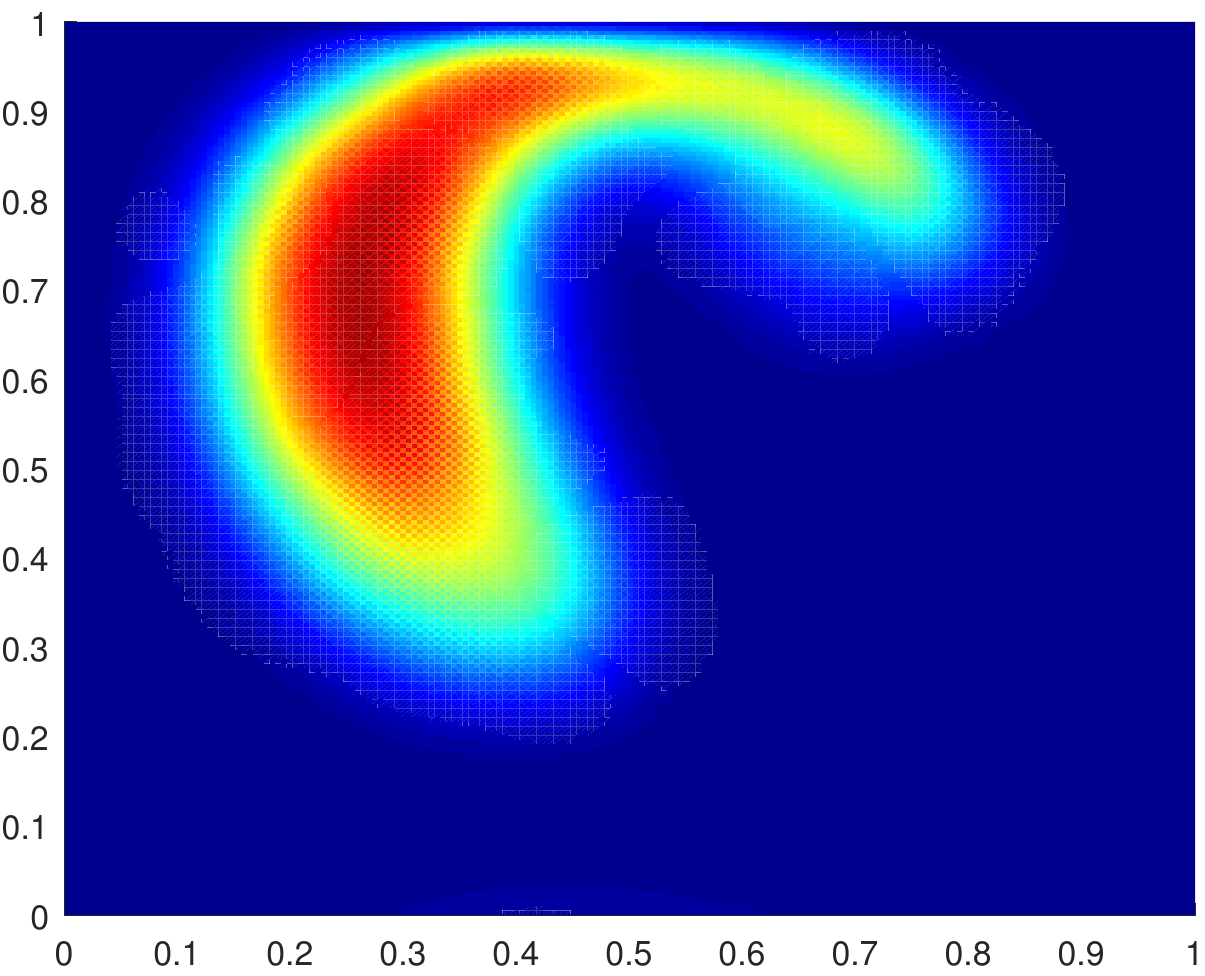}
	\includegraphics[width=0.3\textwidth]{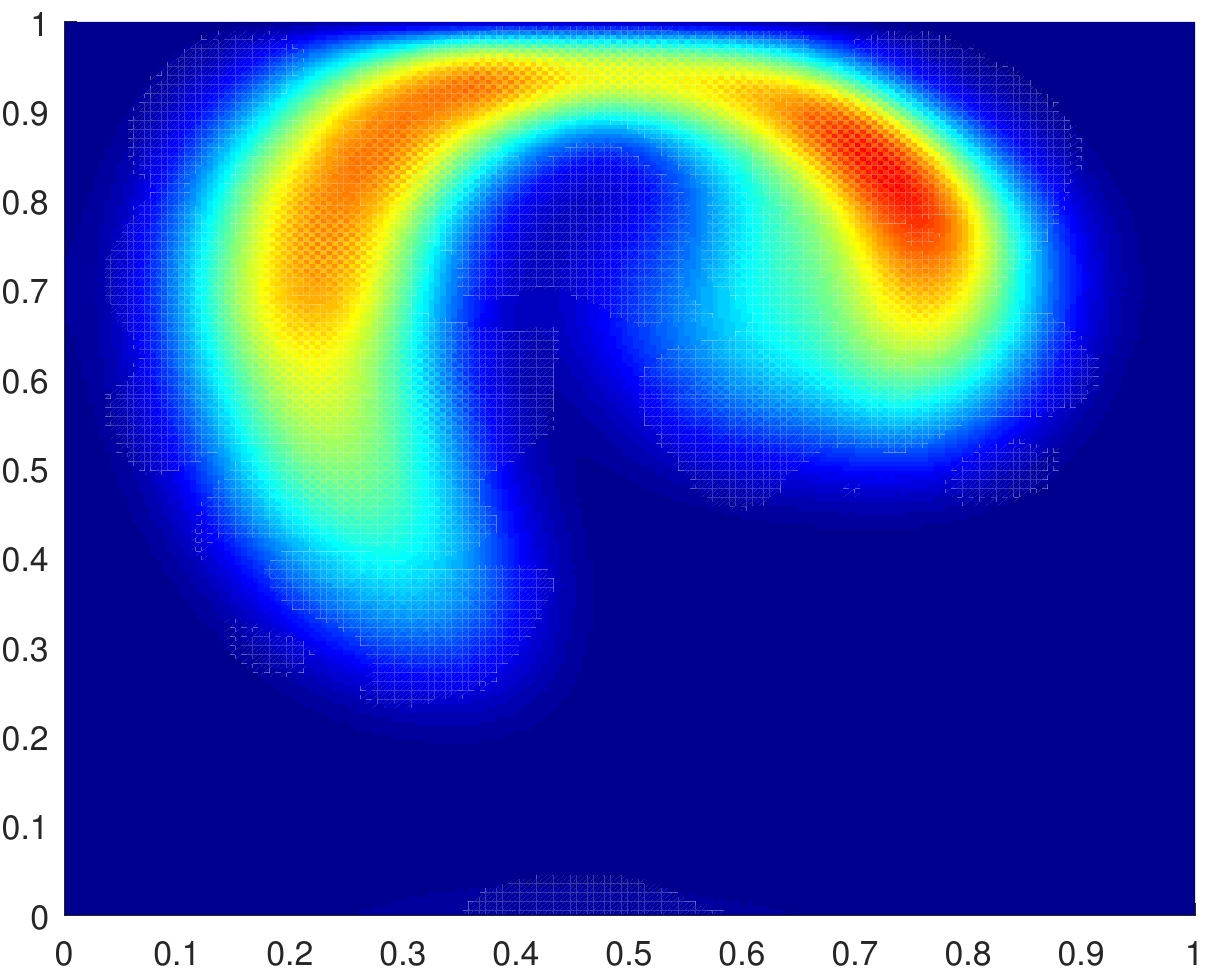}\\
	\includegraphics[width=0.3\textwidth]{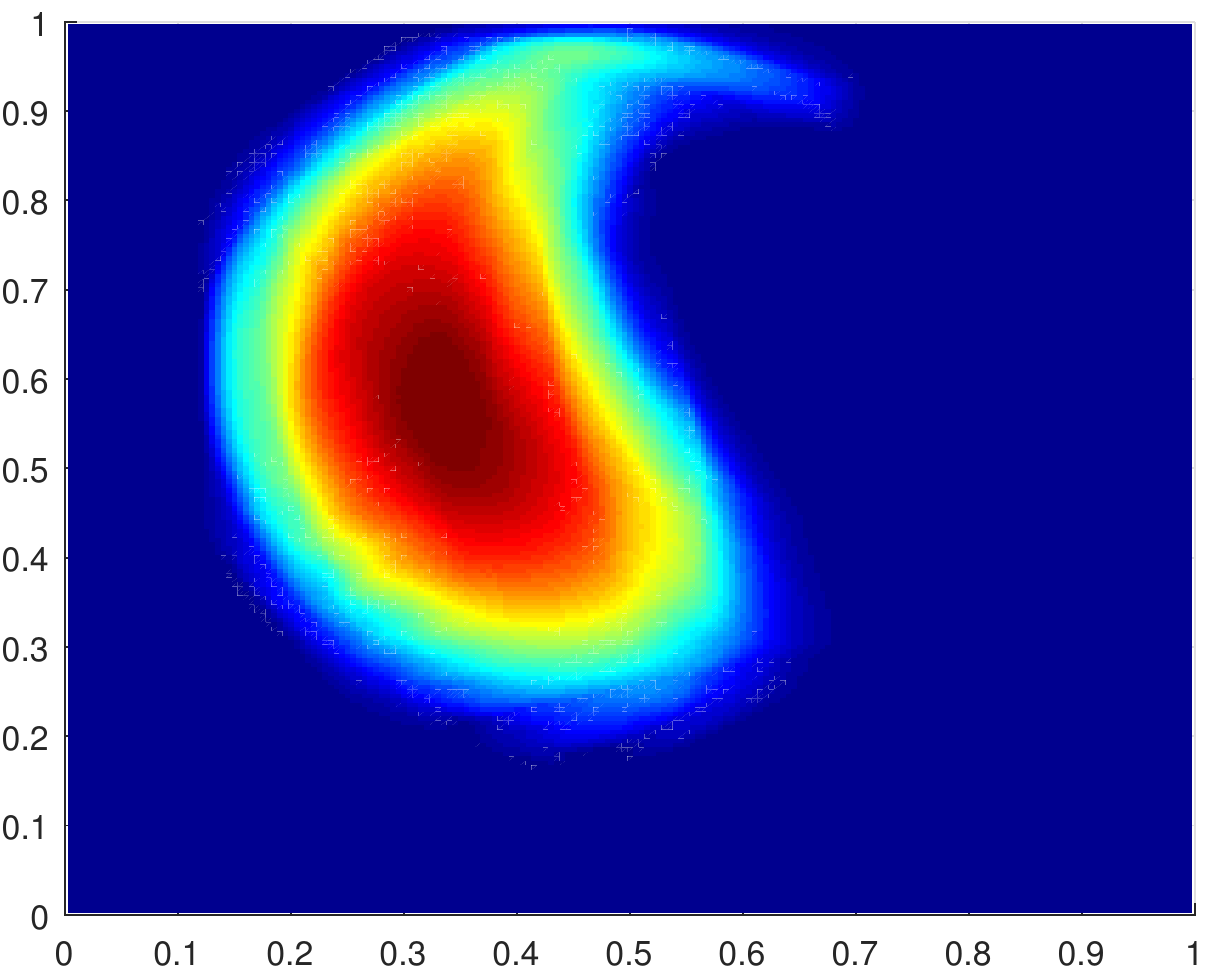}	
	\includegraphics[width=0.3\textwidth]{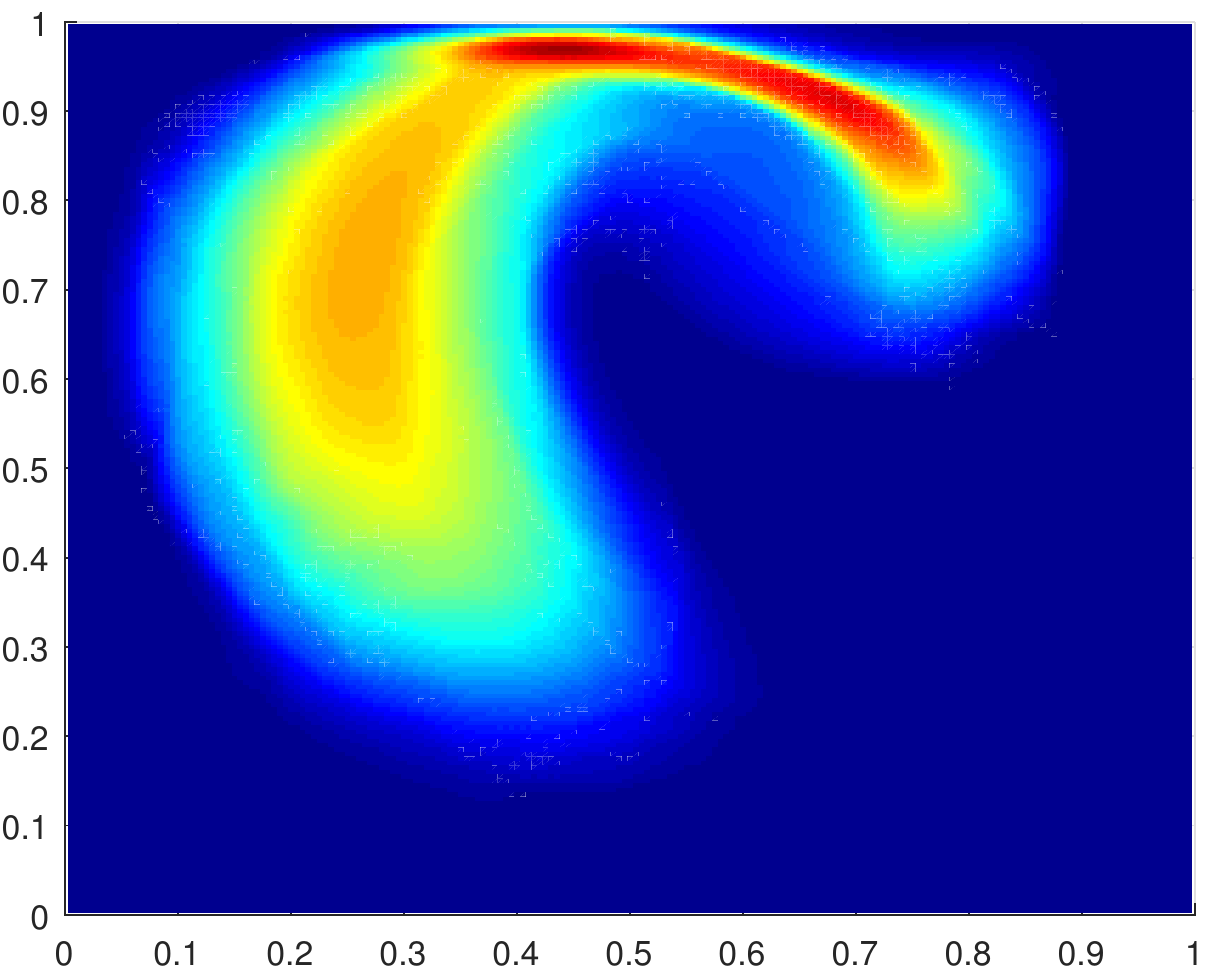}
	\includegraphics[width=0.3\textwidth]{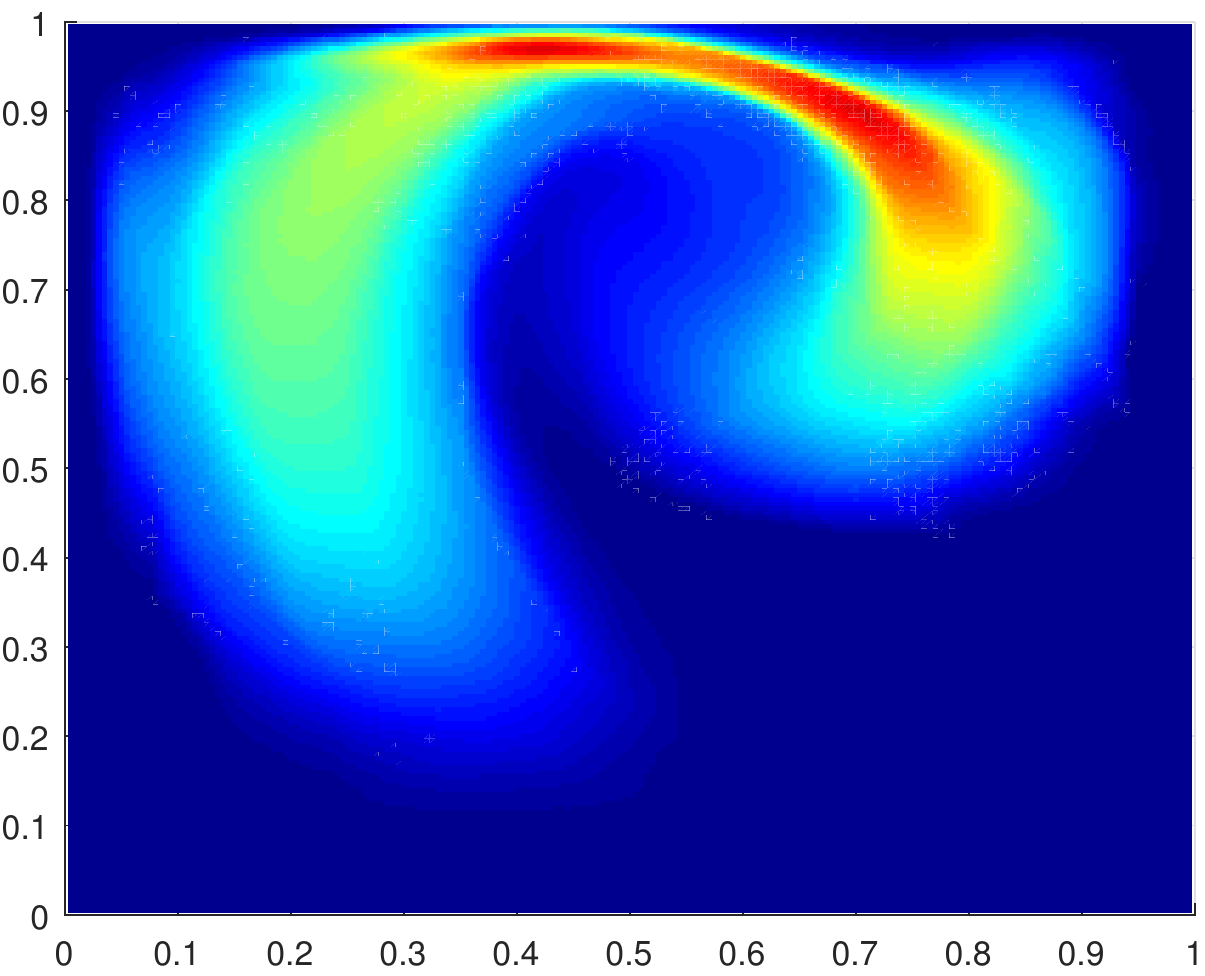}\\
	\centering\begin{tikzpicture}
		\begin{axis}[
		    hide axis,
		    scale only axis,
		    height=0pt,
		    width=0pt,
		    colormap/jet,
		    colorbar horizontal,
		    point meta min=0,
		    point meta max=6,
		    colorbar style={
		        width=10cm,
		        xtick={0,1,2,...,6}
		    }]
		    \addplot [draw=none] coordinates {(0,0)};
		\end{axis}
		\end{tikzpicture}
	\caption{Microscopic result is plotted on top, smoothed microscopic results are plotted in the middle and the macroscopic result is plotted on the bottom. From left to right the results at time $T=2$, $3.5$, $5$ are shown.}
	\label{fig:cavitylam_mit}
\end{figure}
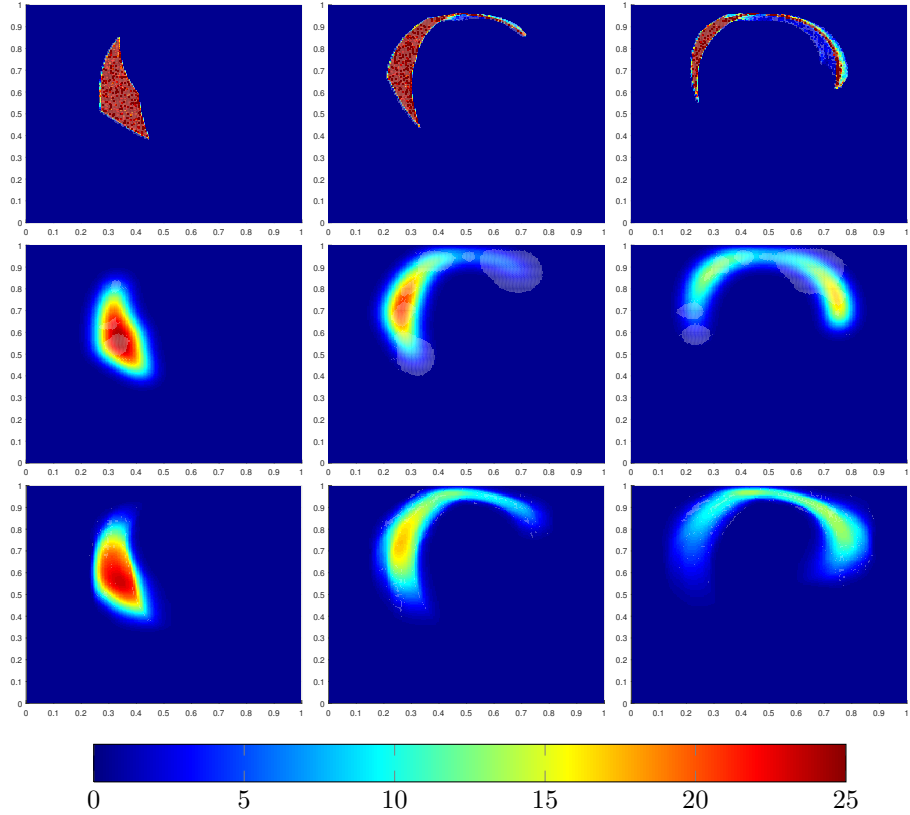

The identical test case with interaction is considered in Figure \ref{fig:cavitylam_mit}.

\begin{figure}[h!]
	\includegraphics[width=0.45\textwidth]{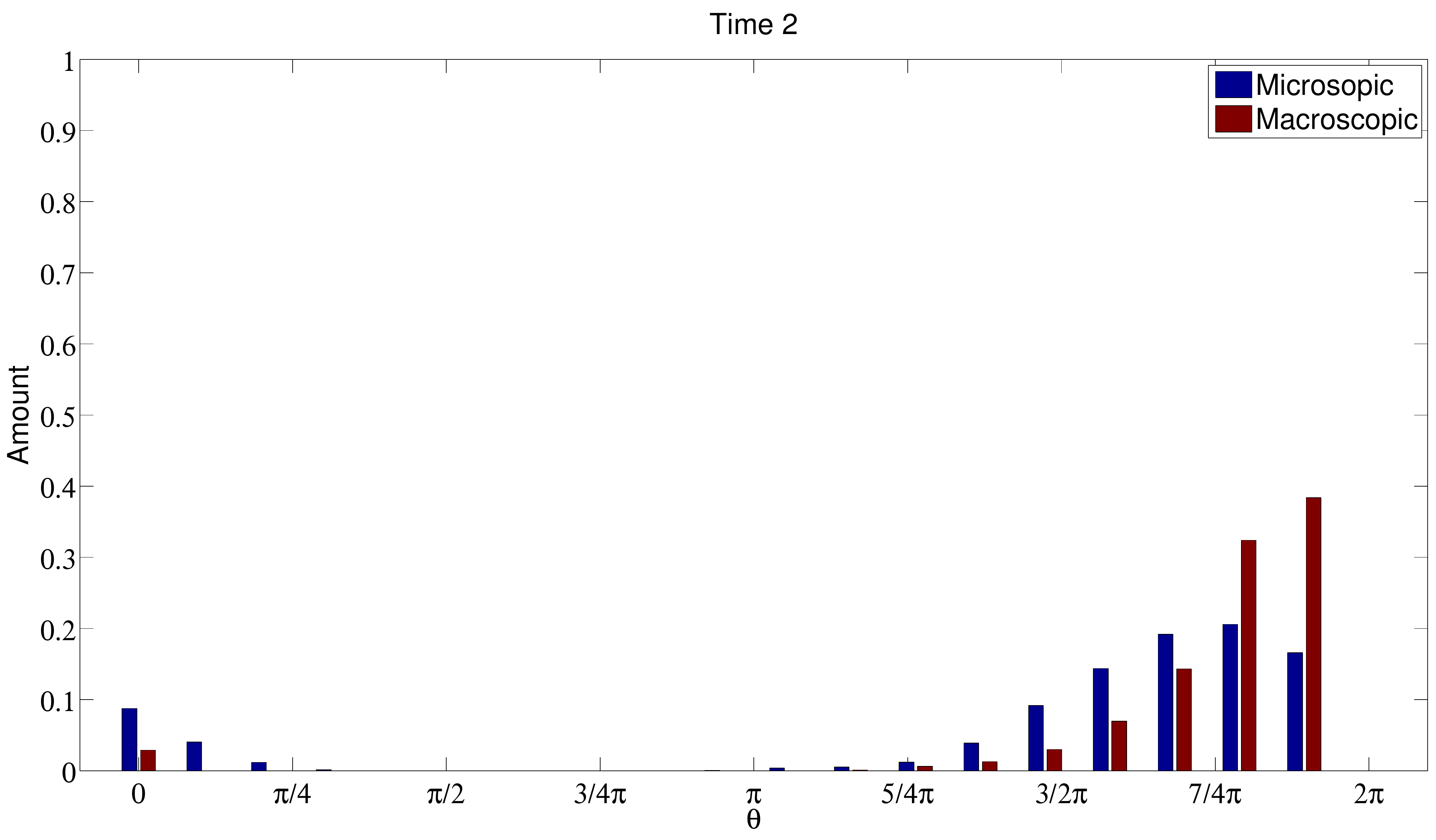}
	\includegraphics[width=0.45\textwidth]{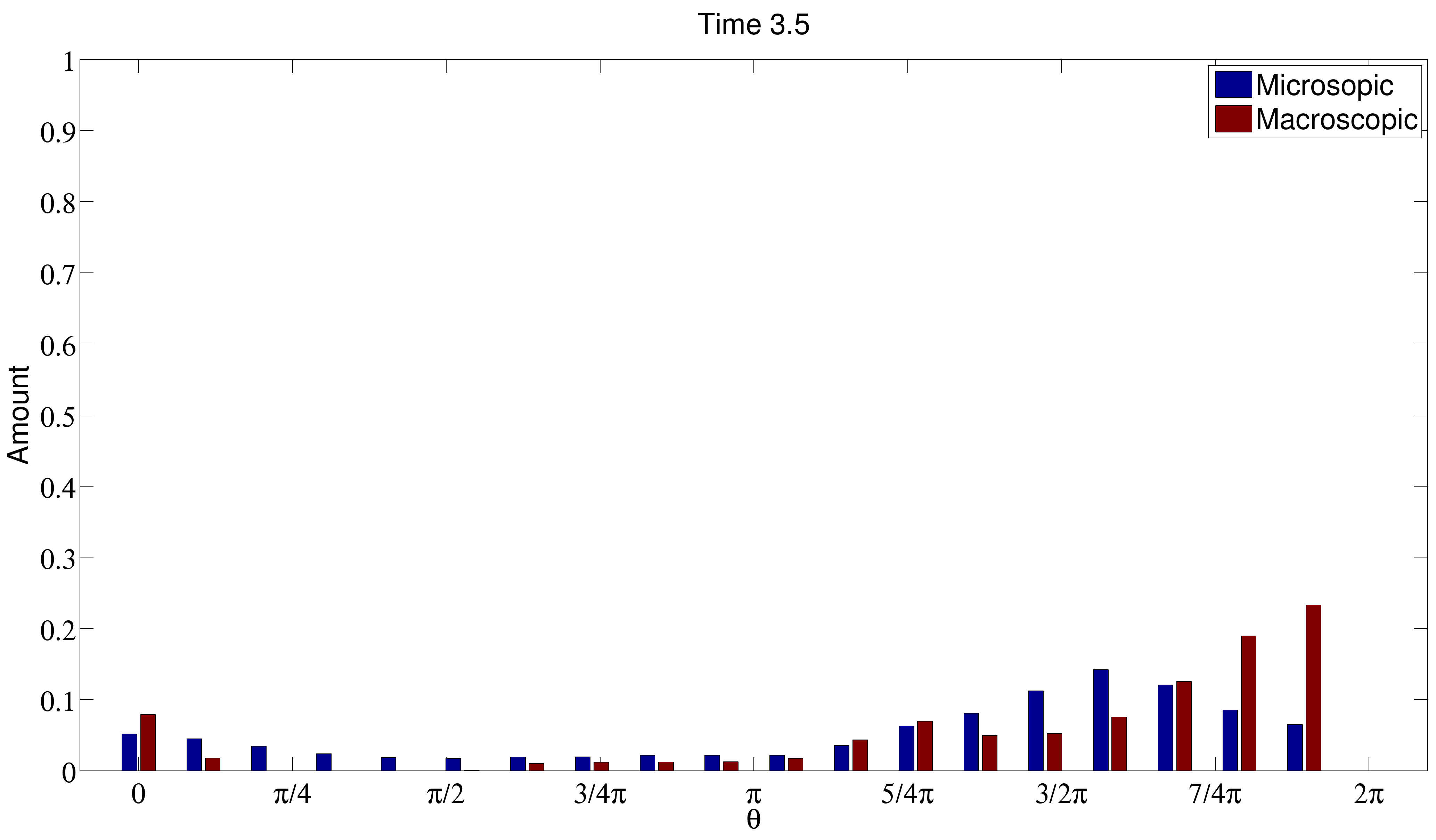}\\
	\center \includegraphics[width=0.45\textwidth]{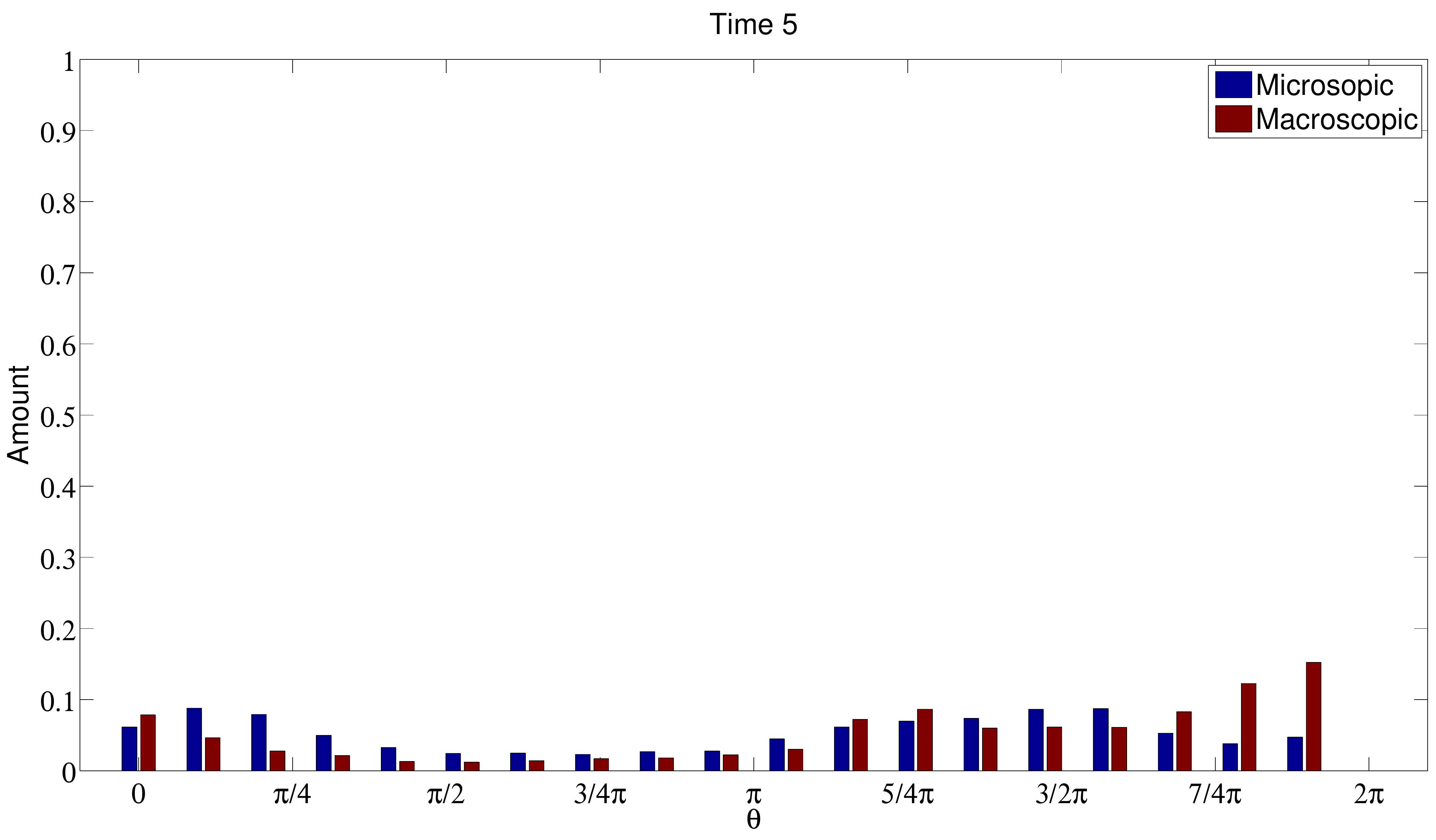}
	\caption{Angular distribution for the microscopic (blue) and the macroscopic results, at time $T=2$, $3.5$, $5$.}
	\label{fig:cavitylam_mit_winkel}
\end{figure}

Again the particles follow the surrounding fluid, but show a more diffusive behavior due to the interaction.
Again the different models show very similar results. Taking a closer look at the figures, one notices that the macroscopic solution approaches the upper boundary faster than the rest.
In the angular distributions in Figure \ref{fig:cavitylam_mit_winkel}, both models show similar results as well.

% % % % % % % % % % % % % % % %
% % %  Concluding Remarks % % %
% % % % % % % % % % % % % % % %
\section{Concluding remarks}
\label{sec:conclusion}
In this paper, we have investigated a spatially two-dimensional model  for a system of interacting ellipsoidal particles immersed in a fluid.
The particles are influenced by the  surrounding incompressible fluid and the mutual interaction modeled with a modified interaction potential of Berne \cite{berne1972gaussian}. We derived a kinetic mean field equation and closed the equations for different marginals of the distribution function using two simple closure procedures.
In the numerical tests, we investigated the quality of the mean field approximations. In particular, we compared a situation, where the distribution of the ellipsoids converges towards a stationary solution. Moreover, we investigated a top-bottom, a rotational and a driven cavity flow surrounding the ellipsoids.
In comparison with the microscopic solution, the hydrodynamic closures provide much better results than the diffusive approximation.
If one is interested in an accurate approximation of the angular distributions, closures for eqautions involving the 
$(r,\theta)$-marginals have to be considered.
In future work, more sophisticated closures might be investigated, e.g., moment closures that lead to an additional equation for the variance/temperature of the distribution function, which would allow for a wider range of applications.

\section*{Acknowledgment}
Funding by the Deutsche Forschungsgemeinschaft (DFG) within the RTG GrK 1932 "Stochastic Models for 
Innovations in the Engineering Sciences", project area P1, is gratefully acknowledged.

\bibliographystyle{plain}
\bibliography{literatur1,literatur2}

\end{document}